\documentclass{amsart}

\usepackage{amsmath,amsthm,amsfonts,amssymb,bm,graphicx, mathrsfs}

\usepackage
{hyperref}
\hypersetup{colorlinks=true,citecolor=blue,linkcolor=blue,urlcolor=blue,
pdfstartview=FitH }

\theoremstyle{plain}
\newtheorem{theorem}{Theorem}

\newtheorem{lemma}{Lemma}

\newtheorem{prop}[lemma]{Proposition}
\numberwithin{equation}{section}

\theoremstyle{definition}

\renewcommand{\geq}{\geqslant}
\renewcommand{\leq}{\leqslant}

\newcommand{\changed}[1]{{\color{black} #1}}

\usepackage{calc}
\newsavebox\CBox
\newcommand\hcancel[2][0.5pt]{%
  \changed{\ifmmode\sbox\CBox{$#2$}\else\sbox\CBox{#2}\fi%
  \makebox[0pt][l]{\usebox\CBox}%
  \rule[0.5\ht\CBox-#1/2]{\wd\CBox}{#1}}}

\newcommand{\E}{{\tt E}}

\makeatletter
\DeclareRobustCommand\widecheck[1]{{\mathpalette\@widecheck{#1}}}
\def\@widecheck#1#2{%
    \setbox\z@\hbox{\m@th$#1#2$}%
    \setbox\tw@\hbox{\m@th$#1%
       \widehat{%
          \vrule\@width\z@\@height\ht\z@
          \vrule\@height\z@\@width\wd\z@}$}%
    \dp\tw@-\ht\z@
    \@tempdima\ht\z@ \advance\@tempdima2\ht\tw@ \divide\@tempdima\thr@@
    \setbox\tw@\hbox{%
       \raise\@tempdima\hbox{\scalebox{1}[-1]{\lower\@tempdima\box
\tw@}}}%
    {\ooalign{\box\tw@ \cr \box\z@}}}
\makeatother

\begin{document}

\author{Valentin Blomer}

\address{Mathematisches Institut, Bunsenstr. 3-5, 37073 G\"ottingen, Germany} \email{vblomer@math.uni-goettingen.de}

 \author{Rizwanur Khan}
\address{Department of Mathematics, University of Mississippi, University, MS 38677, USA}
\email{rrkhan@olemiss.edu}

 \title{Twisted moments of $L$-functions and spectral reciprocity}

\thanks{First author partially supported by the DFG-SNF lead agency program grant BL 915/2-1 and BL 915/2-2.}
 
\keywords{Spectral reciprocity, moments of $L$-functions, amplification, subconvexity}

\begin{abstract}  A reciprocity formula is established that expresses  the fourth moment of automorphic  $L$-functions of level $q$ twisted by the $\ell$-th Hecke eigenvalue as the fourth moment of automorphic  $L$-functions of level $\ell$ twisted by the $q$-th Hecke eigenvalue.  Direct corollaries include subconvexity bounds for $L$-functions in the level aspect and a short proof of an upper bound for the fifth moment of automorphic $L$-functions. 
\end{abstract}

\subjclass[2010]{Primary: 11M41, 11F72}

\setcounter{tocdepth}{2}  \maketitle 

\maketitle

 \section{Introduction}

\subsection{Spectral reciprocity} Let $q, \ell$ be two  distinct odd primes. Gau{\ss}' celebrated law of quadratic reciprocity expresses the quadratic Legendre symbol $(\ell/q)$ in terms of $(q/\ell)$. This is a very remarkable statement, because it connects the arithmetic of two  unrelated finite fields $\Bbb{F}_{q}$ and $\Bbb{F}_{\ell}$. Hilbert re-interpreted quadratic reciprocity as a local-to-global statement: the local Hilbert symbols have to satisfy a global consistency relation (their product equals 1). 
The present paper establishes a   reciprocity formula between  the spectrum of the Laplacian on two different arithmetic hyperbolic surfaces $\Gamma_0(q) \backslash \Bbb{H}$ and $\Gamma_0(\ell)\backslash \Bbb{H}$ featuring a product of $L$-functions of total degree 8. One of the key steps involves the additive reciprocity formula \eqref{rec} below, which is used as a vehicle to shift information between different places of $\Bbb{Q}$, so that again a local-to-global principle works in the background. We will go into more detail later and proceed by describing our main result more precisely. 

Let $F$ be an automorphic form for the group ${\rm SL}_3(\Bbb{Z})$ and in the following summations  let generally denote by $f$ a cusp form for a congruence subgroup of ${\rm SL}_2(\Bbb{Z})$. Then roughly speaking, we will prove a formula of the shape
$$\sum_{f \text{ of level } q} L(s, f \times F) L(w, f) \lambda_f(\ell) \rightsquigarrow \sum_{f \text{ of level } \ell} L(s', f\times F) L(w', f) \lambda_f(q)$$
where
\begin{equation}\label{new}
\textstyle  s' = \frac{1}{2}(1 + w - s), \quad w' = \frac{1}{2}(3s + w - 1).
\end{equation}
 In particular, specializing $s = w = 1/2$ and $F = \E_0$ the minimal parabolic Eisenstein series with trivial spectral parameters whose Hecke eigenvalues are the ternary divisor function $\tau_3(n)$, we obtain a reciprocity formula for the twisted fourth moment 
\begin{equation}\label{twisted}
\sum_{f \text{ of level } q} L(1/2, f)^4 \lambda_f(\ell) \rightsquigarrow \sum_{f \text{ of level } \ell} L(1/2, f )^4 \lambda_f(q).
\end{equation}
In addition to its aesthetic features linking arithmetic data on different congruence quotients of the upper half plane, this formula is useful for applications, because we can trade the level of the forms for the twisting Hecke eigenvalue. We will demonstrate this is in the next subsection. We also mention that such a formula is implicit in the work of K\i ral and Young \cite{KY}, although in a coarser and less conceptual form. After having stated the main result in Theorem \ref{thm1} below, we will also discuss other versions of reciprocity of central $L$-values. \\ 

For the rest of the paper let $q, \ell$ be  two   positive integers,  $s, w$ two complex numbers with positive real part, and  $F$  a cuspidal or non-cuspidal automorphic form for the group ${\rm SL}_3(\Bbb{Z})$ with Fourier coefficients $A(n_1, n_2)$. 
  We need to   (a)  parametrize the spectrum of level $q$ and identify suitable test functions on the various components, (b)   define correction factors at the ramified primes $p \mid q\ell$, and (c) add additional ``main terms''. 

Irreducible automorphic representations of $Z(\Bbb{A}_{\Bbb{Q}}) {\rm GL}_2(\Bbb{Q})\backslash {\rm GL}_2(\Bbb{A}_{\Bbb{Q}})/K_0(q)$ with $K_0(q) = \{\left(\begin{smallmatrix} a& b\\ c& d \end{smallmatrix}\right)\in {\rm GL}_2(\widehat{\Bbb{Z}}) \mid   c\equiv 0 \, (\text{mod } q \widehat{\Bbb{Z}})\}$ can be generated by three types of functions:
 \begin{itemize}
 \item Cuspidal holomorphic newforms $f$  of weight $k = k_f  \in 2\Bbb{N}$, level $q' \mid q$ and Hecke eigenvalues $\lambda_f(n) \in \Bbb{R}$; we denote the   set of such forms by $\mathcal{B}_{\text{hol}}(q)$ and the set of newforms of exact level $q$ by $\mathcal{B}^{\ast}_{\text{hol}}(q)$; 
 \item Cuspidal Maa{\ss} newforms $f$ of spectral parameter $t = t_f\in \Bbb{R} \cup [-i \vartheta, i \vartheta]$,  level $q' \mid q$ and Hecke eigenvalues $\lambda_f(n) \in \Bbb{R}$, where at the current state of knowledge  $\vartheta = 7/64$ can be taken \cite{KiS} (though $\vartheta = 0$ is expected); we denote the  set of such forms by $\mathcal{B}(q)$ and the set of newforms of exact level $q$ by $\mathcal{B}^{\ast}(q)$; 
 \item Unitary Eisenstein series $E_{t, \chi}$ for $\Gamma_0(q)$, where $t \in \Bbb{R}$ and $\chi$ is a primitive Dirichlet character of conductor $c_{\chi}$ satisfying $c_{\chi}^2 \mid q$. Their $n$-th Hecke eigenvalue  is $\lambda_{t, \chi}(n) = \sum_{ab = n} \chi(a)\bar{\chi}(b) (b/a)^{it}$ for $(n, q) = 1$.
\end{itemize}

We write $\mathcal{T}_0 := (\Bbb{R} \cup [-i \vartheta, i \vartheta] ) \times 2\Bbb{N}$ 
and consider the set of functions $\mathfrak{h} = (h, h^{\text{hol}}) : \mathcal{T}_0 \rightarrow \Bbb{C}$, where
\begin{equation}\label{weakly}
  h(-t) = h(t) \ll (1 + |t|)^{-15}, \,\, t\in \Bbb{R} \cup [-i \vartheta, i \vartheta], \quad \quad h^{\text{hol}}(k) \ll k^{-15}, \,\,k \in 2\Bbb{N}. 
\end{equation}
 We call such functions \emph{weakly admissible}.    We will define in Section \ref{admissible} the notion of an \emph{admissible} function $\mathfrak{h}$ which imposes    additional analytic conditions on $\mathfrak{h}$;  see Lemma \ref{lem1}  for explicit examples. 
 
For a cuspidal (holomorphic or Maa{\ss}) newform $f$ or an Eisenstein series $E_{t, \chi}$ and $w \in \Bbb{C}$ define
\begin{equation}\label{Lambda}
\Lambda_f(\ell; w) := \sum_{ab = \ell} \frac{\mu(a)}{a^w} \lambda_f(b), \quad \Lambda_{(t, \chi)}(\ell; w) := \sum_{ab = \ell} \frac{\mu(a)}{a^w} \lambda_{t, \chi}(b).
\end{equation}
 We also need  local correction factors $\tilde{L}_q(s, w, f \times F)$ resp.\ $\tilde{L}_q(s, w, E_{t, \chi} \times F)$ that are defined in \eqref{local-q-cusp} and \eqref{local-q-eis} and satisfy the following properties.
 \begin{lemma}\label{lem2}
For $f \in \mathcal{B}^{\ast}(q) \cup \mathcal{B}^{\ast}_{\text{\rm hol}}(q)$ we have 
$ \tilde{L}_q(s, w, f \times F) =1$. 
 For arbitrary $f \in \mathcal{B}(q) \cup \mathcal{B}_{\text{\rm hol}}(q)$, $t \in \Bbb{R}$, $\chi$ a primitive Dirichlet character of conductor $c_{\chi}$ with $c_{\chi}^2 \mid q$, $\varepsilon > 0$ and $\Re s, \Re w \geq 1/2$ we have the uniform bound
\begin{equation}\label{local-bound}
 \tilde{L}_q(s, w, f \times F), \,\, \tilde{L}_q(s, w, E_{t, \chi} \times F)  \ll_{\varepsilon} q^{\theta +\varepsilon}.
 \end{equation}
  \end{lemma} 
   Here and henceforth, $\theta$ (not to be confused with $\vartheta$) denotes an admissible exponent towards the Ramanujan conjecture for $F$. At the current state of knowledge, $\theta \leq 5/14$ (cf.\ \cite{BB}) is known in general, and for $F = \E_0$ we have $\theta = 0$.  \\ 

For a weakly admissible function $\mathfrak{h} = (h,  h^{\text{hol}})$, complex numbers $s, w$ with $\Re s, \Re w \geq  1/2$ and coprime integers $q, \ell$  we define the twisted spectral mean values
\begin{equation}\label{moment}
\begin{split}
&  \mathcal{M}^{\text{Maa{\ss}}, \pm}_{q, \ell}(s, w; \mathfrak{h}) :=   \frac{\phi(q)}{q^2} \sum_{d_0 \mid q} \sum_{f \in \mathcal{B}^{\ast}(d_0)} \epsilon_f^{(1 \mp 1)/2}\frac{L(s,  f \times F) L(w, f)}{L(1, \text{Ad}^2 f)}\tilde{L}_q(s, w,   f \times F)\frac{   \Lambda_f(\ell; w)}{\ell^{w}}     h(t_f),\\
&  \mathcal{M}^{\text{hol}}_{q, \ell}(s, w; \mathfrak{h}) :=    \frac{\phi(q)}{q^2}   \sum_{d_0 \mid q}   \sum_{f \in \mathcal{B}^{\ast}_{\text{hol}}(d_0)} \frac{L(s,    f \times F) L(w, f)}{L(1, \text{Ad}^2 f)}\tilde{L}_q(s, w,  f \times F)  \frac{ \Lambda_f(\ell; w)}{\ell^{w}}   h^{\text{hol}}(k_f),\\
&  \mathcal{M}^{\text{Eis}}_{q, \ell}(s, w; \mathfrak{h}) :=   \frac{\phi(q)}{q^2}\sum_{\chi : c_{\chi}^2 \mid q}  \int_{\Bbb{R}}   \frac{L(s + it, F \times \chi) L(s - it, F \times \bar{\chi}) L(w + it, \chi)L(w - it, \bar{\chi})}{L(1 + 2 it, \chi^2)L(1 - 2it, \bar{\chi}^2)} \\
& \quad\quad\quad\quad\quad\quad\quad \quad\quad\quad\quad \tilde{L}_q(s, w, E_{t, \chi} \times F) \frac{\Lambda_{t, \chi}(\ell; w)}{\ell^{w}}  h(t) \frac{dt}{2\pi},
\end{split}
\end{equation}  
where $\epsilon_f \in \{\pm 1\}$ is the parity of the Maa{\ss} form $f \in \mathcal{B}(q)$ and $\Re s, \Re w \not= 1$ in the expression for $\mathcal{M}^{\text{Eis}}_{q, \ell}(s, w; \mathfrak{h})$. The Dirichlet series expansions of    $L(s, \text{Ad}^2 f)$ and $L(s, f \times F)$ in $\Re s> 1$ (which serve as a definition for these functions) are given in  \eqref{Lad} and \eqref{L23}. We write 
\begin{equation}\label{moment-together}
\mathcal{M}^{\pm}_{q, \ell}(s, w, \mathfrak{h}) =  \mathcal{M}^{\text{Maa{\ss}}, \pm}_{q, \ell}(s, w; \mathfrak{h}) + \mathcal{M}^{\text{Eis}}_{q, \ell}(s, w; \mathfrak{h}) +  \delta_{\pm = +} \mathcal{M}^{\text{hol}}_{q, \ell}(s, w; \mathfrak{h}).
\end{equation}

Absolute convergence  follows from \eqref{weakly}, Weyl's law,  the convexity bound for the respective $L$-functions and the lower bounds \eqref{hl}. It is not hard to see (cf.\ Lemma \ref{analyticcont}) that under suitable circumstances $\mathcal{M}^{\text{Eis}}_{q, \ell}(s, w; \mathfrak{h})$ for $\Re s, \Re w > 1$ can be continued   to an $\varepsilon$-neighbourhood of $\Re s, \Re w \geq 1/2$, and this continuation equals $\mathcal{M}^{\text{Eis}}_{q, \ell}(s, w; \mathfrak{h})$ for $1/2 < \Re s, \Re w <1$ plus some polar terms. 



\begin{theorem}\label{thm1} Let $\mathfrak{h}$ be an admissible function, $(q, \ell) = 1$ and let $s, w\in \Bbb{C}$ be such that  
\begin{equation}\label{final-region}
1/2 \leq \Re s \leq \Re w  < 3/4.
\end{equation}
Let $F$ be either cuspidal or the minimal parabolic Eisenstein series\footnote{More general Eisenstein series could be treated in the same way, but we do not have applications in this case.} $\E_0$. Then
\begin{equation}\label{final-formula}
\begin{split}
\mathcal{M}^{+}_{q, \ell}(s, w; \mathfrak{h})  = \mathcal{N}_{q, \ell}(s, w; \mathfrak{h}) +  \sum_{\pm} \mathcal{M}^{\pm}_{\ell, q}(s', w'; \mathscr{T}_{s', w'}^{\pm}\mathfrak{h}), 
\end{split}
\end{equation}
where $s', w'$ are as in \eqref{new}, the ``main term'' $\mathcal{N}_{q, \ell}(s, w; \mathfrak{h})$  and the integral transform $\mathscr{T}^{\pm}_{s', w'}\mathfrak{h}$  
are given in \eqref{main-term-final} resp.\   \eqref{final-trafo}. 
The function $\mathscr{T}^{\pm}_{s', w'}\mathfrak{h}$  is weakly admissible, and for $\varepsilon > 0$ we have
 \begin{equation}\label{required-bound}
 \mathcal{N}_{q, \ell}(s, w; \mathfrak{h})\ll_{s, w, \mathfrak{h}, \varepsilon} \ell^{\theta - 1 + \varepsilon} + q^{\theta - 1 + \varepsilon}. 
 \end{equation}
  \end{theorem}


In the special case $q = \ell = 1$, $F = \E_0$, a beautiful   formula for the fourth moment was envisaged by Kuznetsov  in the late 1980s and completed by Motohashi \cite{Mo}.  This formula, however, has  different features. Our proof decomposes the fourth moment as $4 = 3+1$, whereas Kuznetsov and Motohashi use a decomposition of the form $4=2+2$. While more symmetric on the surface, their setup leads to 
  convergence problems that are overcome by considering a carefully designed \emph{difference} of the holomorphic and the Maa{\ss} spectrum. In particular, everywhere positive test functions $\mathfrak{h}$ as in Lemma \ref{lem1}c) are in general excluded. (As long as $\ell = q = 1$, one can exploit the fact that there are no holomorphic forms of small weight to circumvent this problem, but for more general values of $q$ and $\ell$ this device is not available.) This is probably the reason why this formula seems to not have been used for applications. 

Our formula should be seen as a contribution to the rich theory of automorphic forms  on ${\rm GL}(2) \times {\rm GL}(4)$ in the special situation where the ${\rm GL}(4)$ function is an isobaric $3+1$-sum. A corresponding reciprocity formula for $\mathcal{M}^-_{q, \ell}(s, w; \mathfrak{h})$ could also be established, in fact in an analytically simpler fashion (cf.\ \cite[end of Section 2]{Mo}). The structure of the argument indicates that a linear combination of these two formulae is an exact involution; it would be very interesting to provide a formal proof of this statement.

 A spectral reciprocity formula in  the  \emph{archimedean}  aspect -- i.e.\ the level is kept fixed, but the spectral parameter is varying -- for generic automorphic forms $F$ on ${\rm GL}(4)$   with an additional twist by the parity
$$ \sum_{t_f \leq T} \epsilon_f L(1/2, f \times F) \rightsquigarrow T \sum_{t_f \ll 1}  \epsilon_f L(1/2, f \times \tilde{F})$$
(where $\tilde{F}$ denotes the dual form) was considered 
by a  different technique   in   \cite{BLM} along with applications to non-vanishing of $L$-functions.   (In our set-up this would correspond to the easier, but excluded case $q= \ell$ in which case the Hecke eigenvalue is essentially the root number, at least as long as $q$ is squarefree.) Finally we mention that there is ``lower dimensional'' version of reciprocity for certain $L$-functions on ${\rm GL}(1) \times {\rm GL}(2)$:
$$\sum_{\chi \text{(mod }q)} |L(1/2, \chi)|^2 \chi(\ell) \rightsquigarrow \sum_{\chi \text{(mod }\ell)} |L(1/2, \chi)|^2 \chi(q).$$ 
 This was first observed by Conrey \cite{Co}, and extended by Young \cite{Yo} and Bettin \cite{Be}. It would be very interesting to investigate if there is a ``master formula'' that contains all previously mentioned spectral reciprocity formulae as special cases.

\subsection{Applications}\label{app} The most interesting applications of our reciprocity formula concern the twisted fourth moment \eqref{twisted}, which has a long history in the theory of $L$-functions. The investigation started with seminal work of Duke-Friedlander-Iwaniec \cite{DFI}, who used  their work on the binary additive divisor problem to obtain\footnote{This  is not stated explicitly, but at least for prime level follows very easily from their Corollary 2, for instance.}
\begin{equation}\label{DFI-bound} 
   \sum_{f } L(s, f)^4 \lambda_f(\ell) \ll_{k, s, \varepsilon} q^{1+\varepsilon} \ell^{-1/2} + q^{11/12+\varepsilon} \ell^{3/4}
\end{equation}
for $\Re s = 1/2$, where the sum is taken over holomorphic newforms of fixed weight $k$ and prime level $q$. 
Together with the amplification method they deduced the first subconvexity result 
\begin{equation}\label{subconvex}
L(s, f) \ll_{k, s} q^{\frac{1}{4} - \delta}, \quad \delta < 1/192, \quad \Re s = 1/2
\end{equation}
for this family. 
The then newly developed machinery was taken up by Kowalski-Michel-VanderKam \cite{KMV}, who  derived a precise asymptotic formula for the left hand side of \eqref{DFI-bound} in the case of weight $k=2$ and prime level $q$ with error term $O_{\varepsilon}(q^{11/12+\varepsilon} \ell^{3/4})$. This allowed them to deduce, via mollification, various non-vanishing results for central $L$-values. Quite recently, Balkanova and Frolenkov \cite{BF} improved the error term somewhat. On heuristic grounds, one would expect an error term $O(q^{1/2+\varepsilon})$, and one would expect that the limit of current machinery should be an error term $O(q^{1/2+\varepsilon} \ell^{1/2} )$, at least under the Ramanujan conjecture. This is precisely what we prove:
\begin{theorem}\label{thm2}  Let $q$ be prime, $\varepsilon > 0$, $(\ell, q) = 1$ and  $\mathfrak{h} = (h, h^{\text{{\rm hol}}})$ an admissible function. Then
$$\sum_{f \in \mathcal{B}^{\ast}(q)} \frac{L(1/2,  f)^4  }{L(1, \text{{\rm Ad}}^2 f)}   \lambda_f(\ell)  h(t_f) +   \sum_{f \in \mathcal{B}^{\ast}_{\text{{\rm hol}}}(q)} \frac{L(1/2,    f  )^4 }{L(1, \text{{\rm Ad}}^2 f)}  \lambda_f(\ell) h^{\text{{\rm hol}}}(k_f)
\ll_{\mathfrak{h}, \varepsilon}
(q \ell)^{\varepsilon}\left(q \ell^{-1/2} + q^{1/2 + \vartheta} \ell^{1/2}\right).$$
\end{theorem}
A slightly more general result  is given in Proposition \ref{prop-new} in Section \ref{sec-new}. 
The bound is optimal for $\ell$ up to size $q^{1/2}$ (modulo the Ramanujan conjecture).  Since this is the length of an approximate functional equation for $L(1/2, f)$, we can sum this bound trivially to obtain a bound for the fifth moment. 
\begin{theorem}\label{cor1} Let $q$ be prime, $\varepsilon > 0$.  Then
$$
\sum_{ f \in \mathcal{B}^{\ast} (q)  } L(1/2, f)^5 e^{-t_f^2} \ll_{\varepsilon} q^{1+ \vartheta + \varepsilon}.$$
\end{theorem}
For holomorphic forms of weight $k \in \{4, 6, 8, 10, 14\}$ an analogous result was recently obtained by K\i ral and Young \cite{KY} in an impressive 80-page tour de force argument. Our   version for Maa{\ss} forms comes   as a special case of a   more general framework. 
Since $L(1/2, f)$ is non-negative \cite{KaS}, Theorem  \ref{cor1} implies immediately the strong subconvexity bound
\begin{equation}\label{1/5subconvex}
L(1/2, f) \ll_{\varepsilon} q^{(1+\vartheta + \varepsilon)/5}.
\end{equation}
This is an example of ``self-amplification'', where the fourth moment is amplified by a fifth copy of the $L$-function itself. Working with a traditional amplifier, however, is a more robust technique, since one can exploit positivity more strongly. This gives
 \begin{theorem}\label{cor2}  Let $q$ be squarefree with $(6, q) = 1$, $\varepsilon > 0$  and $f \in \mathcal{B}^{\ast}_{\text{\rm hol}}(q) \cup  \mathcal{B}^{\ast}(q)$. Then
$$L(1/2, f) \ll_{\varepsilon} q^{\frac{1}{4} - \frac{1 - 2\vartheta}{24}+\varepsilon}$$
where  the implicit constant depends on the archimedean parameter of $f$, i.e.\ $t_f$ or $k_f$. 
\end{theorem}
Note that this holds both for Maa{\ss} forms and for holomorphic forms of arbitrary (fixed) weight, including weight $k=2$.  The assumptions on $q$ could be weakened   with more work. While asymptotically weaker than \eqref{1/5subconvex} for   $\vartheta\rightarrow 0$, with the current value of $\vartheta = 7/64$ the exponent $0.217$  is numerically a little stronger than $0.221$ in   \eqref{1/5subconvex} and appears to be the current record. The generic $1/24$-saving, which is 16.6 percent of the convexity bound, is the natural limit in situations where the Iwaniec-Sarnak amplifier based on the relation $\lambda_f(p)^2 = 1+ \lambda_f(p^2)$ is applied, cf.\ for instance the corresponding bound \cite{IS} in the sup-norm problem. As we cannot exclude the possibility of many small Hecke eigenvalues, a more efficient amplifier is not available unconditionally, but an optimal amplifier (whose existence is potentially a much deeper problem than the subconvexity problem discussed here) would yield the bound 
$$L(1/2, f) \ll_{\varepsilon} q^{\frac{1}{4} - \frac{1 - 2\vartheta}{16}+\varepsilon}$$
of Burgess quality, which again marks the limit of current technology.

\subsection{The methods}

The general strategy of the proof is simple to describe. We apply the Petersson and/or Kuznetsov spectral summation formula, followed by the ${\rm GL}(3)$ Voronoi formula. In this way, the Kloosterman sums become ${\rm GL}(1)$-exponentials of the form $e(n \bar{d}/c)$, to which we can apply the additive reciprocity formula 
\begin{equation}\label{rec}
   e\Bigl(\frac{n\bar{d}}{c}\Bigr)= e\left(-\frac{n\bar{c}}{d}\right) e\left(\frac{n}{cd}\right).
 \end{equation}  
    Now we reverse all transformations, i.e.\ we apply the ${\rm GL}(3)$ Voronoi formula in the other direction and the Kuznetsov formula backwards to arrive at a ``dualized'' spectral sum. The application of \eqref{rec} makes this procedure non-involutory and  yields the desired reciprocity formula. Note that in our application $c$ will be divisible by $q$ and $d$ by $\ell$, so the formula \eqref{rec} moves information from the primes dividing $q$  to the primes dividing $\ell$ as well as the infinite prime, because the right most exponential is treated as an archimedean weight function. 

Experience has shown that this type of simple back-of-an-envelope heuristics can lead to enormous technical challenges. 
In this paper we have tried to present an approach that exploits structural features (including higher rank tools) as much as possible and circumvents many of the well-known technical problems.

On the analytic side, we avoid  stationary phase arguments and intricate asymptotic analysis completely, but let instead analytic continuation (along with Stirling's formula) do the work for us. We do not start with an approximate functional equation (which also relieves us from any root number considerations), but work always with infinite Dirichlet series in the region of absolute convergence and postpone analytic continuation to the center of the critical strip to the last moment. 

On the arithmetic side, we design the set-up carefully to take care of the combinatorial difficulties of higher rank Hecke algebras. This was already a key device in \cite{BLM}. The ${\rm GL}(3)$ Voronoi formula (or -- roughly equivalently -- three instances of Poisson summation in the Eisenstein case) is conceptually well-understood, but introduces in practice a whole alphabet of auxiliary variables for additional divisibility and coprimality conditions. The multiple Dirichlet series that we introduce in Section \ref{multiple} is tailored to the features of the ${\rm GL}(3)$ Hecke algebra. The combinatorial difficulties do not completely disappear in this way; a shadow  remains in the rather complicated local correction factors $\tilde{L}_q(s, w, f \times F)$ and $\tilde{L}_q(s, w, E_{t, \chi} \times F)$, but the advantage of our approach is complete symmetry before and after the application of the arithmetic reciprocity formula \eqref{rec}. We would like to emphasize the strength of genuine higher rank tools (such as combinatorial and analytic properties of the  ${\rm GL}(3)$ Voronoi formula) even for applications that involve only ${\rm GL}(2)$ objects (such as Theorems \ref{thm2} -- \ref{cor2}), and it seems that the present paper is the first instance where higher rank tools are used for the analysis of the twisted fourth moment.

One of the major problems in earlier approaches \cite{KMV, KY} is the presence of very complicated main terms coming from various sources, the most difficult being the various zero frequencies in the multiple applications of the Poisson summation formula. These terms are complicated by nature, but in our approach they arise in a very simple and conceptual way twice as a single residue of a certain multiple Dirichlet series. 


It may be instructive to compare the general strategy with other options used in earlier works. As mentioned before, the analysis of \cite{DFI} and \cite{KMV} is based on a decomposition $4 = 2+2$ and a treatment of the binary additive divisor problem. The formula of Kuznetsov and Motohashi \cite{Mo} also decomposes $4 = 2+2$ and dualizes both ${\rm GL}(2)$ components. Unfortunately this leads to an essentially self-dual deadlock situation that is resolved -- as mentioned earlier -- by subtracting the holomorphic spectrum from the Maa{\ss} spectrum in a carefully designed way. K\i ral and Young \cite{KY} instead 
use the decomposition $5=3+1+1$ in the context of the fifth moment. 
They first dualize one of the ${\rm GL}(1)$ components, apply arithmetic reciprocity \eqref{rec} and then dualize the ${\rm GL}(3)$ component. A major technical difficulty is the fact that this leads to Kloosterman sums with multiplicative inverses in the arguments, so that the Kuznetsov formula at general cusps for possibly non-squarefree levels is required. In contrast, our approach never touches the ${\rm GL}(1)$ components, but dualizes the ${\rm GL}(3)$ component twice. This more symmetric set-up avoids complications with the Kuznetsov formula and makes a classical amplification procedure possible, which seems to be hard to implement in the work of \cite{KY}. The idea of coupling a $4=3+1$ structure along with arithmetic reciprocity \eqref{rec} goes back to important work of X.\ Li \cite{Li} and was first used for varying levels by the second author in \cite{Kh}. 

\subsection{A roadmap} Sections \ref{basic} -- \ref{voronoi-formula} contain essentially known material, in particular the Kuznetsov formula and the Voronoi summation formula tailored in suitable form for later reference. An important technical ingredient in Section \ref{kuz-formula} is a detailed Fourier expansion of a complete orthonormal basis (including oldforms) in terms of Hecke eigenvalues.  Lemma \ref{final-decay} will be used at the very end of the argument to verify that the integral transform $\mathscr{T}^{\pm}_{s', w'}\mathfrak{h}$ is weakly admissible. To treat holomorphic functions of small weight as in Theorem \ref{cor2}, we use a combination of the Petersson and the Kuznetsov formula and certain special functions that are also introduced in Section \ref{kuz-formula}. 

Section \ref{multiple} features a multiple Dirichlet series in three variables that contains both ${\rm GL}(3)$ Fourier coefficients and Kloosterman sum. It is carefully designed as the combinatorial hinge between the Kuznetsov formula and the Voronoi formula, so that both formulas can be applied without extra technical complications. If $F= \E_0$ is an Eisenstein series, we take some time to compute the Laurent expansion of the triple pole in one of the variables that will contribute to the main term $\mathcal{N}_{q, \ell}(s, w; \mathfrak{h})$ in Theorem \ref{thm1}.  These are related to the ``fake main terms'' in \cite{KY}. After a somewhat involved computation it turns out that the Laurent coefficients (as a function of the two other variables) have a beautiful formulation as a quotient of Riemann zeta-functions and their derivatives. 

Section \ref{admissible} gives a precise definition of admissible functions and establishes some technical properties of these functions. Section \ref{int-trafo} is also of analytic nature and studies in detail a certain integral transform that will later become the central part of the   transform  $\mathscr{T}_{s', w'}^{\pm}$ in \eqref{final-trafo}. The key point here is to obtain analytic continuation and suitable decay properties which is achieved by careful contour shifts. 

Section \ref{prelim-formula} features a prototype of the reciprocity formula based on the triad  Voronoi-reciprocity-Voronoi. At this point we still work with points $s, w$ whose real part is sufficiently large. The use of Kuznetsov formula at the beginning and at the end to obtain the full 5-step procedure outlined in the previous subsection, as well as the analytic continuation to $\Re w \geq \Re s \geq 1/2$, is postponed to Section \ref{proof1} where the proof of Theorem \ref{thm1} is completed. 
The necessary local computations to deal with the ramified primes $p \mid \ell q$ are provided in Section \ref{local}.  It is then a simple task to derive Theorems \ref{thm2} -- \ref{cor2} in the final two sections. 

Finally we mention that while the abstract procedure Kuznetsov-Voronoi-reciprocity-Voronoi-Kuznetsov is completely symmetric, there are certain differences on a technical level before and after the reciprocity formula. In the latter case, some local factors at primes $p \mid \ell q$ and the archimedean place converge absolutely only slightly to the right of $\Re w \geq \Re s \geq 1/2$. To get the desired analytic continuation, some extra maneuvers are necessary, as can be seen in the proof of Lemma \ref{pole2} and in particular in Lemma \ref{sec-prop} that will be applied both in Sections \ref{prelim-formula} and \ref{proof1}.

\subsection{Notation and conventions} Throughout, the letter $\varepsilon$ denotes an arbitrarily (and sufficiently) small positive real number, not necessarily the same at each occurrence. All implied constants may depend on $\varepsilon$ (where applicable), but this is  suppressed from the notation. We will frequently encounter multiple sums and integrals of holomorphic functions in one or more variables. All expressions of this type will be absolutely convergent by which we mean in addition locally uniformly  convergent in the auxiliary variables, so that they represent again holomorphic functions. Also implied constants depending on complex variables are always understood to depend locally uniformly on them. 
 By an $\varepsilon$-neighbourhood of a strip $c_1 \leq \Re s \leq c_2$ we mean the open set $c_1 -\varepsilon < \Re s < c_2 + \varepsilon$, and similarly for multidimensional tubes. Occasionally we will  encounter meromorphic functions in multidimensional tubes  which will be holomorphic outside a finite set of polar divisors given by affine hyperplanes.   
By $v_p$ we denote the usual $p$-adic valuation. We write $a \mid b^{\infty}$ to mean that $a$ divides some power of $b$. We write $\Bbb{N}_0 = \Bbb{N}\cup \{0\}$.

\textbf{Acknowledgements.} The authors would like to thank Matthew Young and Guohua Chen as well as the referees for useful comments and suggestions that helped improving and correcting various aspects of teh paper. 

\section{Basic analysis}\label{basic}

\subsection{Mellin transform}  Throughout this paper we denote by 
\begin{equation}\label{mellin1}
\widehat{W}(s) = \int_0^{\infty} W(x) x^s \frac{dx}{x}
\end{equation}
 the Mellin transform of a   function $W : [0, \infty) \rightarrow \Bbb{C}$ of class $C^J$ for some $J \in \Bbb{N}_0$ satisfying $x^j W^{(j)}(x) \ll_{J, a, b} \min(x^{-a}, x^{-b})$ for some $-\infty < a < b < \infty$   and all $0 \leq j \leq J $. The function $\widehat{W}$ is then (initially) defined in $a < \Re s < b$ as an absolutely convergent integral and satisfies $\widehat{W}(s) \ll (1 + |s|)^{-J}$ in this region (in some cases it may be continued meromorphically to a larger region). 

The inverse Mellin transform of a  function $\mathcal{W}$ that is holomorphic in a strip containing $a \leq \Re s \leq b$ and bounded by $\mathcal{W}(s) \ll (1 +|s|)^{-r}$ for some $r > n \in \Bbb{N}$, is given by
\begin{equation}\label{mellin2}
\widecheck{\mathcal{W}}(x) = \int_{(c)} \mathcal{W}(s) x^{-s} \frac{ds}{2\pi i},
\end{equation}
where here and in the following $(c)$ denotes the line $\Re s = c$ with $a \leq c \leq b$. We have 
\begin{equation}\label{mellin}
x^j\widecheck{\mathcal{W}} {}^{(j)}(x) \ll \min(x^{-a}, x^{-b}),
\end{equation}
for $j = 0, 1, \ldots, n-1$.

\subsection{The gamma function} We recall the reflection, recursion and duplication formula $$\Gamma(s)\Gamma(1-s)= \pi \sin(\pi s)^{-1}, \quad \Gamma(s+1) = s\Gamma(s), \quad \Gamma(s)\Gamma(s+1/2) = \sqrt{\pi} 2^{1-2s} \Gamma(2s).$$ 
For fixed $\sigma \in \Bbb{R}$, real $|t| \geq 3$, and any $M > 0$ we recall Stirling's formula
\begin{equation}\label{stir}
  \Gamma(\sigma + it) = e^{-\frac{\pi}{2}|t|} |t|^{\sigma-\frac{1}{2}} \exp\left(i t \log \frac{|t|}{e}\right)g_{\sigma, M}(t) + O_{\sigma, M}(|t|^{-M}),
\end{equation}
where
\begin{equation}\label{stir1}
  g_{\sigma, M}(t) = \sqrt{2\pi} \exp\left(\frac{\pi}{4}(2\sigma-1) i \, \text{sgn}(t) \right) +  O_{\sigma, M}\left(|t|^{-1}\right)
\end{equation}
and also
\begin{equation}\label{above}
|t|^j  
g^{(j)}_{\sigma, M}(t) \ll_{j, \sigma, M} 1
\end{equation}
for all fixed $j \in \Bbb{N}_0$. This implies in particular bounds of the shape
$$\frac{\Gamma(as + b)}{\Gamma(cs + d)} \ll_{a, b, c, d, \sigma} (1+ |t|)^{(a-c)\sigma + b-d}$$
for $a, b, c, d \in \Bbb{R}$, $s = \sigma + it$ and $\min_{n \in \Bbb{N}_0} |as + b + n| \geq 1/10$.

We define $\Gamma_{\Bbb{R}}(s) = \pi^{-s/2} \Gamma(s/2)$ and for $j \in \{0, 1\}$ we write
$$ G_j(s) = \frac{\Gamma_{\Bbb{R}}(s + j)}{\Gamma_{\Bbb{R}}(1-s + j)} = 2(2\pi)^{-s}\Gamma(s) \begin{cases}  \cos(\pi s/2), & j = 0,\\  \sin(\pi s/2), & j = 1,\end{cases}$$
where the last identity follows from the reflection and duplication formula of the gamma function. We also need the linear combination
\begin{equation}\label{Gpm}
G^{\pm}(s) = \frac{1}{2}G_0(s) \pm \frac{i}{2} G_1(s) =  \Gamma(s) (2\pi)^{-s} \exp(\pm i\pi s/2), 
\end{equation}
which by \eqref{stir} and \eqref{above}  satisfies for $|t| \geq 3$ 
 the bound 
\begin{equation}\label{bound-g}
  G^{\pm}(s) \ll (1 + |s|)^{\Re s - 1/2} e^{-\pi \max(0, \pm \Im s)/2}
\end{equation}
and the asymptotic formula
\begin{equation}\label{asymp-g}
  G^{\pm}(\sigma + it) = |t|^{\sigma - 1/2}   \exp\Bigl(i t \log \frac{|t|}{2\pi e}\Bigr) \tilde{g}_{\sigma, M}(t) + O_{\sigma, M}\left(|t|^{-M}\right)
\end{equation}
with $\tilde{g}_{\sigma, M}$ satisfying \eqref{above}. Note that $G^{\pm}(\sigma + it)$ is exponentially decaying for $\pm t \rightarrow \infty$ by \eqref{bound-g}. 

From the theory of hypergeometric  integrals we quote the following formula \cite[(2.2.1.2)]{PBM}: if $a, b, c, d\in \Bbb{C}$ satisfy $\Re(c+d)- 1> \Re(a+b) > 0$, then 
\begin{equation}\label{hyper}
\int_{(\sigma)} \frac{\Gamma(a + s) \Gamma(b-s)}{\Gamma(c+s)\Gamma(d-s)} \frac{ds}{2\pi i} = \frac{\Gamma(a+b)\Gamma(c+d-1-a-b)}{\Gamma(c-a)\Gamma(d-b)\Gamma(d+c-1)}
\end{equation}
where the integral on the left   is absolutely convergent and the path of integration separates the poles, i.e.\ $\Re b > \sigma > -\Re a$, for instance $\sigma = \frac{1}{2}(b-a)$. 

\subsection{Integration by parts} We quote   a useful integration by parts lemma  \cite[Lemma 8.1]{BKY}, which follows immediately from \cite[(8.6)]{BKY}.

\begin{lemma} \label{integrationbyparts}
 Let $Y \geq 1$, $X, Q, U, R > 0$, $B \in \Bbb{N}$, 
and suppose that $w$ 
 is a smooth function with support on some interval $[\alpha, \beta]$, satisfying
\begin{equation*} 
w^{(j)}(t) \ll X U^{-j}
\end{equation*}
for $0 \leq j \leq B$. 
Suppose $H$ 
  is a smooth function on $[\alpha, \beta]$ such that
\begin{equation*}
 |H'(t)| \gg R, \quad 
H^{(j)}(t) \ll Y Q^{-j} \,\, \text{for } j=2, 3, \dots, B. 
\end{equation*}
Then 
\begin{equation*}
I = \int_{\Bbb{R}} w(t) e^{i H(t)} dt
 \ll_B (\beta - \alpha) X \left[(QR/\sqrt{Y})^{-B/2} + (RU)^{-B/2}\right].
\end{equation*}
\end{lemma}

\section{The Kuznetsov formula}\label{kuz-formula}

Let $N \in \Bbb{N}$ and write
\begin{equation}\label{0}
 N \nu(N) := [\Gamma_0(1) : \Gamma_0(N)], \quad \text{i.e.} \quad \nu(N) = \prod_{p \mid N} \left(1 + \frac{1}{p}\right). 
\end{equation} 
We equip $\Gamma_0(N) \backslash \Bbb{H}$ with the inner product
\begin{equation}\label{inner}
\langle f, g\rangle := \int_{\Gamma_0(N)\backslash \Bbb{H}} f(z) \bar{g}(z) \frac{dx\, dy}{y^2}.
\end{equation}
The Kuznetsov/Petersson formulae require a sum over an orthonormal basis of automorphic forms of a given level $N$. While vectors in different representation spaces of $Z(\Bbb{A}_{\Bbb{Q}}) {\rm GL}_2(\Bbb{Q})\backslash {\rm GL}_2(\Bbb{A}_{\Bbb{Q}})/K_0(N)$ are always orthogonal, such a space may contain more than one    $L^2$-normalized vector of given $K_{\infty}$-type. Classically this corresponds to the existence of oldforms.  In the following we give explicit formulae for the Fourier coefficients of a complete orthonormal basis. 

\subsection{Fourier expansion of Eisenstein series}   The unitary Eisenstein spectrum of $\Gamma_0(N)$ can be parametrized by a continuous parameter  $s = 1/2 + it$  together with pairs $(\chi, M)$, where $\chi$ is a primitive Dirichlet character of conductor $c_{\chi}$ and $M \in \Bbb{N}$ satisfies $c_{\chi}^2 \mid M \mid N$. For our purposes this adelic parametrization is more convenient than the classical parametrization by cusps. For fixed $t$ and $\chi$ the Eisenstein series for various $M$  belong to the same representation space.  The corresponding spectral decomposition is proved in \cite{GJ}, and the Fourier expansion of these Eisenstein series is computed explicitly in \cite[Section 5]{KL} (see \cite[Section 2.7]{BH} for similar calculations over number fields). In the following, we quote from \cite{KL}. 
In the notation of \cite{KL} we have $M = \prod p^{i_p}$. 
We define as in \cite[(5.22)]{KL}
\begin{equation}\label{1}
\mathfrak{n}^2(M) :=  \frac{1}{M} \prod_{\substack{p \mid N\\ p \nmid (M, N/M)}} \frac{p}{(p+1)} \prod_{ p \mid (M, N/M)} \frac{p-1}{p+1} = \frac{1}{M} \tilde{ \mathfrak{n}}^2(M),
\end{equation}
say. 
Next,  as in \cite[Section 5.5]{KL} we define
\begin{equation*}
N_1 = \prod_{p \mid N/M} p^{v_p(N)}.  
\end{equation*}
The normalized Eisenstein series $E_{\chi, M, N}(z, s)$ of level $N$ corresponding to   $(\chi, M)$ has the Fourier expansion $$E_{\chi, M, N}(z, 1/2 + it) = \rho^{(0)}_{\chi, M, N}(t, y) + \frac{2 \pi^{1/2 + it} y^{1/2}}{\Gamma(1/2 + it)} \sum_{n \not= 0} \rho_{\chi, M, N}(n, t) K_{it}(2 \pi |n| y) e(nx),$$ where for $n \not= 0$ we have 
$$\rho_{\chi, M, N}(n, t)  = \frac{C(\chi, M)}{(N \nu(N))^{1/2}\mathfrak{n}(M)  L^{(N)}(1 + 2it, \chi^2)} \cdot \frac{|n|^{it}}{M^{1 + 2it}} \sum_{\substack{c \mid n \\ (c, N_1) = 1}} \frac{\chi(c)}{c^{2it}} \sum_{\substack{d \, (\text{mod }M)\\ (d, M) = 1}} \chi(d) e\left(\frac{d n/c}{M}\right),$$
for a constant $|C(\chi, M)| = 1$. As usual, the superscript $L^{(N)}$ denotes that the Euler factors at primes dividing $N$ are omitted. This follows from \cite[(5.32), (5.34)]{KL} after re-normalizing by $\mathfrak{n}(M)$, cf.\ \cite[(5.22)]{KL}, and taking into account that the adelic inner product in \cite{KL} differs from \eqref{inner} by a factor \eqref{0}. We write 
\begin{equation}\label{3}
   M = c_{\chi} M_1 M_2, \quad \text{where} \quad (M_2, c_{\chi}) = 1, \quad M_1 \mid c_{\chi}^{\infty},  \end{equation}  
    so that $c_{\chi} \mid M_1$ and $(M_1, M_2) = 1$. Then by the Chinese Remainder Theorem the $d$-sum equals
$$r_{M_2}(n/c) \sum_{\substack{d \, (\text{mod } c_{\chi}M_1)\\ (d, c_{\chi})= 1}} \chi(d) e\left(\frac{d \overline{M_2} (n/c)}{c_{\chi} M_1}\right) =  \delta_{cM_1 \mid n} r_{M_2}(n/c)\chi(M_2) \bar{\chi}\left(\frac{n}{cM_1}\right)M_1   
\sum_{\substack{d \, (\text{mod }  c_{\chi}) \\ (d, c_{\chi})= 1}}  \chi(d) e\left(\frac{d}{c_{\chi}}\right), $$
where the Gau{\ss} sum on right hand side has absolute value $c_{\chi}^{1/2}$ and 
\begin{equation}\label{rama}
   r_M(n) = \sum_{\substack{d \, (\text{mod } M)\\ (d, M) = 1}} e\left(\frac{dn}{M}\right) =  \sum_{d \mid (n, M)} d \mu\left(\frac{M}{d}\right)
   \end{equation}
    is the Ramanujan sum. Recalling \eqref{1}, 
we conclude
 $$\rho_{\chi, M, N}(n, t) = \frac{C(\chi, M, t)|n|^{it}}{(N \nu(N))^{1/2}\tilde{ \mathfrak{n}}(M) L^{(N)}(1 + 2it, \chi^2)}  \left(\frac{M_1}{M_2}\right)^{1/2} \sum_{\substack{cM_1 \mid n \\ (c, N_1) = 1}} \frac{\chi(c)}{c^{2it}}   r_{M_2}(n/c)\bar{\chi}\left(\frac{n}{cM_1}\right) $$
 where $|C(\chi, M, t)| = 1$. We note that  the condition $(c, N_1) = 1$ is equivalent to $(c, N/M) = 1$.   Finally we insert the expression on the right of \eqref{rama} for the Ramanujan sum getting the final formula
  \begin{equation}\label{rho-eis}
    \rho_{\chi, M, N}(n, t)  = \frac{C(\chi, M, t)|n|^{it}}{(N \nu(N))^{1/2}\tilde{ \mathfrak{n}}(M) L^{(N)}(1 + 2it, \chi^2)}  \left(\frac{M_1}{M_2}\right)^{1/2}  \sum_{\delta \mid M_2} \delta \mu(M_2/\delta)  \bar{\chi}(\delta)  \sum_{\substack{cM_1\delta f =  n \\ (c, N/M) = 1}}\frac{\chi(c)}{c^{2it}}   \bar{\chi}(f) 
    \end{equation}
  with the notation \eqref{0}, \eqref{1},   \eqref{3}. In particular, using the third bound of \eqref{hl} below, we find for $t \in \Bbb{R}$ that 
  \begin{equation}\label{rho-bound}
    \rho_{\chi, M, N}(n, t)\ll ((1 + |t|)N|n|)^{\varepsilon} \frac{(M_1M_2)^{1/2}}{N^{1/2}} \leq ((1 + |t|)N|n|)^{\varepsilon} .
  \end{equation}

\subsection{Fourier expansion of cusp forms} 
The cuspidal spectrum is parametrized by  pairs $(f, M)$ of $\Gamma_0(N)$-normalized newforms $f$ of level $N_0 \mid N$ and integers $ M \mid N/N_0$. This comes from orthonormalizing the set $\{z \mapsto f(Mz) : M \mid N/N_0\}$ (which belongs to the same representation space) by Gram-Schmidt\footnote{Of course, there are many ways to apply the Gram-Schmidt procedure. For concreteness, we choose the one used in \cite{BM} and referenced below; this fixes our basis uniquely}. 
If $f$ is a Maa{\ss} form, we write the Fourier expansion of the pair $(f, M)$ as
$$ (2 \cosh(\pi t) y)^{1/2} \sum_{n \not= 0} \rho_{f, M, N}(n)   K_{it}(2 \pi |n| y) e(nx).$$
A standard Rankin-Selberg computation shows (verbatim as in \cite[(2.8) - (2.9)]{Bl}) that for $M = 1$, i.e.\ for newforms, one has 
\begin{equation}\label{rho1}
|\rho_{f, 1, N}(1)|^2 = \frac{1}{L(1, \text{Ad}^2 f)N\nu(N)} \prod_{p \mid N_0} \left(1 - \frac{1}{p^2}\right)
\end{equation} and $\rho_{f, 1, N}(n) = \rho_{f,1, N}(1) \lambda_f(n)$ for $n \in \Bbb{N}$. Here we use the notation 
\begin{equation}\label{Lad}
L(s, \text{Ad}^2 f) = \zeta^{(N)}(2s) \sum_n \frac{\lambda_f(n^2)}{n^s} 
\end{equation}
(which may differ from the corresponding Langlands $L$-function by finitely many Euler factors). 
 For general $M$, the Fourier coefficients of the orthonormal basis were computed in\footnote{The conditions $(d, N_0) = 1$ in the definition of $r_f$ and  $(b, N_0) = 1$ in the definitions of $\alpha$ and $\beta$ were erroneously missing there, see {\tt www.uni-math.gwdg.de/blomer/corrections.pdf}.}   \cite[Lemma 9]{BM}. Define
\begin{displaymath}
\begin{split}
& r_{f}(c) = \sum_{b \mid c} \frac{\mu(b) \lambda_f(b)^2}{b} \Bigl(\sum_{\substack{d \mid b\\ (d, N_0) = 1}}\frac{1}{d} \Bigr)^{-2}, \quad \alpha(c) = \sum_{\substack{b \mid c \\ (b, N_0) = 1}} \frac{\mu(b)}{b^2}, \quad \beta(c) = \sum_{\substack{b \mid c\\ (b, N_0) = 1} } \frac{\mu^2(b)}{b^2},\end{split}
\end{displaymath}
define $\mu_f(c)$ as the Dirichlet series coefficients of $L(s, f)^{-1} $ and let
$$\xi'_{f}(M, d) = \frac{\mu(M/d) \lambda_f(M/d)}{r_{f}(M)^{1/2}   \beta(M/d)}, \quad \xi''_{f}(M, d) = \frac{\mu_f(M/d) }{r_{f}(M)^{1/2}   \alpha(M)^{1/2}}.$$
Write $M = M_1 M_2$ where $M_1$ is squarefree, $M_2$ is squarefull and $(M_1, M_2)=1$. Then
\begin{equation}\label{xi-arithmetic}
\xi_{f}(M, d) := \xi_{f}'(M_1, (M_1, d)) \xi_{f}''(M_2, (M_2, d)) \ll M^{\varepsilon}(M/d)^{\vartheta}
\end{equation}
and 
\begin{equation}\label{rho-cusp}
 \rho_{f, M, N}(n) = \frac{1}{L(1, \text{Ad}^2 f)^{1/2}(N\nu(N))^{1/2}} \prod_{p \mid N_0} \left(1 - \frac{1}{p^2}\right)^{1/2}  \sum_{d \mid M} \xi_{f}(M, d) \frac{d}{M^{1/2}} \lambda_f(n/d)  
 \end{equation}
 for $n \in \Bbb{N}$ 
with the convention $\lambda_f(n) = 0$ for $n \not \in \Bbb{Z}$. For $-n \in \Bbb{N}$ we have $\rho_{f, M, N}(n) = \epsilon_f \rho_{f, M, N}(-n)$ where $\epsilon_f \in \{\pm 1\}$ is the parity of $f$. Using a standard Rankin-Selberg bound as well as the standard lower bound \eqref{hl} below, we obtain
\begin{equation}\label{rho-cusp-bound}
  \sum_{n \leq x}  |\rho_{f, M, N}(n) |^2 \ll \frac{1}{N}(Nx(1 + |t_f|))^{\varepsilon} x .
\end{equation}

If $f$ is a holomorphic newform of weight $k$ and level $N_0 \mid N$, we write the Fourier expansion of the pair $(f, M)$ as
$$\left(\frac{2\pi^2}{\Gamma(k)}\right)^{1/2}\sum_{n > 0} \rho_{f, M, N}(n) (4\pi n)^{(k-1)/2} e(nz).$$
If $M=1$, then  \eqref{rho1} remains true (cf.\ e.g.\  \cite[(2.1)]{BKY}), and so does \eqref{rho-cusp} for  $n \in \Bbb{N}$ and arbitrary $M \mid N/N_0$ as well as \eqref{rho-cusp-bound} with $k_f$ in place of $(1 + |t_f|)$.  For negative $n$ we define $\rho_{f, M, N}(n)  = 0$. 
 
\subsection{Versions of the Kuznetsov formula}
For $x > 0$ we define the integral kernels 
\begin{displaymath}
\begin{split}
& \mathcal{J}^+(x, t) := \frac{\pi  i}{ \sinh(\pi t)} (J_{2 it}(4 \pi x) - J_{-2it}(4 \pi x)), \\
& \mathcal{J}^-(x, t) := \frac{\pi  i}{ \sinh(\pi t)} (I_{2 it}(4 \pi x) - I_{-2it}(4 \pi x)) = 4 \cosh(\pi t) K_{2it}(4\pi x),\\
& \mathcal{J}^{\text{hol}}(x, k) := 2\pi i^k J_{k-1}(4\pi x) =  \mathcal{J}^+(x, (k-1)/(2i)) , \quad k \in 2\Bbb{N}. 
 \end{split}
 \end{displaymath}
 For future reference we record the Mellin transforms of these kernels:
 \begin{equation}\label{mellin-j}
\begin{split}
 \widehat{\mathcal{J}^+(., t)}(u) & = \frac{\pi i (2\pi)^{-u} }{2\sinh(\pi t)} \left(\frac{\Gamma(u/2 + it)}{  \Gamma(1 - u/2 + it)}  - \frac{\Gamma(u/2 - it)}{  \Gamma(1 - u/2 - it)} \right) \\
&  = (2\pi)^{-u} \Gamma(u/2 + it)\Gamma(u/2 - it) \cos(\pi u /2), \\
 \widehat{\mathcal{J}^-(., t)}(u) & =   (2\pi)^{-u} \Gamma(u/2 + it)\Gamma(u/2 - it) \cosh(\pi t),\\
 \widehat{\mathcal{J}^{\text{hol}}(., k)}(u) & = i^k(2\pi)^{-u} \pi \Gamma( (u + k-1)/2) \Gamma( (1+k-u)/2)^{-1}. 
 \end{split}
 \end{equation}
 These formulae follow from \cite[17.43.16 \& 18]{GR} together with the reflection formula for the gamma function. 
 
Let $h$ be an even  function $h$ satisfying  $h(t) \ll (1+|t|)^{-2-\delta}$ for $t \in \Bbb{R} \cup [-i\vartheta, i\vartheta]$, and let $h^{\text{hol}} : 2\Bbb{N} \rightarrow \Bbb{C}$ be a function satisfying $h^{\text{hol}}(k) \ll k^{-2-\delta}$ for some $\delta > 0$.  
 Then for $q\in \Bbb{N}$, $n, m \in \Bbb{Z} \setminus \{0\}$ we define 
  \begin{equation*}
\begin{split}
& \mathcal{A}_q^{\text{Maa{\ss}}}(n, m; h) := \sum_{q_0M \mid q}   \sum_{f \in \mathcal{B}^{\ast}(q_0)} \rho_{f, M, q}(n) \overline{ \rho_{f, M, q}(m)} h(t_f),\\
&  \mathcal{A}_q^{\text{Eis}}(n, m; h) := \sum_{c_{\chi}^2\mid M \mid q} \int_{\Bbb{R}} \rho_{\chi, M, q}(n, t) \overline{ \rho_{\chi, M, q}(m, t)} h(t) \frac{dt}{2\pi},\\
& 
\mathcal{A}_q^{\text{hol}}(n, m; h^{\text{hol}}) :=  \sum_{q_0M \mid q}   \sum_{ f \in \mathcal{B}_{\text{hol}}^{\ast}(q_0) } \rho_{f, M, q}(n) \overline{ \rho_{f, M, q}(m)} h^{\text{hol}}(k_f).
 \end{split}
\end{equation*}
(Note that by definition $A_q^{\text{hol}}(n, m; h^{\text{hol}}) = 0$ if $n$ or $m$ are negative.) If in addition $h$ is holomorphic in an $\varepsilon$-neighbourhood of $|\Im t| \leq 1/2$ and still satisfies the bound $h(t) \ll (1+|t|)^{-2-\delta}$ in this region, then for $n, m \in \Bbb{N}$ the Bruggeman-Kuznetsov formula states\footnote{Classical references parametrize the Eisenstein series by cusps, but the same proof works of course with the spectral expansion of Eisenstein series given by \cite{GJ}.}    (e.g.\ \cite[Theorem 16.3]{IK}, cf.\ \cite[(16.19)]{IK} for the normalization there)
  \begin{equation}\label{kuz1} 
\begin{split}
 \mathcal{A}_q^{\text{Maa{\ss}}}(n, m; h) +  \mathcal{A}_q^{\text{Eis}}(n, m; h)
&  =  \delta_{n, m} \int_{-\infty}^{\infty} h(t) \frac{ t \tanh(\pi t) dt}{2\pi^2}  + \sum_{q \mid c} \frac{S(n, m, c)}{c}\mathscr{K}h \Bigl(\frac{\sqrt{nm}}{c}\Bigr),
\end{split}
\end{equation} 
where
\begin{equation}\label{H}
 \mathscr{K}h(x) :=   \int _{-\infty}^{\infty} \mathcal{J}^+(x, t) h(t)   t \tanh(\pi t) \frac{dt}{2\pi^2}
  \end{equation}
  and $S(m, n, c) \ll (m, n, c)^{1/2} \tau(c) c^{1/2}$ is the Kloosterman sum. Absolute convergence of the $c$-sum in \eqref{kuz1} follows  from shifting the contour in \eqref{H} to $\Im t = \pm 3/8$, say.  
The formula \eqref{kuz1} is complemented by the Petersson formula \cite[Proposition 14.5]{IK}
 \begin{equation}\label{pet}
 \begin{split}
& 
\mathcal{A}_q^{\text{hol}}(n, m; \delta_{k = k_0}) = \frac{k_0-1}{2\pi^2} \Bigg(\delta_{n, m} +  2\pi i^{-k_0} \sum_{q\mid c} \frac{S(m, n, c)}{c} J_{k_0-1}\Big(\frac{4\pi \sqrt{mn}}{c}\Big)\Bigg)
 \end{split}
  \end{equation}
for $k_0 \in 2\Bbb{N}$, $n, m \in \Bbb{N}$.   Sometimes it is useful to apply the Petersson formula \eqref{pet} and   the Kuznetsov formula \eqref{kuz1} simultaneously. With this in mind, for a pair of functions $\mathfrak{h} = (h,  h^{\text{hol}})$ satisfying the above conditions, we define 
\begin{equation}\label{kast}
\mathscr{K}^{\ast} \mathfrak{h}(x) :=  \mathscr{K} h(x)  +  \sum_{k\in 2\Bbb{N}}i^{-k}  \frac{k-1}{\pi}h^{\text{hol}}(k) J_{k-1}(4\pi x).
\end{equation}
Then
  \begin{equation}\label{kuz-all} 
\begin{split}
\mathcal{A}_q(n, m; \mathfrak{h}) & :=   \mathcal{A}_q^{\text{Maa{\ss}}}(n, m; h) +  \mathcal{A}_q^{\text{Eis}}(n, m; h)+   \mathcal{A}_q^{\text{hol}}(n, m; h^{\text{hol}}) \\
&  =   \delta_{n, m} \mathscr{N}\mathfrak{h}  + \sum_{q \mid c} \frac{S(n, m, c)}{c}\mathscr{K}^{\ast}\mathfrak{h} \Bigl(\frac{\sqrt{nm}}{c}\Bigr),
\end{split}
\end{equation} 
for $n, m \in \Bbb{N}$ with
$$\mathscr{N}\mathfrak{h} :=  \int_{-\infty}^{\infty} h(t)   t \tanh(\pi t) \frac{dt}{2\pi^2} +\sum_{k\in 2\Bbb{N}} \frac{k-1}{2\pi^2} h^{\text{hol}}(k).$$  

Conversely,  if $H \in C^3((0, \infty))$ satisfies\footnote{There are various assumptions in the literature, e.g. \cite[(16.38)]{IK}, \cite[Theorem 2]{Ku}, \cite[(2.4.6)]{Mo}. We follow the latter, although the precise exponents make no difference for our argument.}  $x^j H^{(j)}(x) \ll \min(x, x^{-3/2})$ for $0 \leq j \leq 3$ and $n, m, q \in \Bbb{N}$, then  we have \cite[Theorem 16.5]{IK}
 \begin{equation}\label{kuz2} 
\begin{split}
\sum_{q \mid c} &\frac{S(\pm n, m, c)}{c}H \left(\frac{\sqrt{nm}}{c}\right) =  \mathcal{A}_q(\pm n, m; \mathscr{L}_{\pm} H) 
 \end{split}
\end{equation}
where 
\begin{equation}\label{Hback}
\begin{split}
\mathscr{L}_{\pm}H = (\mathscr{L}^{\pm}H, \mathscr{L}^{\text{hol}}H), \quad \mathscr{L}^{\diamondsuit}H = \int_{0}^{\infty} \mathcal{J}^{\diamondsuit}(x, .) H(x) \frac{dx}{x}. 
\end{split}
\end{equation}
for $\diamondsuit \in \{+, -, \text{hol}\}$. 

The formulas \eqref{kuz-all} and \eqref{kuz2} are inverses to each other, and so are the corresponding integral transforms: for $H \in C^2((0, \infty))$ with $H^{(j)}(x) \ll \min(x^{1/2}, x^{-5/2})$ for $j = 0, 1, 2$ we have the Sears-Titchmarsh inversion formula (cf.\ \cite[(4.9)]{ST} or \cite[(A.4)]{Ku}) 
\begin{equation}\label{inv}
\mathscr{K}^{\ast} \mathscr{L}_+ H = \mathscr{K}^{\ast}(\mathscr{L}^+  H, \mathscr{L}^{\text{hol}}H) = H.
 \end{equation}

To treat  holomorphic cusp forms of small weight, we use the following special functions, borrowed from \cite[Section 2]{BHM}.  For integers $3 < b < a$ with $a \equiv b$ (mod 2) let 
\begin{equation}\label{defHab}
H(x)=  H_{a, b}(x) = i^{b-a} J_a(4\pi x) (4\pi x)^{-b}.
\end{equation}
 By \cite[(2.21)]{BHM} we have 
 $$\mathscr{L}^{+}H(t) = \frac{b!}{2^b} \prod_{j=0}^b \left(t^2 + \Big(\frac{a+b}{2} - j\Big)^2\right)^{-1}, \quad \mathscr{L}^{\text{hol}}H(k) = \frac{b!}{2^b} \prod_{j=0}^b \left(\Big(\frac{(1-k)i}{2}\Big)^2 + \Big(\frac{a+b}{2} - j\Big)^2\right)^{-1}.$$
Obviously we have
\begin{equation}\label{pos}
\mathscr{L}^{+}H(t) >0, \quad \mathscr{L}^{+}H(t) \asymp (1 + |t|)^{-2b-2}, \quad t \in \Bbb{R} \cup [-i\vartheta, i \vartheta]
\end{equation}
and
\begin{equation*}
 \mathscr{L}^{\text{hol}}H(k) > 0 \,\,\, \text{for}\,\,\, 2 \leq k \leq a-b, \quad |\mathscr{L}^{\text{hol}}H(k)| \asymp k^{-2b-2}. 
 \end{equation*}
We choose a constant $c(a, b)$ such that 
\begin{equation}\label{hposhol}
h_{\text{pos}}^{\text{hol}} (k)  := \mathscr{L}^{\text{hol}}H_{a, b}(k) + \delta_{k > a-b} c(a, b) k^{-2b-1} > 0
\end{equation}
for all $k \in 2\Bbb{N}$, and we put $h_{\text{pos}} :=  \mathscr{L}^{+}H_{a, b}$. 
Then by  \eqref{pos} and \eqref{hposhol}, the pair $\mathfrak{h}_{\text{pos}} = (h_{\text{pos}},h^{\text{hol}}_{\text{pos}} )$ (which depends on $a, b$, but this is suppressed from the notation) is strictly positive on the entire spectrum $\mathcal{T}_0$, and we have by \eqref{inv}, \eqref{defHab} and \eqref{hposhol} the identity  
 \begin{equation}\label{pos-id}
   \mathscr{K}^{\ast} \mathfrak{h}_{\text{pos}}(x) = i^{b-a} J_a(4\pi x) (4\pi x)^{-b} + \sum_{a-b < k \in 2\Bbb{N}}  \frac{i^{-k}(k-1)}{\pi}\frac{c(a, b)}{k^{2b+1}} J_{k-1}(4\pi x).
   \end{equation}
 
For future reference we state the following lemma, for which we recall the notation \eqref{mellin2}.

\begin{lemma}\label{final-decay}  Let $0 \leq \vartheta \leq 7/64$ and  $A \geq 5$ an  integer.  Let $\Phi$ be a function that is holomorphic in   $-2\vartheta-\varepsilon < \Re u < A$ for some $\varepsilon > 0$ and satisfies $\Phi(u) \ll (1 + |u|)^{-A}$ in this region. \\
{\rm a)} We have  $$\mathscr{L}^{\pm}\widecheck{\Phi}(t) \ll_A (1  + |t|)^{-A}, \quad  t \in \Bbb{R} \cup [-i\vartheta, i\vartheta],\quad\quad \mathscr{L}^{\text{{\rm hol}}}\widecheck{\Phi}(k) \ll_A k^{-A}, \quad k \in 2\Bbb{N}.$$ 
{\rm b)}  If  $0 \leq \tau < 1$ and  $\Phi$ is in addition meromorphic in $\Re u \geq -2\tau - \varepsilon$ with finitely many poles at   $u_1, \ldots, u_n \in \Bbb{C}$, then $\mathscr{L}^{\pm}\widecheck{\Phi}(t)$ has meromorphic continuation to $|\Im t| < \tau$ with poles at most at $t = \pm i u_j/2$, $j = 1, \ldots, n$.  
\end{lemma}

\textbf{Proof.} a) Using the definitions \eqref{Hback} and \eqref{mellin2} and exchanging integrals, we have for $t \in \Bbb{R} \cup [-i\vartheta, i\vartheta]$ that 
\begin{equation}\label{start}
\mathscr{L}^{\pm}\widecheck{\Phi}(t) = \int_{(2\vartheta + \varepsilon/2)}\widehat{ \mathcal{J}^{\pm}(., t)}(u) \Phi(-u) \frac{du}{2\pi i}
\end{equation}
as an absolutely convergent integral (by Stirling's formula), where $ \widehat{\mathcal{J}^{\pm}(., t)}(u)$ is given in \eqref{mellin-j}. This contour is within the region of holomorphicity of $\Phi(-u)$, but to the right of the poles of $\widehat{\mathcal{J}^{\pm}(., t)}$ for $t \in \Bbb{R} \cup [-i\vartheta, i\vartheta]$. To deduce the required bound, we may assume without loss of generality that $|t| \geq 1$. We shift the contour to the left to $\Re u = -A + 1/2$.  On the way we pick up poles at $u = -2n \pm 2 it$, $n \in \{0, 1, \ldots, [\frac{1}{2}(A-\frac{1}{2})]\}$ with residues of the shape
$$\pm 2 (2\pi)^{2n \pm 2 it} \cosh(\pi t) \Gamma(-n \mp 2it) \Phi(2n \pm 2it) \ll (1 + |t|)^{-A-1/2}$$
with various sign combinations.  
We estimate the remaining integral by Stirling's formula \eqref{stir} as
$$\ll \int_{(-A+1/2)} \big((1 + |\Im u + 2t|)(1 + |\Im u - 2t|)\big)^{-(2A+1)/4} (1 + |u|)^{-A} |du| \ll (1 + |t|)^{-A-1/2},$$
which completes the proof  for $\diamondsuit \in \{+, -\}$. 

The holomorphic case is similar (and a little easier). Here we have nothing to show for $k \leq A$, and for $k \geq A$, the contour shift to the left $\Re u = -A+1 \in \Bbb{N}$ does not produce any poles. The remaining integral   can be estimated in the same way using the stronger bound $$ \widehat{ \mathcal{J}^{\text{hol}}(., k)}(-A+1 + iw) \ll  \frac{\Gamma((-A + k   + i w)/2)}{\Gamma(  (k + A- iw)/2) } \ll  (k + |w|)^{-A}$$ for $A \in \Bbb{N}$ by  the repeated application of the recursion  formula for the gamma function.

b) We fix some $t_0 \in \Bbb{C}$ with $|\Im t_0| < \tau$. By assumption we can bend the contour in \eqref{start} to be to the right of $2t_0$ (within the region of meromorphicity of $\Phi$), but the unbounded part still runs at real part equal\footnote{to ensure absolute convergence since a priori we have no growth condition for $|\Im t| > 2\vartheta + \varepsilon$ available.} to $2\vartheta + \varepsilon/2$. We may pick up poles with residue $\widehat{\mathcal{J}^{\pm}(., t)}(-u_j)$, which are meromorphic functions with poles at most at $t = \pm i u_j/2$, while the remaining integral is holomorphic in a neighbourhood of $t_0$.

\section{The ${\rm GL}(3)$ Voronoi summation formula}\label{voronoi-formula}

Let $F$ be a   cuspidal  automorphic form for the group ${\rm SL}_3(\Bbb{Z})$. As in the introduction, let $\theta \leq 5/14$ be an admissible bound for the Ramanujan-Petersson conjecture for $F$.  We denote its archimedean Langlands parameters by $\mu = (\mu_1, \mu_2, \mu_3)$ with 
\begin{equation}\label{mu}
  \mu_1 + \mu_2 + \mu_3 = 0,\quad |\Re \mu_j| \leq \theta, 
\end{equation}  
   and we denote its  Fourier coefficients by $A(n, m)$ as in \cite{Go}. They satisfy $A(n, m) = \overline{A(m, n)}$, 
\begin{equation}\label{hecke}
A(n, m) = \sum_{d \mid (n, m)} \mu(d) A(n/d, 1) A(1, m/d)
\end{equation}
(as follows from   \cite[Theorem  6.4.11]{Go} and M\"obius inversion), and the Rankin-Selberg bound
\begin{equation}\label{gl3-RS}
\sum_{n m^2 \leq x} |A(n, m)|^2 \ll x.
\end{equation}
Individually, we only know 
\begin{equation}\label{indiv}
   A(n, m) \ll (nm)^{\theta+\varepsilon}.
 \end{equation}  
We will always regard $F$ as fixed, and all implied constants may depend on $F$, in particular on $\mu$, but we suppress this from the notation.   

The Voronoi summation formula \cite[Theorem 1.8]{MS} states the following: for   $c, d, m \in \Bbb{Z}$ with $(c, d) = 1$, $c, m > 0$  and a smooth compactly supported weight function $w : \Bbb{R}_{>0} \rightarrow \Bbb{C}$ we have
\begin{equation}\label{Voronoi-ms}
\sum_{n} A(m, n) e\left(\frac{n \bar{d}}{c}\right) w(n) = c\sum_{\pm}  \sum_{n_1 \mid cm} \sum_{n_2 > 0} \frac{A(n_2, n_1)}{n_2n_1} S(   md, \pm n_2, mc/n_1)W^{\pm} \left(\frac{ n_2n_1^2}{c^3m}\right),
\end{equation}
where
$$W^{\pm}(x) =  \int_{(1)}   x^{-s}\mathcal{G}_{\mu}^{\pm} (s+1)  \widehat{w}(-s) \frac{ds}{2\pi i}  $$
with 
\begin{equation}\label{defgmu}
\begin{split}
\mathcal{G}_{\mu}^{\pm}(s) & = \frac{1}{2} \prod_{j=1}^3G_0(s+\mu_j) \pm \frac{1}{2i} \prod_{j=1}^3G_1(s+\mu_j)\\
& = 4(2\pi)^{-3s} \prod_{j=1}^3 \Gamma(s + \mu_j)\Bigg(\prod_{j=1}^3 \cos\Big(\frac{\pi(s + \mu_j)}{2}\Big) \pm \frac{1}{i}  \prod_{j=1}^3\sin\Big(\frac{\pi(s + \mu_j)}{2}\Big)\Bigg) .
\end{split}
\end{equation}
(To compare this with the formula in \cite{MS}, write $w = w_{\text{even}} + w_{\text{odd}}$ with $w_{\text{even}} (x)= \frac{1}{2}( w(x) + w(-x))$, $w_{\text{odd}} (x)= \frac{1}{2}( w(x) - w(-x))$, and observe that $M_{0}w_{\text{even}} = M_{1}w_{\text{odd}} = \widehat{w}$.)  The following lemma summarizes the analytic properties of $\mathcal{G}_{\mu}^{\pm}$. 

\begin{lemma}\label{lemG} The functions   $\mathcal{G}_{\mu}^{\pm}$ are meromorphic on $\Bbb{C}$ with poles at $s= -n - \mu_k$, $n \in \Bbb{N}_0$, $k \in \{1, 2, 3\}$. In particular, they are holomorphic in $\Re s > \theta$. Away from poles they satisfy the bound\footnote{Recall that all implied constants may depend on $\mu$.}
\begin{equation}\label{bound-G}
\mathcal{G}_{\mu}^{\pm}(s) \ll_{\Re s} (1 + |\Im s|)^{3 \Re s - \frac{3}{2}} e^{-\pi \max(0, \pm \Im s/2)}.
\end{equation}
In particular,  $\mathcal{G}_{\mu}^{\pm}(s)$ is exponentially decaying for $\pm \Im s \rightarrow \infty$. 
For $| t| \geq 3$ sufficiently large   we have the asymptotic formula
\begin{equation}\label{asymp-G}
\begin{split}
\mathcal{G}_{\mu}^{\pm}(\sigma + i t) = |t|^{3\sigma - \frac{3}{2}}
& \exp\left(  3i t\log \frac{|t  |}{2\pi e}\right) w_{\sigma, M, \mu}(t) + O_{\sigma, M}\left(|t|^{-M}\right)
\end{split}
\end{equation}
with
$$|t|^j w_{\sigma, M, \mu}^{(j)}(t) \ll_{j, \sigma, M}  1$$
for all $j, M \in \Bbb{N}_0$. Finally, for $n \in \Bbb{N}_0$ and $\epsilon_2 \in \{\pm 1\}$ we have 
\begin{equation}\label{new-formula}
 \sum_{\epsilon_1 \in \{\pm 1\}} e^{\epsilon_1 i\pi n/2} \mathcal{G}_{-\mu}^{-\epsilon_1}(1-v + n) \mathcal{G}_{\mu}^{\epsilon_1\epsilon_2}(v)  = \frac{1+\epsilon_2}{2} (2\pi)^{-3 n}  \prod_{k=1}^3   \prod_{j=1}^n (v-j+\mu_k).  
\end{equation}
\end{lemma}

\textbf{Proof.} Equations \eqref{bound-G} and \eqref{asymp-G} follow directly from  from \eqref{mu} and Stirling's formula \eqref{stir}. For the proof of \eqref{new-formula} we use the recursion and reflection formula of the gamma function to see that $\mathcal{G}_{-\mu}^{-\epsilon_1}(1-v + n) \mathcal{G}_{\mu}^{\epsilon_1\epsilon_2}(v) $ equals 
\begin{displaymath}
\begin{split}
    16(2\pi)^{-3-3n} &\Big[  \prod_{k=1}^3   \frac{\pi \prod_{j=1}^n (-v+j-\mu_k) }{  \sin(\pi(v+\mu_k))} \Big] 
 \Big[ \prod_{k=1}^3 \cos\left(\frac{\pi(v+\mu_k)}{2}\right) +\frac{\epsilon_1\epsilon_2}{i} \prod_{k=1}^3 \sin\left(\frac{\pi(v+\mu_k)}{2}\right)\Big] \\ 
& \Big[\prod_{k=1}^3 \cos\left(\frac{\pi(1-v+n-\mu_k)}{2}\right) - \frac{\epsilon_1}{i} \prod_{k=1}^3 \sin\left(\frac{\pi(1-v+n-\mu_k)}{2}\right)\Big].
\end{split}
\end{displaymath}
Hence $e^{\epsilon_1 i \pi n/2}\mathcal{G}_{-\mu}^{-\epsilon_1}(1-v + n) \mathcal{G}_{\mu}^{\epsilon_1\epsilon_2}(v) $ equals 
\begin{displaymath}
\begin{split}
    2(2\pi)^{ -3n} &\Big[  \prod_{k=1}^3   \frac{ \prod_{j=1}^n (-v+j-\mu_k) }{  \sin(\pi(v+\mu_k))} \Big] 
 \Big[ \prod_{k=1}^3 \cos\left(\frac{\pi(v+\mu_k)}{2}\right) +\frac{\epsilon_1\epsilon_2}{i} \prod_{k=1}^3 \sin\left(\frac{\pi(v+\mu_k)}{2}\right)\Big] \\ 
& (-1)^n \Big[\epsilon_1 i \prod_{k=1}^3 \cos\left(\frac{\pi( v+\mu_k)}{2}\right) +   \prod_{k=1}^3 \sin\left(\frac{\pi(v+\mu_k)}{2}\right)\Big].
\end{split}
\end{displaymath}
Summing over $\epsilon_1 \in \{\pm 1\}$, we can drop all terms depending linearly on $\epsilon_1$, and  \eqref{new-formula} follows easily from the addition theorem of the $\sin$-function. 
 \\

We rewrite the summation formula in terms of Dirichlet series as follows. For positive   integers $c, d, m$ with $(c, d) = 1$ let
\begin{equation}\label{defPhi}
\Phi(c, \pm d, m; v) := \sum_{n > 0} A(m, n) e\left(\pm \frac{n \bar{d}}{c}\right) n^{-v}.
\end{equation}
By \eqref{hecke} -- \eqref{indiv} this is absolutely convergent in $\Re v > 1$ and satisfies the uniform bound 
\begin{equation}\label{phibound}
  \Phi(c, \pm d, m; v) \ll \alpha(m) :=  m^{\varepsilon}\max_{d \mid m} |A(d, 1)|, \quad \Re v \geq 1+\varepsilon.
  \end{equation} Moreover, let
 \begin{equation}\label{defXi}
 \Xi(c, \pm d, m; v) := c \sum_{n_1 \mid cm} \sum_{n_2 > 0} \frac{A(n_2, n_1)}{n_2n_1} S(\pm  md, n_2, mc/n_1) \left(\frac{ n_2n_1^2}{c^3m}\right)^{-v}.
 \end{equation}
This is absolutely convergent in $\Re v > 0$ and satisfies the uniform bound
\begin{equation}\label{xibound}
\Xi(c, \pm d, m; v) \ll (mc^3)^{1/2 + \Re v+\varepsilon} ,  \quad \Re v \geq \varepsilon,
\end{equation}
using Weil's bound for Kloosterman sums and again \eqref{hecke} -- \eqref{indiv}. 
  
By converting additive characters into multiplicative characters we see that the functions $\Phi(c, \pm d, m; .) $ and $\Xi(c, \pm d, m; .)$ are, up to finitely many Euler factors at primes dividing $m$ that are holomorphic in  $\Re v > \theta$ resp.\ $\Re v \geq \theta-1$, linear combinations of $L$-functions corresponding to $F \times \chi$ for Dirichlet characters $\chi$. This shows that  $\Phi(c, \pm d, m; .)$ is of finite order and analytic in (an $\varepsilon$-neighbouhood of) $\Re v \geq 1/2$ and $\Xi(c, \pm d, m; .)$ is of finite order and analytic in $\Re v \geq -1/2$. 
 The Voronoi formula is equivalent to the  vector-valued functional equation
\begin{equation}\label{func1}
\left(\begin{matrix}\Phi(c, d, m; v) \\ \Phi(c, -d, m; v)\end{matrix}\right) =  \left(\begin{matrix} \mathcal{G}_{\mu}^+(1-v) & \mathcal{G}_{\mu}^-(1-v) \\ \mathcal{G}_{\mu}^-(1-v) & \mathcal{G}_{\mu}^+(1-v) \end{matrix}\right)  \left(\begin{matrix} \Xi(c, d, m; -v) \\ \Xi(c, -d, m; -v)\end{matrix}\right).
\end{equation}
Inverting the ``scattering matrix'' (using \eqref{new-formula} with $n=0$), we obtain
\begin{equation}\label{func2}
\left(\begin{matrix}\Xi(c, d, m; v) \\ \Xi(c, -d, m; v)\end{matrix}\right) =  \left(\begin{matrix} \mathcal{G}_{-\mu}^-(-v) & \mathcal{G}_{-\mu}^+(-v) \\ \mathcal{G}_{-\mu}^+(-v) & \mathcal{G}_{-\mu}^-(-v) \end{matrix}\right)  \left(\begin{matrix} \Phi(c, d, m; -v) \\ \Phi(c, -d, m; -v)\end{matrix}\right).
\end{equation}
In particular, both $\Xi(c, \pm d, m; .) $ and $\Phi(c, \pm d, m; .)$ are of finite order and entire. 
 By the Phragm\'en-Lindel\"of principle, \eqref{phibound}, \eqref{xibound} and \eqref{bound-G} we obtain 
\begin{equation}\label{xi}
\begin{split}
& \Phi(c, d, m; v) \ll_{\Re v} \alpha(m)\left(mc^3 (1 + |\Im v|)^3\right)^{ \max\left(\frac{1}{2} - \Re v, \frac{1}{2} (1 - \Re v),  0\right)+\varepsilon},\\
&\Xi(c, d, m; v) \ll_{\Re v} \alpha(m) \left(mc^3\right)^{ \max\left(\frac{1}{2} + \Re v, \frac{1}{2} (1 + \Re v),  0\right)+\varepsilon} (1 + |\Im v|^3)^{\max(0, -\Re v-\frac{1}{2}, - \frac{1}{2}\Re v)+\varepsilon}.
\end{split}
\end{equation}

The Voronoi formula is a consequence of the functional equation of the twisted $L$-functions $L(s, F \times \chi)$ and relations in the unramified ${\rm GL}(3)$ Hecke algebra. Therefore the functional equations \eqref{func1} and \eqref{func2} must continue to hold if $F$ is an Eisenstein series, except that in this case the functions $\Xi(c, \pm d, m; .) $ and $\Phi(c, \pm d, m; .)$ are not necessarily entire any more (since $L(s, F \times \chi)$ will have a pole if $F$ is Eisenstein and $\chi$ is trivial). We are only interested in the case when $F= E_0$ is the minimal parabolic Eisenstein series with trivial parameters $\mu= (0, 0, 0)$. In this case, $A(n, 1) = A(1, n) = \tau_3(n)$ is the ternary divisor function, and in general $A(n, m)$ is given by \eqref{hecke}. As above, we see that $\Phi(c, \pm d, m; .)$ can only have a (triple) pole at $v=1$ in $\Re v \geq 1/2$ and $\Xi(c, \pm d, m; .) $ can only have a (triple) pole at $v=0$ in $\Re v \geq -1/2$.  (A classical Voronoi formula for $\E_0$ analogous to \eqref{Voronoi-ms} with extra polar terms has been worked out  in \cite{Li2}.)  The corresponding Laurent expansions are computed in Lemma \ref{lem6} and \ref{pole1} below.  

\section{A multiple Dirichlet series}\label{multiple}

Fix   $\ell, q \in \Bbb{N}$ and   $s, w, u \in \Bbb{C}$. If 
   \begin{equation}\label{variables}
   \begin{split}
 &\Re (w+u/2)> 1, 
   \quad  \Re(3s + u/2) > 2, \quad \Re (3s - u/2) > 4, \quad  \Re u < -1/2,  \end{split}
 \end{equation}
 we define 
 \begin{equation}\label{defD}
 \mathcal{D}^{\pm}_{q, \ell}(s, u, w) := \sum_{(c, d) = 1} \sum_{\substack{q \mid mc\\ \ell \mid md}}  \frac{  \Xi(c, \pm d, m; -1 + s + u/2)}{c^{3s+u/2 -1} m^{s+w}d^{w+u/2}}
 \end{equation}
 and 
 \begin{equation}\label{deftildeD}
\widetilde{\mathcal{D}}^{\pm}_{q, \ell}(s, u, w) := \sum_{(c, d) = 1} \sum_{\substack{q \mid mc\\ \ell \mid md}} \frac{\Phi(c, \pm d, m; 1 - s - u/2)}{c^{3s + u/2 - 1}m^{w+s} d^{w + u/2}}
\end{equation} 
 with $\Phi$ as in  \eqref{defPhi} and  $\Xi$ as in \eqref{defXi}, satisfying the functional equations
 \begin{equation}\label{Dfunc}
 \begin{split}
  &  \mathcal{D}^{\pm}_{q, \ell}(s, u, w) = \mathcal{G}^{\mp}_{-\mu}(1-s-u/2) \widetilde{\mathcal{D}}^{+}_{q, \ell}(s, u, w) +  \mathcal{G}^{\pm}_{-\mu}(1-s-u/2) \widetilde{\mathcal{D}}^{-}_{q, \ell}(s, u, w),\\
& \widetilde{\mathcal{D}}^{\pm}_{q, \ell}(s, u, w) = \mathcal{G}^{\pm}_{\mu}(s+u/2) \mathcal{D}^{+}_{q, \ell}(s, u, w) +  \mathcal{G}^{\mp}_{\mu}(s+u/2) \mathcal{D}^{-}_{q, \ell}(s, u, w)
 \end{split}
 \end{equation} 
 by \eqref{func2} and \eqref{func1}.   Recall that the numerator of \eqref{defD} and \eqref{deftildeD} is holomorphic in the region \eqref{variables} if $F$ is cuspidal. 
 
 If $F = \E_0$, the only polar divisor of $\mathcal{D}^{\pm}_{q, \ell}(s, u, w)$ in $\Re(-1+s+u/2) \geq -1/2$ can occur at $-1+s+u/2 = 0$, and the only polar divisor of $\tilde{\mathcal{D}}^{\pm}_{q, \ell}(s, u, w)$ in $\Re(1-s-u/2) \geq 1/2$ can be at $1-s-u/2 = 1$.  We will see in a moment in Lemmas \ref{pole1} and \ref{pole2}   that the Laurent expansions at these poles are independent of the $\pm$ sign. Since $\mathcal{G}_{\mu}^+(s) + \mathcal{G}^-_{\mu}(s)$ has a triple zero at $s=1$, the functional equations \eqref{Dfunc} imply that no other polar divisors of $\mathcal{D}^{\pm}_{q, \ell}(s, u, w)$ and $\widetilde{\mathcal{D}}^{\pm}_{q, \ell}(s, u, w)$ can exist in the domain \eqref{variables} of definition.

By \eqref{xi} the triple sums in \eqref{defD} and \eqref{deftildeD} are absolutely convergent: the first condition in \eqref{variables} ensures absolute convergence of the $d$-sum, the other three conditions ensure absolute convergence of the $c$-sum; the $m$-sum requires $\Re(s+w) > 1$, $\Re(w - \frac{1}{2}u) > 1/2$ and $\Re (\frac{1}{2}s + w - \frac{1}{4}u) > 1$ for absolute convergence, which also follows from \eqref{variables}.  Again by \eqref{xi} we have  the bounds
  \begin{equation}\label{boundD}
  \begin{split}
 &  \mathcal{D}^{\pm}_{q, \ell}(s, u, w) \ll_{s, w} (1 + |u|)^{3\max(0, \frac{1}{2} - \Re (s + \frac{1}{2}u),  \frac{1}{2} - \Re( \frac{1}{2}s + \frac{1}{4}u))+\varepsilon},\\
   & \widetilde{\mathcal{D}}^{\pm}_{q, \ell}(s, u, w) \ll_{s, w} (1 + |u|)^{3\max(0, -\frac{1}{2} + \Re (s + \frac{1}{2}u),    \Re( \frac{1}{2}s + \frac{1}{4}u))+\varepsilon}
   \end{split}
   \end{equation}
  for $s, u, w$ satisfying \eqref{variables} (away from the pole if $F = \E_0$).  
 If in addition $\Re(s + u/2) > 1$, we can insert the definition \eqref{defXi} into \eqref{defD}. This gives the following alternative representation.
 
 \begin{lemma}\label{lem6} Let   $s, u, w \in \Bbb{C}$ satisfy 
 \begin{equation}\label{lem1-cond}
 \Re(w + u/2) > 1, \quad \Re(s + u/2) > 1, \quad \Re u < -1/2.
 \end{equation}
 Then 
 \begin{equation}\label{alt}
 \mathcal{D}^{\pm}_{q, \ell}(s, u, w) = \sum_{\ell \mid r} \sum_{q \mid n_1c} \sum_{n_2} \frac{A(n_2, n_1)S(\pm r, n_2, c)}{n_2^{s + u/2 }n_1^{2s} c^{1-u}r^{w+u/2}}.
 \end{equation}
 \end{lemma}

\textbf{Proof.} The conditions \eqref{lem1-cond} imply \eqref{variables}. By absolute convergence we can re-arrange the sums to see that $\mathcal{D}^{\pm}_{q, \ell}(s, u, w)$ equals
 \begin{displaymath}
\begin{split}
 &  \sum_{  d} \sum_{\ell \mid md}  \frac{1}{(md)^{w+u/2}} \sum_{n_1} \sum_{\substack{q \mid mc\\ (c, d) = 1\\ \frac{n_1}{(n_1, m)} \mid c}} \frac{1}{(mc)^{1-u}}\sum_{n_2} \frac{A(n_2, n_1)}{n_2^{s + u/2 } n_1^{2s + u-1}} S\left(\pm md, n_2, \frac{mc}{n_1}\right)  \\
&= \sum_{  d} \sum_{\ell \mid md} \frac{1}{(md)^{w+u/2}}  \sum_{fg = m} \sum_{(\nu_1, g) = 1} \sum_{\substack{q \mid mc\\ (c, d) = 1\\ \nu_1 \mid c}}\frac{1}{(mc)^{1-u}} \sum_{n_2} \frac{A(n_2, \nu_1f)}{n_2^{s + u/2 } (\nu_1f)^{2s + u-1}} S\left(\pm md, n_2, \frac{gc}{\nu_1}\right) \\
& = \sum_{  d}   \sum_{ \ell \mid fgd}\frac{1}{(fgd)^{w+u/2}}   \sum_{(\nu_1, g) = 1} \sum_{\substack{ (\nu_1 \gamma, d) = 1\\ q \mid fg\nu_1\gamma}} \frac{1}{(fg\nu_1 \gamma)^{1-u}}\sum_{n_2} \frac{A(n_2, \nu_1f)S(\pm fgd, n_2, g\gamma ) }{n_2^{s + u/2 } (\nu_1f)^{2s + u-1}} .   
   \end{split}
\end{displaymath}
There is a bijection between  integer quintuples 
$$(d, f, g, \nu_1, \gamma) \quad  \text{satisfying} \quad   \ell \mid fgd, \,q \mid fg\nu_1\gamma, \,(\nu_1, g) = 1 , \,(\nu_1\gamma, d) = 1$$
and  integer triples
$$(n_1, r, c) \quad \text{satisfying} \quad q\mid n_1 c, \, \ell \mid r$$
given by 
$$(n_1, r, c) = (\nu_1f, f gd , g\gamma)$$
with inverse map
$$(d, f, g, \nu_1, \gamma) = \left(\frac{r}{(n_1c, r)}, (n_1, r), \frac{(n_1c, r)}{(n_1, r)}, \frac{n_1}{(n_1, r)}, \frac{c(n_1, r)}{(n_1c, r)} \right). $$
This shows the desired formula \eqref{alt}.\\


Next we compute the Laurent expansion of  $\mathcal{D}^{\pm}_{q, \ell}(s, u, w)$ at $u = 2-2s$ if $F = \E_0$. We write
\begin{equation}\label{resR}
\mathcal{D}_{q, \ell}^{\pm}(s, u, w) = \sum_{j=1}^3 \frac{\mathcal{R}_{q,\ell;j}(s, w)}{(u - (2-2s))^j} + O(1).
\end{equation}
We will see in the proof of the next lemma that    $\mathcal{R}_{q,\ell; j}(s, w)$ is independent of the $\pm$ sign. 

\begin{lemma}\label{pole1} Let $F = \E_0$ and suppose $(q, \ell) = 1$. Then the Laurent coefficients $\mathcal{R}_{q, \ell;j}(s, w)$ are  meromorphic in $\Bbb{C}\times \Bbb{C}$. In an $\varepsilon$-neighbourhood  of the region $\Re w \geq \Re s\geq 1/2$, the function
\begin{equation}\label{hol}
\big((s - \textstyle\frac{1}{2})(w-s)(s+w-1)\big)^4\mathcal{R}_{q, \ell;j}(s, w)
\end{equation}
 is holomorphic and  bounded by
$O_{s, w} \left((q\ell)^{-1+\varepsilon}\right)$. 
\end{lemma}

\textbf{Proof.} In the region \eqref{lem1-cond} we have (recall \eqref{hecke} and $A(n, 1) = A(1, n) = \tau_3(n)$)
\begin{equation}\label{dql}
\begin{split}
\mathcal{D}_{q, \ell}^{\pm}(s, u, w) & = 
\sum_{\ell \mid r}  \sum_{q \mid n_1 c} \sum_{d,    b_1, b_2, b_3} \frac{\tau_3(n_1) \mu(d)   S(\pm r, db_1b_2b_3, c)}{n_1^{2s} d^{3s + u/2}(b_1 b_2 b_3)^{s + u/2} c^{1-u} r^{w+u/2}}\\ 
& =  
 \sum_{\ell \mid r} \sum_{q \mid n_1c} \sum_{d} \sum_{ b_1, b_2, b_3\, (\text{mod }c)} \frac{\tau_3(n_1)\mu(d)   S(\pm r, db_1b_2b_3, c)}{n_1^{2s} d^{3s + u/2}  c^{1+3s+u/2} r^{w+u/2}} 
  \prod_{j=1}^3\zeta\left(s + u/2, \frac{b_j}{c}\right)
\end{split}
\end{equation} 
where
\begin{equation}\label{hurwitz}
\zeta(s, \alpha) = \sum_{(n+\alpha) > 0} \frac{1}{(n+\alpha)^s} = \frac{1}{s-1} - \psi(\alpha) - \gamma(\alpha)(s-1) + \ldots
\end{equation}
is the Hurwitz zeta function, and   $\psi$, $\gamma$ are suitable functions denoting the Taylor coefficients. We note that 
\begin{equation}\label{averagezeta}
\begin{split}
\sum_{b \, (\text{mod } m)} \zeta(s, b/m)& = m^s \zeta(s)\\
& = \frac{m}{s-1} + m(\log m + \gamma) + \left(\frac{1}{2}m \log^2 m + \gamma m \log m - m \gamma_1\right)(s-1) + \ldots 
\end{split}
\end{equation}
for Euler's constant $\gamma = 0.577\ldots$  and another constant $\gamma_1 \in \Bbb{R}$. 

Inserting \eqref{hurwitz} into \eqref{dql}, we compute   the Laurent coefficients as 
 $$\mathcal{R}_{q ,\ell;j}(s, w) = 
 \sum_{\ell \mid r} \sum_{q \mid n_1c}  \sum_{d} \sum_{ b_1, b_2, b_3\, (\text{mod }c)} \frac{\tau_3(n_1) \mu(d)   S(\pm r, db_1b_2b_3, c)}{n_1^{2s} d^{2s + 1}  c^{2+2s} r^{w+1-s}}
 \rho_j\left(\frac{b_1}{c}, \frac{b_2}{c}, dcr\right)$$
 where
 \begin{equation}\label{rho}
 \begin{split}
  & \rho_1(n_1, n_2, m) = 8, \quad \rho_2(n_1, n_2, m) = -12\psi(n_1) - 4\log m, \\
   & \rho_3(n_1, n_2, m) = 6\psi(n_1)\psi(n_2) - 3\gamma(n_1) + (\log m)^2 + 6 \psi(n_1) \log m.
 \end{split}
 \end{equation}
  We open the Kloosterman sum and evaluate the $b_3$-sum getting
\begin{displaymath}
\begin{split}
 \mathcal{R}_{q,\ell; j}(s, w)& =  
 \sum_{\ell \mid r}  \sum_{q \mid n_1c} \sum_{d} \sum_{\substack{ b_1, b_2\, (\text{mod }c)\\ c\mid db_1b_2}} \frac{\tau_3(n_1) \mu(d)   S(\pm r, 0, c)}{n_1^{2s}d^{2s + 1}  c^{1+2s} r^{w+1-s}}
  \rho_j\left(\frac{b_1}{c}, \frac{b_2}{c}, dcr\right)\\
&  =  
 \sum_{q \mid n_1c} \sum_{d} \sum_{a \mid c}  \sum_{\ell \mid ra}\frac{\mu(c/a)}{a^{w-s}} \sum_{\substack{ b_1, b_2\, (\text{mod }c)\\ c\mid db_1b_2}} \frac{\tau_3(n_1)  \mu(d)
  }{n_1^{2s}d^{2s + 1}  c^{1+2s} r^{w+1-s}}\rho_j\left(\frac{b_1}{c}, \frac{b_2}{c}, dca r\right).
\end{split}
\end{displaymath}
At this point it is clear that the $\pm$-sign plays no role.  We order this by $\beta_i = (b_i, c)$ getting
\begin{equation}\label{1new}
\mathcal{R}_{q,\ell; j}(s, w) = 
  \sum_{q \mid n_1c}\sum_{d} \sum_{a \mid c}  \sum_{\ell \mid ra}\sum_{\substack{\beta_1, \beta_2 \mid c\\ c \mid d \beta_1\beta_2}} \frac{\tau_3(n_1) \mu(d)   \mu(c/a) 
  }{n_1^{2s}a^{w-s} d^{2s + 1}  c^{1+2s} r^{w+1-s}}R_j\left(\frac{c}{\beta_1}, \frac{c}{\beta_2}, dca r\right)
\end{equation}
where
\begin{equation*}
\begin{split}
& R_3(n_1, n_2, m) = 8\phi(n_1)\phi(n_2), \\
& R_2(n_1, n_2, m) = -12 \Psi^{\ast}(n_1)\phi(n_2) - 4 \phi(n_1)\phi(n_2) \log m, \\
& R_1(n_1, n_2, m) = 6\Psi^{\ast}(n_1) \Psi^{\ast}(n_2) - 3 G^{\ast}(n_1)\phi(n_2) + \phi(n_1)\phi(n_2) (\log m)^2 + 6\Psi^{\ast}(n_1)\phi(n_2) \log m 
\end{split}
\end{equation*}
with
$$\Psi^{\ast}(n) = \sum_{\substack{b \, (\text{mod } n)\\ (b, n) = 1}}  \psi(b/n), \quad G^{\ast}(n) =    \sum_{\substack{b \, (\text{mod } n)\\ (b, n) = 1}}  \gamma(b/n).$$
Let 
\begin{equation}\label{open}
\phi_j(n) =  \sum_{ab= n} \mu(a)b (\log b)^j
\end{equation}
denote the Dirichlet series coefficient of $(-1)^j \zeta^{(j)}(s-1)/\zeta(s)$. In particular, $\phi_0  = \phi$ is the Euler phi-function. 
Then
$$\Psi(m) :=  \sum_{b\, (\text{mod }m)}  \psi(b/m) = -m(\gamma + \log m)$$
by \eqref{averagezeta}, so that 
by M\"obius inversion
$$\Psi^{\ast}(n) = \sum_{ab = n} \mu(a) \Psi(b) = -\gamma \phi_0(n)- \phi_1(n).  $$
Similarly, we have 
$$G(m) :=  \sum_{b\, (\text{mod }m)}  \gamma(b/m) =  - \frac{1}{2} m \log^2 m-\gamma m \log m + m \gamma_1,$$
so that 
$$G^{\ast}(n) = \sum_{ab = n} \mu(a) G(b) = \gamma_1 \phi_0(n) - \gamma \phi_1(n) - \frac{1}{2} \phi_2(n) .$$
Altogether, 
$$R_j(n_1, n_2, m) =  \sum_{\nu+\mu+\kappa \leq 3-j} C_{\nu, \mu, \kappa, j}\phi_{\nu}(n_1) \phi_{\mu}(n_2) (\log m)^{\kappa}$$
for certain constants $C_{\nu, \mu, \kappa, j}$.  Substituting back into \eqref{1new} we obtain a linear combination of 
 \begin{equation*} 
\sum_{q \mid n_1c}  \sum_{d}  \sum_{a \mid c} \sum_{\ell \mid r a} \sum_{\substack{\beta_1, \beta_2 \mid c\\ c \mid d \beta_1\beta_2}} \frac{\tau_3(n_1)\mu(d)   \mu(c/a)
 }{n_1^{2s}a^{w-s} d^{2s + 1}  c^{1+2s} r^{w+1-s}}\phi_{\nu}(c/\beta_1) \phi_{\mu}(c/\beta_2) (\log dca r)^{\kappa} 
 \end{equation*}
with $\nu + \mu + \kappa \leq 3-j$. 
We make several changes of variables. We write $\beta_i  \gamma_i = c$, switch to the co-divisor $\gamma_i$ and  open the $\phi_{\nu}$-functions as in \eqref{open}, writing $\gamma_i = a_ib_i$. Next we replace $c$ with $ac$ and recast the previous display as
$$  
\sum_{q \mid n_1a c}  \sum_{d}    \sum_{\substack{a_1b_1, a_2b_2 \mid a c\\ a_1b_1a_2b_2 \mid acd  }}  \sum_{\ell \mid  ra} \frac{\tau_3(n_1)\mu(d)   \mu(c) \mu(a_1)\mu(a_2) b_1b_2 
}{n_1^{2s}(c  d)^{2s + 1}  a^{1+s+w} r^{w+1-s}}
(\log b_1)^{\nu} (\log b_2)^{\mu} (\log da^2c r)^{\kappa}.
$$

We introduce the following generalized function
\begin{equation}\label{defZ}
Z_{q, \ell}(x_1, x_2, x_3) := \sum_{q \mid n_1a c}  \sum_{d}  \sum_{\substack{a_1b_1, a_2b_2 \mid a c\\ a_1b_1a_2b_2 \mid acd  }}  \sum_{\ell \mid  ra} \frac{\tau_3(n_1)\mu(d)   \mu(c) \mu(a_1)\mu(a_2) b_1b_2 
}{n_1^{2s}(c d)^{2s + 1}  a^{1+s+w} r^{w+1-s}}
b_1^{x_1}  b_2^{x_2} (da^2c  r)^{x_3}
\end{equation}
so that $\mathcal{R}_{q, \ell; j}(s, w)$ is a linear combination of 
\begin{equation*}
\frac{\partial^{\nu}}{\partial x_1^{\nu}}\frac{\partial^{\mu}}{\partial x_2^{\mu}}\frac{\partial^{\kappa}}{\partial x_3^{\kappa}}Z_{q, \ell}(x_1, x_2, x_3)|_{x_1 = x_2 = x_3 = 0}
\end{equation*}
with $\nu + \mu + \kappa = 3 - j \leq 2$. 
The function $Z_{q, \ell}(x_1, x_2, x_3)$ can be written as an Euler product, and for a generic prime $p \nmid \ell q $ the $p$-Euler factor equals
$$(1 - p^{-2s})^{-3}(1 - p^{-(1+w-s-x_3)})^{-1} \sum_{ \alpha  = 0}^{\infty} \sum_{\alpha_1, \alpha_2, \gamma, \delta = 0}^1\sum_{\substack{\alpha_1+\beta_1, \alpha_2+\beta_2\leq \alpha+\gamma\\ \alpha_1 + \beta_1 + \alpha_2 + \beta_2 \leq \alpha+ \gamma + \delta}} \frac{(-1)^{\delta+\gamma + \alpha_1 + \alpha_2}p^{ \beta_1(1+x_1) + \beta_2(1+x_2)}}{p^{(\gamma + \delta)(2s+1-x_3) + \alpha(1 + s + w - 2x_3) }} .$$
This can be expressed in closed form using geometric series,  for instance by  distinguishing the 5 disjoint cases (i) $\alpha = \gamma = 0$, (ii) $\alpha + \gamma \geq 1$, $\alpha_1 = \alpha_2 = \beta_1 = \beta_2 = 0$, (iii)  $\alpha + \gamma \geq 1$, $\alpha_1  = \beta_1 = 0$, $\alpha_2 +  \beta_2  \geq 1$, (iv) $\alpha + \gamma \geq 1$, $\alpha_1  + \beta_1 \geq 1$, $\alpha_2 +  \beta_2 =  0$ and (v) $\alpha + \gamma \geq 1$, $\alpha_1  + \beta_1 \geq 1$, $\alpha_2 +  \beta_2 \geq  1$. After a lengthy, but completely straightforward computation we obtain the beautiful expression
 \begin{equation}\label{euler}
 \frac{(1 - p^{-2s})^{-3}(1 - p^{-(1+w-s-x_3)})^{-1}(1-p^{-(s+w-x_1-2x_3)})^{-1}(1-p^{-(s+w-x_2-2x_3)})^{-1}}{(1 - p^{-(1+s+w-2x_3)})^{-1}(1 - p^{-(2s - x_1-x_3)})^{-1}(1 - p^{-(2s - x_2-x_3)})^{-1}}.
 \end{equation}
The Euler factor at primes   $p \mid q\ell$ can again be computed by a more complicated, but finite computation with geometric series, and it is clear that it is a rational function in $p^{-s}$, $p^{-w}$ and $p^{-x_j}$, $j = 1, 2, 3$, in particular meromorphic. 
Since the product of \eqref{euler} over all primes is a quotient of zeta-functions, it follows that $\mathcal{R}^{\pm,j}_{q, \ell}(s, w)$ is a meromorphic function. 

In an $\varepsilon$-neighbourhood of the region $\Re w \geq \Re s \geq 1/2$ and $x_1 = x_2 = x_3 = 0$,  the Euler factors at primes $p \mid \ell q$ converge absolutely. In particular, they are holomorphic, and we see after taking derivatives and putting $x_1 = x_2 = x_3 = 0$ that \eqref{hol} is holomorphic. To get the desired bound, we estimate the Euler factors at $p \mid   q$ trivially as
$$4\underset{\nu_1 + \alpha + \gamma \geq v_p(q)}{\sum_{ \nu_1, \alpha, \rho = 0}^{\infty} \sum_{\gamma, \delta  = 0}^1}  \sum_{\beta_1 + \beta_2 \leq \alpha + \delta + \gamma}  \frac{(\nu_1+2)(\nu_1+1)/2 }{p^{\nu_1 +2(\alpha + \delta+\gamma)+\rho - \beta_1 - \beta_2 + O(\varepsilon)}} \ll (p^{v_p(q)})^{-1+O(\varepsilon)}.$$
The same argument with the condition $\alpha + \rho \geq v_p(\ell)$ instead of $\nu_1 + \alpha + \gamma \geq v_p(q)$  applies for Euler factors at $p \mid \ell$. This completes the proof of the lemma. \\

We also need to study the   Laurent expansion of the companion function
\begin{equation}\label{tildeD}
\widetilde{\mathcal{D}}^{\pm}_{q, \ell}(s, u, w)  = \sum_{j=1}^3 \frac{\tilde{\mathcal{R}}_{q,\ell; j}(s, w)}{(u + 2s)^j} + O(1)
\end{equation}
at $s + u/2 = 0$ in the case when $F = \E_0$ (again the Laurent coefficients are independent of the $\pm$ sign). 
\begin{lemma}\label{pole2}  Let $F = \E_0$ and $(q, \ell) = 1$. Then the Laurent coefficients $\tilde{\mathcal{R}}_{q, \ell; j}(s, w)$ are  meromorphic in $\Bbb{C}\times \Bbb{C}$. In an $\varepsilon$-neighbourhood  of the region $\Re w \geq \Re s\geq 1/2$, the function
\begin{equation}\label{hol1}
\big((s - \textstyle\frac{1}{2})(w-s)(s+w-1)\big)^3\tilde{\mathcal{R}}_{q, \ell;j}(s, w)
\end{equation}
 is holomorphic and  bounded by
 $O_{s, w}\left((\ell q)^{\varepsilon} q^{-1}\right)$. 
\end{lemma}

\textbf{Proof.} The proof is similar, so we can be brief. We have
\begin{displaymath}
\begin{split}
\Phi(c, \pm d, m, v)&  = \sum_n  A(m, n) e\left(\pm  \frac{n \bar{d}}{c}\right) n^{-v} = \sum_{r_1r_2 = m} \frac{\mu(r_1) A(r_2, 1)}{r_1^v} \sum_n \frac{A(1, n) e(\pm n r_1 \bar{d}/c)}{n^v}\\
& =  \sum_{r_1r_2 = m} \frac{\mu(r_1) \tau_3(r_2)}{r_1^v}  \frac{1}{c^{3v}} \sum_{a_1, a_2, a_3 \, (\text{mod }c)} e\left( \pm \frac{a_1a_2a_3 r_1 \bar{d}}{c}\right) \prod_{j=1}^3\zeta\left(v, \frac{a_j}{c}\right),
\end{split}
\end{displaymath}
so that
$$\widetilde{\mathcal{D}}^{\pm}_{q, \ell}(s, u, w) =  \sum_{(c, d) = 1} \sum_{\substack{q \mid r_1r_2c\\ \ell \mid r_1r_2d}}  \frac{\tau_3(r_2) \mu(r_1) }{r_1^{1+w-u/2}r_2^{s+w}c^{2-u} d^{w + u/2} }\sum_{a_1, a_2, a_3 \, (\text{mod }c)} e\left( \pm \frac{a_1a_2a_3 r_1 \bar{d}}{c}\right) \prod_{j=1}^3\zeta\left(1-s-\frac{u}{2}, \frac{a_j}{c}\right)  .$$
We conclude
$$
\tilde{\mathcal{R}}_{q,\ell;j}(s, w) =   \sum_{(c, d) = 1} \sum_{\substack{q \mid r_1r_2c\\ \ell \mid r_1r_2d}}    \frac{\tau_3(r_2) \mu(r_1) }{r_1^{1+w+s}r_2^{s+w}c^{2+2s} d^{w -s} }\sum_{a_1, a_2, a_3 \, (\text{mod }c)} e\left( \pm \frac{a_1a_2a_3 r_1 \bar{d}}{c}\right)\tilde{\rho}_j\left(\frac{a_1}{c}, \frac{a_2}{c}, \frac{d}{r_1c^{2}}\right)$$
with $\tilde{\rho}_j(n_1, n_2, m) = (-1)^j \rho_j(n_1, n_2, m)$ as in \eqref{rho}. We evaluate the $a_3$-sum getting
$$
\tilde{\mathcal{R}}_{q,\ell;j}(s, w) =   \sum_{(c, d) = 1} \sum_{\substack{q \mid r_1r_2c\\ \ell \mid r_1r_2d}}      \frac{\tau_3(r_2)\mu(r_1) }{r_1^{1+w+s}r_2^{s+w}c^{1+2s} d^{w -s} }\sum_{\substack{a_1, a_2 \, (\text{mod }c)\\  c \mid a_1a_2r_1}} \tilde{\rho}_j\left(\frac{a_1}{c}, \frac{a_2}{c}, \frac{d}{r_1c^{2}}\right).$$
Again we order by $\alpha_i = (a_i, c)$ getting
$$\tilde{\mathcal{R}}_{q,\ell;j}(s, w) =  \sum_{(c, d) = 1} \sum_{\substack{q \mid r_1r_2c\\ \ell \mid r_1r_2d}}     \frac{\tau_3(r_2)\mu(r_1) }{r_1^{1+w+s}r_2^{s+w}c^{1+2s} d^{w -s} } \sum_{\substack{\alpha_1, \alpha_2 \mid c\\ c \mid \alpha_1\alpha_2 r_1}}   \tilde{R}_j\left(\frac{\alpha_1}{c}, \frac{\alpha_2}{c}, \frac{d}{r_1c^{2}}\right)$$
with $$\tilde{R}_j(n_1, n_2, m) = (-1)^j R_j(n_1, n_2, m) =  \sum_{\nu+\mu+\kappa \leq 3-j} \tilde{C}_{\nu, \mu, \kappa, j}\phi_{\nu}(n_1) \phi_{\mu}(n_2) (\log m)^{\kappa}$$
for certain constants $\tilde{C}_{\nu, \mu, \kappa, j}$. Removing the coprimality condition $(c, d) = 1$ by M\"obius inversion and changing variables similarly as in the preceding proof we are left with terms of the form
\begin{displaymath}
\begin{split}
&   \sum_{\substack{q \mid r_1r_2b c\\ \ell \mid r_1r_2b d}}      \frac{\tau_3(r_2)\mu(r_1) \mu(b) }{r_1^{1+w+s}r_2^{s+w}c^{1+2s} b^{1+w+s}d^{w -s} } \sum_{\substack{\alpha_1, \alpha_2 \mid b c\\ \alpha_1\alpha_2 \mid b c r_1}}    \phi_{\nu}(\alpha_1) \phi_{\mu}(\alpha_2) (\log  d(b r_1)^{-1} c^{-2})^{\kappa}.
\end{split}
\end{displaymath}
As before we define the generalized function
\begin{equation}\label{tildeZ}
\tilde{Z}_{q, \ell}(x_1, x_2, x_3) =        \sum_{\substack{q \mid r_1r_2b c\\ \ell \mid r_1r_2b d}}  \sum_{\substack{a_1b_1, a_2b_2  \mid b c\\ a_1b_1a_2b_2 \mid b c r_1}}    \frac{\tau_3(r_2)\mu(r_1) \mu(b) \mu(a_1) \mu(a_2)b_1b_2 }{r_1^{1+w+s}r_2^{s+w}c^{1+2s} b^{1+w+s}   d^{w -s} } b_1^{x_1} b_2^{x_2} (  d (b r_1)^{-1} c^{-2})^{x_3}. 
\end{equation}
We compute the generic Euler factor at primes $p \nmid \ell q$ to be
\begin{equation*} 
 \frac{(1 - p^{-s-w})^{-3}(1 - p^{-(w-s-x_3)})^{-1}(1-p^{-(2s-x_1+2x_3)})^{-1}(1-p^{-(2s-x_2+2x_3)})^{-1}}{(1 - p^{-(s+w+x_3-x_1)})^{-1}(1 - p^{-(s+w + x_3-x_2)})^{-1}(1 - p^{-(2s+1+2x_3)})^{-1}}.
 \end{equation*}
This is the same computation with different exponents after renaming the variables $(a, c, d, n_1, r) \rightarrow (r_1, r_2, b, c, d)$ in \eqref{defZ}. 
As before, we conclude that $\tilde{\mathcal{R}}_{q, \ell;j}(s, w)$ is meromorphic. 

The estimation of \eqref{hol1} requires slightly more care, because the region $\Re w \geq \Re s \geq 1/2$, $x_1 = x_2 = x_3 = 0$ just fails to be inside the region of absolute convergence of the individual Euler factors. However, the only problem is caused by the $d$-sum in \eqref{tildeZ}. 

Let first $p\mid q$. Summing over $d$ (restricted to powers of $p$) first,  the same computation as before shows that  the $p$-Euler factor equals $(1 - p^{-(w-s-x_3)})^{-1}$ times a holomorphic expression that is bounded by $\ll (p^{v_p(q)})^{-1+\varepsilon}$ in an $\varepsilon$-neighbourhood of $\Re w \geq \Re s \geq 1/2$, $x_1 = x_2 = x_3 = 0$.  

Suppose now that $p\mid \ell$. 
We split the sum into $v_p(\ell) + 1$ terms with $p^{\lambda} \mid d$ for $0 \leq \lambda \leq v_p(\ell)$. In each of these we sum over $d$ first, and then estimate the rest as before. 
This gives an expression of the form $p^{\lambda(w-s-x_3)}(1 - p^{-(w-s-x_3)})^{-1}$ times an  absolutely convergent sum that is uniformly bounded in an $\varepsilon$-neighbourhood of $\Re w \geq \Re s \geq 1/2$, $x_1 = x_2 = x_3 = 0$. 

This yields the desired bound for \eqref{hol1} and  completes the proof.

 \section{Admissible functions}\label{admissible}

We call a function $H : \Bbb{R}_{> 0} \rightarrow \Bbb{C}$ \emph{admissible of type $(A, B)$} for some  $A, B>5$ if it is of one of the following two types:
\begin{itemize}
\item {[first kind]} we have $H(x) = \mathscr{K}h(x)$ (with $\mathscr{K}$ as in \eqref{H}) where $h$ is even and holomorphic in $|\Im t| < A$, such that   $h(t) \ll (1 + |t|)^{-B-2}$   and   $h$ has zeros at $\pm (n+1/2)i$, $n \in \Bbb{N}_0$, $n+1/2 < A$;
\item  {[second kind]} we have 
\begin{equation}\label{secondtype}
  H(x) = \alpha_0 J_a(4\pi x)(4\pi x)^{-b} + \sum_{a-b < k \in 2\Bbb{N}} \alpha_k J_{k-1}(4\pi x)
\end{equation}  
   for a constants $\alpha_0 \in \Bbb{C}$,  $a, b \in \Bbb{N}$ with $a - b \geq A$, and $\alpha_k \ll k^{-B-2}$.  
\end{itemize}
If $\{H_j\} = \{\mathscr{K}h_j\}$ is a family of admissible functions of type $(A, B)$ (of the first kind), we call it \emph{uniformly admissible of type $(A, B)$} if the implied constant in the required bound $h_j(t) \ll (1 + |t|)^{-B-2}$   
 can be chosen independently of $j$. 
 
We call a weakly admissible (i.e.\ satisfying \eqref{weakly}) pair $\mathfrak{h} = (h, h^{\text{hol}})$ \emph{admissible} if $h(t)$ is holomorphic in an $\varepsilon$-neighbourhood of $|\Im t| \leq 1/2$  and 
 $  \mathscr{K}^{\ast}\mathfrak{h}$ (with $\mathscr{K}^{\ast}$ as in \eqref{kast}) 
is admissible of type\footnote{We made no effort to optimize the requirements on $\mathfrak{h}$; most likely the number 500 could be reduced to 50, say, with no   additional effort.} (500, 500), i.e.\ of one of the above two kinds. Such a pair satisfies in particular the assumptions of \eqref{kuz-all}.   
We call a family of pairs of the  special shape $\{\mathfrak{h}_j = (h_j, 0)\}$ uniformly admissible if the family $\{h_j\}$ is uniformly admissible of type $(500, 500)$. 

From now on, we agree on the \textbf{convention} that all implied constants may depend on admissible weight functions $\mathfrak{h}$ where applicable; however, if $\{\mathfrak{h}_j = (h_j, 0)\}$ is a uniformly admissible family,   implied constants can always be chosen independently of $j$ (this applies also to \eqref{required-bound}, for instance). We will consider a uniformly  admissible family only once in this paper, in the proof of Theorem \ref{cor1}.

 \begin{lemma}\label{lem1} A pair $\mathfrak{h} = (h, h^{\text{{\rm hol}}})$  is admissible if \\
{\rm a)}   $\mathfrak{h} = (0, \delta_{k = k_0})$ for some $k_0 > 500$;\\
{\rm b)}  $\mathfrak{h} = (h, 0)$, where $h$ is  even and holomorphic in $|\Im t| < 500$, such that   $h(t) \ll (1 + |t|)^{-502}$    and   $h$ has zeros at $\pm (n+1/2)i$, $n \in \Bbb{N}_0$, $n+1/2 < 500$;\\
{\rm c)} $\mathfrak{h} = (h_{\text{{\rm pos}}}, h_{\text{{\rm pos}}}^{\text{{\rm hol}}})$ as in \eqref{pos-id} with $a = 1000$, $b = 400$. This function is strictly positive on   $\mathcal{T}_0 = (\Bbb{R} \cup [-i\vartheta, i \vartheta]) \times 2\Bbb{N}$. 
\end{lemma}
This is clear from the definition, cf.\ also \eqref{pos}. A typical  function satisfying the hypotheses of Lemma \ref{lem1}b) is  $$h(t) = e^{-t^2} \prod_{n \leq 500} \big(t^2 + (n+ \textstyle \frac{1}{2})^2\big).$$

We need two technical properties of admissible functions that are presented in the following two lemmas.  
 
\begin{lemma}\label{admis1} Let $H$ be an admissible function of type $(A, B)$.  \\
{\rm a)} The Mellin transform $\widehat{H}$ is holomorphic in $\Re u > -A$ and in this region bounded by 
\begin{equation}\label{assumption-c}
\widehat{H} (u) \ll_{\Re u} (1 + |u|)^{\Re u - 1} .
\end{equation} 
{\rm b)} For any $M > 0$,  we have the asymptotic formula
\begin{equation}\label{asymp-j}
\widehat{H}(\sigma + it) = |t|^{\sigma - 1} \exp\left(it \log \frac{|t|}{4\pi e}\right) j_{\sigma}(t) + O_{\sigma}(|t|^{\sigma-1 -B})
\end{equation}
for $\sigma > -A$, $|t| \geq 30$  and a smooth function $j_{\sigma}$ satisfying for all $\nu\leq B$ the bound
\begin{equation}\label{ad3}
|t|^{\nu} j^{(\nu)}_{\sigma}(t) \ll_{\sigma} 1.
\end{equation}
\end{lemma} 
 
As per our convention,  if $\{H_j\} = \{\mathscr{K}h_j\}$ is uniformly admissible of type $(A, B)$, then all implied constants in the previous lemma can be chosen independently of $j$. \\

\textbf{Proof.} We treat the two kinds of admissible functions separately. 

1) Suppose that $H = \mathscr{K}h$ is of the first kind. For $  \Re u > 0$ we have by  \eqref{mellin-j} and \eqref{H} that
 \begin{equation}\label{mel}
 \begin{split}
 \widehat{H}(u) 
 & = \frac{i}{4\pi} \int_{-\infty}^{\infty}    (2\pi)^{-u}\left( \frac{\Gamma(u/2 + i\tau)}{\Gamma(1 - u/2 + i\tau)}  -  \frac{\Gamma(u/2 - i\tau)}{\Gamma(1 - u/2 - i\tau)}\right) \frac{h(\tau) \tau}{\cosh(\pi \tau)} \, d\tau\\
 & = \frac{i}{2\pi} \int_{-\infty}^{\infty}    (2\pi)^{-u} \frac{\Gamma(u/2 + i\tau)}{\Gamma(1 - u/2 + i\tau)}  \frac{h(\tau) \tau}{\cosh(\pi \tau)} \, d\tau ,
 \end{split}
 \end{equation}
 as an absolutely convergent integral. 
 We shift the contour to the line $\Im \tau = -A + 1/8$. This does not cross poles of $1/\cosh(\pi \tau)$ because of the zeros of $h$. In this way see that the integral is holomorphic in $ \Re u> -2A + 1/4$ and bounded by $O((1+|u|)^{\Re u-1})$, with a uniform implied constant if $\{\mathscr{K}h_j\}$ is a uniform family.  This confirms a).   

 Moreover, for $u = \sigma + it$ with $\sigma > -A$, $|t| \geq  30$ and $r= (\sigma-1)/2$ (say) we have
\begin{displaymath}
\begin{split} 
 &\widehat{H}(\sigma+it) |t|^{-\sigma+1} \exp\left(-it  \log\frac{|t|}{4\pi e}\right) \\
 &= \frac{i}{2\pi} \int_{-\infty}^{\infty} (2\pi)^{-\sigma} |t|^{-\sigma-1} \exp\left(-it  \log \frac{|t|}{2 e}\right) \frac{\Gamma(\frac{1}{2}\sigma - r + i(\frac{1}{2}t + \tau))}{\Gamma(1 - \frac{1}{2}\sigma - r + i(\tau - \frac{1}{2}t))} \frac{h(\tau+ir) (\tau+ir)}{\cosh(\pi(\tau + ir))} d \tau,
 \end{split}
 \end{displaymath}
 having shifted the contour  in \eqref{mel} to imaginary part $r$ (as we need to obtain the analytic continuation to $\Re u = \sigma$). 
 We want to show that except for an error of $O(|t|^{-B})$, this function satisfies \eqref{ad3} for $\nu \leq B$ (and fixed $r$ and $\sigma$). By trivial estimates using $h(t) \ll (1 + |t|)^{-B-2}$,  the portion $|\tau | > \frac{1}{4}|t|$ contributes $O(|t|^{-B})$. For the remaining portion  we can insert Stirling's formula \eqref{stir} for the gamma quotient. The error term in \eqref{stir} can   be bounded by any negative power of $|t|$, and the main term   gives us a phase  $\exp(i \phi_{\tau}(t))$ with
 $$ \phi_{\tau}(t) = -t  \log \frac{|t|}{2 e}+ \frac{1}{2}t   \log\frac{|t^2 - (2\tau)^2|}{(2e)^2}  + \tau \log\Big| \frac{t+2\tau}{t -2\tau}\Big|$$
 satisfying
 $$ \phi_{\tau}'(t) = \frac{1}{2} \log\Big|\frac{t^2 - (2\tau)^2}{t^2}\Bigr| \ll \frac{|\tau|}{|t|}, \quad  \phi^{(j)}_{\tau}(t) \ll \frac{|\tau|}{|t|^{j}}\,\,\, (j\geq 2)$$
uniformly in $|\tau| \leq \frac{1}{4}|t|$. Now we  differentiate under the integral sign, and since   
$$|t|^{\nu} \frac{d^{\nu}}{dt^{\nu}} \left( |t|^{-\sigma-1} \exp\left(-it  \log \frac{|t|}{2 e}\right) \frac{\Gamma(\frac{1}{2}\sigma - r + i(\frac{1}{2}t + \tau))}{\Gamma(1 - \frac{1}{2}\sigma - r + i(\tau - \frac{1}{2}t)) \cosh(\pi (\tau + ir))}\right)  \ll 1 + |\tau|,$$
and $\int_{\Bbb{R}} |h(\tau + ir)| (1 + |\tau|)^2 d\tau \ll 1$, we obtain \eqref{ad3} as desired. This confirms b).

 2)  If $H$ is of the form \eqref{secondtype}, then
\begin{equation}\label{second}
\widehat{H}(u) = \alpha_0  \frac{\Gamma(\frac{1}{2}(a+ u -b))}{2(2\pi)^{u}2^b\Gamma(\frac{1}{2}(2+a-u+b ))} + \sum_{a-b < k \in 2\Bbb{N}} \alpha_k   \frac{\Gamma(\frac{1}{2}(k-1 + u ))}{2(2\pi)^{u}\Gamma(\frac{1}{2}(1+k-u ))},
\end{equation}
 and a) is clear. By Stirling's formula, the first term on the right hand side  satisfies \eqref{asymp-j} and \eqref{ad3}. To treat the second term, we write for $k \in 2\Bbb{N}$ and $u = \sigma +it$ the last fraction in the form
 $$\frac{\Gamma(\frac{1}{2}(k-1 + u ))}{ \Gamma(\frac{1}{2}(1+k-u ))} =\frac{\Gamma(\frac{1}{2}(1+\sigma + it))}{\Gamma(\frac{1}{2}(3 - \sigma - it))} \prod_{n=1}^{\frac{k}{2} - 1} \frac{k-1 + \sigma- 2n + it}{k-1 -\sigma - 2n - it}. $$
We use Stirling's formula for the gamma fraction on the right hand side. In order to verify \eqref{asymp-j} and \eqref{ad3}, it remains to show that
\begin{equation}\label{verify} \sum_{k\in 2\Bbb{N}}   \Big| \alpha_k t^{\nu} \frac{d^{\nu}}{dt^{\nu}}   \prod_{n=1}^{k/2 - 1} \left(\frac{k-1 + \sigma- 2n + it}{k-1 -\sigma - 2n - it} \right)\Big| \ll 1
\end{equation}
for $\nu \leq B$. 
For $\nu \geq 1$ we have
$$\frac{d^{\nu}}{dt^{\nu}}  \left(\frac{k-1 + \sigma- 2n + it}{k-1 -\sigma - 2n - it} \right) = 2i^{\nu} \nu!  \frac{k-1-2n}{(k-1-2n-\sigma - i t)^{\nu+1}} \ll_{\sigma} \frac{k}{(1+|t|)^{\nu+1}},$$
and so by Leibniz' rule the left hand side of \eqref{verify} is bounded by $\sum_k |\alpha_k| k^{\nu} \ll 1$ for $  \nu \leq B$. 
 

\begin{lemma}\label{sec-prop} Let $c, d \in \Bbb{C}$, $x \in \Bbb{R} \setminus \Bbb{N}_0$. 
Let $H$ be an admissible function of type $(A, B)$. Suppose that $2x + \Re c > -A$, $\Re d   > \theta+ \max(x, (x+1)/2)$. Then
$$\mathscr{I}^{\pm}(c, d) :=  \int_{(x)} G^{\pm}(\xi) \widehat{H}  (c+2\xi) \mathcal{G}_{-\mu}^{\mp}(d-\xi) \frac{d\xi}{2\pi i}$$
with  $G^{\pm}$ as in \eqref{Gpm} and $\mathcal{G}^{\pm}_{\mu}$ as in \eqref{defgmu}  is holomorphic  for $\Re(c+3d) < 2$  and has a meromorphic continuation to $\Re(c + 3d) < 3$  with a simple polar divisor at most  at $c+3d = 2$.\end{lemma}

\textbf{Proof.} We recall that the spectral parameter $\mu = (\mu_1, \mu_2, \mu_3)$ of the  automorphic form $F$  satisfies  $|\Re \mu_j|\leq \theta \leq 5/14$.  Throughout the proof we always assume   that $c$ and $d$ satisfy
\begin{equation}\label{cond-lemma}
  2x + \Re c > -A,  \quad \Re d   > \theta+ \max(x, (x+1)/2).
  \end{equation} The first condition ensures that we are in the domain of holomorphicity of $\widehat{H}$, the condition  $\Re d   > \theta+ x$ ensures that $ \mathcal{G}_{-\mu}^{\mp}(d-\xi)$ is holomorphic, and the condition $\Re d   > \theta+ (x+1)/2$ will be needed later. 
 
 Holomorphy in $\Re(c+3d) < 2$  follows   immediately from \eqref{bound-g}, \eqref{bound-G} and Lemma \ref{admis1}a, but for $\Re(c+3d) \geq 2$, the integral fails to converge absolutely.   Fix any $\beta > 0$.  Clearly it suffices to continue the restricted  integral with the compact region $|\Im \xi| \leq D :=  2\beta + 2\max_j |\Im \mu_j| + 2$ removed, meromorphically to the region $|\Im d| \leq  \beta$ and $\Re(c + 3d) < 3$.
 
To this end, we will approximate the integrand for $  |\Im \xi| \geq D$ by a simpler expression. The error in this approximation will  decay more quickly by an additional power of $1/|\xi|$, which buys us absolute convergence in $\Re(c + 3d) < 3$. We will then complete the simplified integral by re-inserting the portion $|\Im \xi| \leq D$ and evaluate it explicitly by \eqref{hyper}. 
In the following we treat only  the case $\mathscr{I}^+(c, d)$ and drop the superscript, the other case is analogous. 

 Using the definitions \eqref{Gpm} and \eqref{defgmu}, we have
  $$ G^{+}(\xi) 
  \mathcal{G}_{-\mu}^{-}(d-\xi) =   \Gamma(\xi) (2\pi)^{ 2\xi} \Bigl( \prod_{j=1}^3\Gamma(d - \xi -\mu_j) \Bigr)\sum_{\nu = -3}^4 \gamma_{\nu}(d) e^{i \frac{\pi}{2} \nu \xi}.  $$ 
  for certain holomorphic functions $\gamma_{\nu}(d)$  with $\gamma_{4}(d) = (2\pi)^{-3d}$
   (that in general depend also on the spectral parameter $\mu$ which we suppressed from the notation).  The terms $-3 \leq \nu \leq 3$ lead to exponentially decaying integrals and can easily be continued homomorphically to $\Re(c + 3d)$ arbitrarily large.   
  For notational simplicity let us write $d_j := d - \mu_j$.  We have 
  $$\prod_{j=1}^3\Gamma(d_j - \xi )  = \Gamma(1 - \xi) \Gamma(d_1 + d_2-1  - \xi) \frac{\Gamma(d_1 - \xi) \Gamma(d_2 - \xi)}{ \Gamma(1 - \xi) \Gamma(d_1 + d_2-1  - \xi)}  \frac{\pi}{\sin(\pi(d_3 - \xi)) \Gamma(1 - d_3 + \xi)}$$
by the reflection formula, and by another application of the reflection formula we see that the 
 term corresponding to $\nu = 4$ equals
   \begin{equation*}
    (2\pi)^{2\xi } \frac{\Gamma(d_1 + d_2 - 1 - \xi)}{\Gamma(1 - d_3 + \xi) } \frac{\pi^2 (2\pi)^{-3d} e^{2\pi i \xi}}{\sin(\pi \xi) \sin(\pi (d_3 - \xi))}  \frac{\Gamma(d_1 - \xi) \Gamma(d_2 - \xi)}{ \Gamma(1 - \xi) \Gamma(d_1 + d_2-1  - \xi)}.
    \end{equation*}
 The basic observation   is now that the last gamma fraction is asymptotically constant as $|\Im \xi| \rightarrow \infty$. More precisely, by  \eqref{stir} -- \eqref{stir1}, the difference
  $$   \frac{\Gamma(d_1 - \xi) \Gamma(d_2 - \xi)}{ \Gamma(1 - \xi) \Gamma(d_1 + d_2-1  - \xi)}  -   1,  
 $$
 is holomorphic in $|\Im d| < \beta$, $|\Im \xi| > D$ and bounded by $O_{d, \Re \xi}(|\xi|^{-1})$. Consequently, the difference
 $$  \frac{\pi^2 (2\pi)^{-3d} e^{2\pi i \xi}}{\sin(\pi \xi) \sin(\pi (d_3 - \xi))}  \frac{\Gamma(d_1 - \xi) \Gamma(d_2 - \xi)}{ \Gamma(1 - \xi) \Gamma(d_1 + d_2-1  - \xi)}  - \begin{cases} 0, & \Im \xi > D,\\  4\pi^2 (2\pi)^{-3d}  e^{i\pi d_3},  
 & \Im \xi < - D,\end{cases}$$
 is holomorphic in $|\Im d| < \beta$ (because we stay away from zeros  of the $\sin$-function) and bounded by $O_{d, \Re \xi}(|\xi|^{-1})$. Since
 $$ (2\pi)^{2\xi } \frac{\Gamma(d_1 + d_2 - 1 - \xi)}{\Gamma(1 - d_3 + \xi) } 
 \widehat{H}  (c+2\xi) |\xi|^{-1} \ll_{\Re \xi, d, \Re c} |\xi|^{\Re(c + 3d) - 4}$$
 by \eqref{assumption-c}, this  portion can be continued to $\Re (c+3d) < 3$ as an absolutely (and locally uniformly in $d$) convergent integral over $|\Im \xi| > D$. Hence we are left with 
 $$4\pi^2 (2\pi)^{-3d}  e^{i\pi d_3} \int_{\substack{\Re \xi = x\\ |\Im \xi| > D}}   \widehat{H}  (c+2\xi)   (2\pi)^{2\xi } \frac{\Gamma(d_1 + d_2 - 1 - \xi)}{\Gamma(1 - d_3 + \xi) }  (1 - \text{sgn}(\Im \xi))  \frac{d\xi}{2\pi i},$$
and of course the constant in front of the integral plays no role. In order to remove the unwanted factor $ 1 - \text{sgn}(\Im \xi)$, we use that
 $$-i\cdot  \text{sgn}(\Im \xi) - \frac{\Gamma(1 - d_3 + \xi) \Gamma(d_1 + d_2 - 1/2 - \xi)}{\Gamma(3/2 - d_3 + \xi)\Gamma(d_1 + d_2 -1 - \xi)} $$
is holomorphic in $|\Im d - \beta| < 1$ for $|\Im \xi| > D$ and bounded by $O_{d, \Re \xi}(|\xi|^{-1})$, 
 which follows again from \eqref{stir} -- \eqref{stir1}. Hence our job is reduced to continuing
 $$\int_{\substack{\Re \xi = x\\ |\Im \xi| > D}}   \widehat{H}  (c+2\xi)   (2\pi)^{2\xi } \frac{\Gamma(d_1 + d_2 - 1+ \kappa  - \xi)}{\Gamma(1 + \kappa - d_3 + \xi) }    \frac{d\xi}{2\pi i}$$
 for $\kappa \in \{0, 1/2\}$, and we may now re-insert the compact interval $|\Im \xi| 
 leq D$. At this point the condition $\Re d   > \theta+ (x+1)/2$ comes handy, since the completed  contour is still to the right of all poles. 
 Thus we let 
 $$\mathscr{I}_{\kappa}(c, d) := \int_{(x)}   \widehat{H}  (c+2\xi)   (2\pi)^{2\xi } \frac{\Gamma(d_1 + d_2 - 1+ \kappa  - \xi)}{\Gamma(1 + \kappa - d_3 + \xi) }    \frac{d\xi}{2\pi i}$$
 and  treat the two kinds of admissible functions separately. 

If $H = \mathscr{K}h$ is of the first kind, then by \eqref{mel} we have
$$\mathscr{I}_{\kappa}(c, d )= (2\pi)^{-c-1} i \int_{\Im \tau =-A+1/2  }  \int_{(x)}  \frac{\Gamma(d_1 + d_2 - 1+\kappa - \xi)}{\Gamma(1+\kappa - d_3 + \xi) }    \frac{\Gamma(\frac{c}{2}+\xi + i\tau)}{\Gamma(1-\frac{c}{2} - \xi + i\tau)}   \frac{d\xi}{2\pi i}  \frac{h(\tau) \tau}{\cosh(\pi \tau)} \, d\tau,$$
where the exchange of integration is easily justified by absolute convergence and we are allowed to shift the $\tau$-contour  since $H$ is admissible of type $(A, B)$. The conditions \eqref{cond-lemma} ensure that the arguments of the gamma factors in the numerator are to the right of all poles. 
  By \eqref{hyper} and \eqref{mu}, we can evaluate the $\xi$-integral getting
$$\mathscr{I}_{\kappa}(c, d ) = (2\pi)^{-c-1} i \int_{\Im \tau =-A+1/2 } \frac{\Gamma(d_1+d_2 - 1+\kappa + \frac{c}{2} + i\tau)\Gamma( 2-c-3d)}{\Gamma(2-3d)\Gamma(1-c)\Gamma( 1+\kappa -d_3 - \frac{c}{2} + i\tau )} \frac{h(\tau) \tau}{\cosh(\pi \tau)} \, d\tau.$$
This expression is meromorphic in $\Re (c+3d) < 3$ (intersected with \eqref{cond-lemma})   with its only polar divisor at $c+3d = 2$. (Note that \eqref{cond-lemma} implies that the argument of the first gamma factor in the numerator has positive real part.)

If $H$ is of the second kind, we argue similarly: by \eqref{second}, $\mathscr{I}_{\kappa}(c, d)$  equals
\begin{displaymath}
\begin{split}
 \frac{1}{2} \int_{(x)}&  \Bigl[ \alpha_0  \frac{(2\pi)^{-c}\Gamma(\xi + \frac{1}{2}(a -b+c))}{2^b \Gamma(-\xi+\frac{1}{2}(2+a+b-c ))} + \sum_{a-b < k \in 2\Bbb{N}} \alpha_k   \frac{(2\pi)^{-c} \Gamma(\xi+\frac{1}{2}(k-1 +c ))}{\Gamma(-\xi + \frac{1}{2}(1+k -c))}\Bigr]  \frac{\Gamma(d_1 + d_2 - 1+ \kappa  - \xi)}{\Gamma(1 + \kappa - d_3 + \xi) }   \frac{d\xi}{2\pi i}.
\end{split}
\end{displaymath}
The first term in parentheses yields an integral that is absolutely convergent in $\Re (c+3d) < 2 + b$. For the second term we exchange sum and integration and evaluate the $\xi$-integral explicitly by \eqref{hyper}, which provides the meromorphic continuation to $\Re (c+3d) < 3$ with a  polar divisor only at $c+3d = 2$.

\section{An integral transform}\label{int-trafo}

  Fix  ${\tt s}, {\tt w}\in \Bbb{C}$ and  suppose that $H$ is admissible of type $(A, B)$.  
     Then for all ${\tt u}$ satisfying
  \begin{equation}\label{cond}
  \begin{split}
 &\textstyle -A < \Re(3{\tt s} - {\tt w}   +   {\tt u} -3), 
 \quad  \theta < \Re ({\tt s} + \frac{1}{2}{\tt u}) < \frac{1}{2},  \quad  \Re (\frac{1}{2} {\tt u} +   {\tt w}) > 0 
 \end{split}
    \end{equation}
 the integral transforms
\begin{equation}\label{trafo}
\begin{split}
(\mathscr{V}_{{\tt s}, {\tt w}}^{\pm} \widehat{H})({\tt u}) :=    \int_{\mathcal{C}}   \widehat{H}(3{\tt s} - {\tt w} - 1 + {\tt u} + 2\xi) &\Big[G^+(\xi) \mathcal{G}_{-\mu}^-(1 - {\tt s} - \textstyle\frac{1}{2} {\tt u} - \xi)\mathcal{G}_{\mu}^{\mp}({\tt s} + \frac{1}{2}{\tt u}) \\
&+ G^-(\xi) \textstyle\mathcal{G}_{-\mu}^+(1 - {\tt s} - \frac{1}{2} {\tt u}- \xi)\mathcal{G}_{\mu}^{\pm}({\tt s} + \frac{1}{2}{\tt u})\Big] \displaystyle\frac{d\xi}{2\pi i}  \end{split}
\end{equation}
with  $G^{\pm}$ as in \eqref{Gpm} and $\mathcal{G}^{\pm}_{\mu}$ as in \eqref{defgmu}  define an absolutely convergent integral.  Indeed, the first   condition in \eqref{cond} implies that the argument of $\widehat{H}$ is in the region of holomorphicity as provided in Lemma \ref{admis1}a, while the third  condition ensures,  by  Lemma \ref{lemG}  and \eqref{assumption-c}, absolute convergence. It is important to note that this condition is independent of $\xi$ and in particular not affected by possible contour shifts. Finally, the second condition in \eqref{cond} ensures that as ${\tt u}$ varies no poles are crossed. We would like to continue $(\mathscr{V}_{{\tt s}, {\tt w}}^{\pm} \widehat{H})({\tt u})$ holomorphically to a larger region of $({\tt u}, {\tt s}, {\tt w})$ and also obtain good bounds in terms of ${\tt u}$.    All implied constants in the following may depend on ${\tt s}, {\tt w}, A, B, \mu$. 
 
 We remark already at this point that this integral will come up later with the following slightly deformed contour  consisting of four line segments
 \begin{equation}\label{contour}
 \mathcal{C} = \textstyle (-\frac{3}{5} -i\infty, -\frac{3}{5} - i] \cup [-\frac{3}{5} - i, \frac{1}{10}] \cup [\frac{1}{10}, -\frac{3}{5} + i] \cup [-\frac{3}{5} + i, -\frac{3}{5} + i\infty).
 \end{equation} 
 The first condition in \eqref{cond} has enough elbow room to make sure that the  argument of $\widehat{H}$ is still in the region of holomorphicity.

\begin{lemma}\label{lem9}  Let ${\tt s}, {\tt w}\in \Bbb{C}$,  $A, B > 2$, suppose that $H$ is admissible of type $(A, B)$ and suppose that the set ${\tt u}$ satisfying \eqref{cond}  is non-empty. 
Then $(\mathscr{V}_{{\tt s}, {\tt w}}^{\pm} \widehat{H})({\tt u})$ can be extended holomorphically to all $({\tt u}, {\tt s}, {\tt w}) \in \Bbb{C}^3$ satisfying
\begin{equation}\label{region}
\begin{split}
&\textstyle  \theta   < \Re ({\tt s} + \frac{1}{2}{\tt u}) < \frac{1}{2}A + \frac{1}{2}\Re({\tt s} - {\tt w}), \quad \Re(  \frac{1}{2} {\tt u} + {\tt w}) > 0
\end{split}
\end{equation}
and  satisfies\footnote{Again, if $\{H_j \}$ is uniformly admissible of type $(A, B)$,   the implied constants  can be chosen independently of $j$.}
 \begin{equation}\label{boundJ}
\begin{split}
(\mathscr{V}_{{\tt s}, {\tt w}}^{\pm} \widehat{H})({\tt u}) \ll    
 (1 + |u|)^{3 (|\Re{\tt s} | + |\Re{\tt u} |+ |\Re{\tt w} |) - \min(\frac{1}{8} B, \frac{1}{4} A)}  
 \end{split}
 \end{equation}
%
in this region. Moreover, $(\mathscr{V}_{{\tt s}, {\tt w}}^{\pm} \widehat{H})({\tt u})$ can be extended meromorphically to 
\begin{equation}\label{mero}
\begin{split}
&\textstyle -1+ \theta   < \Re ({\tt s} + \frac{1}{2}{\tt u}) < \frac{1}{4}A + \frac{1}{4}\Re({\tt s} - {\tt w}), \quad \Re(  \frac{1}{2} {\tt u} + {\tt w}) > -1
\end{split}
\end{equation}
with   poles at most at ${\tt u}/2 + {\tt w} = 0$ and $ {\tt u}/2 + {\tt s} + \mu_j = 0$, $j = 1, 2, 3$.
\end{lemma}


\textbf{Proof.} Before we start with the proof, we use \eqref{new-formula} to  compute explicitly the residue of the integrand in \eqref{trafo} at $\xi = -n$, $n \in \Bbb{N}_0$, as
\begin{equation}\label{resi}
\begin{cases} 0 , & \pm = +,\\ \widehat{H}(3{\tt s} - {\tt w} - 1  + {\tt u} - 2n) (2\pi)^{-2n}  (n!)^{-1} \prod_{j=1}^n \prod_{k=1}^3({\tt s} + \frac{1}{2}{\tt u} +\mu_k - j),& \pm = -. \end{cases}
\end{equation}

The proofs of \eqref{boundJ} and \eqref{mero} need a slightly different treatment depending on the sign, and we start with $(\mathscr{V}_{{\tt s}, {\tt w}}^{-} \widehat{H})({\tt u})$. Let initially $({\tt u}, {\tt s}, {\tt w})$ be in the region \eqref{cond}.  We straighten the   $\xi$-contour and shift it to the far left to $\Re \xi = -A_1 \not\in \Bbb{N}$. Since $  \Re ({\tt s} +  {\tt u}/2) > 0$, this does not leave the domain of holomorphicity of $\widehat{H}$   provided 
 \begin{equation}\label{a1}  \Re({\tt s} - {\tt w} - 1 -2A_1) > -A.
 \end{equation} We pick up   possible poles at $\xi = -n$, $n \in \Bbb{N}_0$, $n_0 < A_1$, whose residues are by \eqref{resi} holomorphic functions in ${\tt u}, {\tt s}, {\tt w}$ in the region \eqref{a1}. The remaining integral is holomorphic in the region 
\begin{equation}\label{newregion}
\textstyle \theta < \Re ({\tt s} + \frac{1}{2}{\tt u}) < 1 - \theta + A_1, \quad  \Re (\frac{1}{2}{\tt u} +  {\tt w}) > 0.  
\end{equation}
 This gives a holomorphic continuation of $(\mathscr{V}_{{\tt s}, {\tt w}}^{-} \widehat{H})({\tt u})$ in $({\tt u}, {\tt s}, {\tt w})$ to the intersection of the regions \eqref{a1} and \eqref{newregion}. Choosing $A_1 = \frac{1}{2}(A + \Re ({\tt s}- {\tt w}) - 1) - 1/100$, we obtain a region containing \eqref{region}. 
 
 The same argument combined with Lemma \ref{sec-prop} shows that $(\mathscr{V}_{{\tt s}, {\tt w}}^{-} \widehat{H})({\tt u})$ has a meromorphic continuation to 
 $$\textstyle -1 + \theta < \Re ({\tt s} + \frac{1}{2}{\tt u}) < 1 - \theta + \frac{1}{2}A_1, \quad \Re (\frac{1}{2}{\tt u} +  {\tt w}) > -1$$
 with poles at most at ${\tt u}/2 + {\tt w} = 0$ and $ {\tt u}/2 + {\tt s} + \mu_j = 0$, $j = 1, 2, 3$. 

Now we fix ${\tt s}$ and ${\tt w}$ and  proceed to estimate $(\mathscr{V}_{{\tt s}, {\tt w}}^{-} \widehat{H})({\tt u})$ for ${\tt u}$ satisfying \eqref{region} with the aim of   establishing \eqref{boundJ}. We write $v = \frac{1}{2} \Im{\tt u}$ and may assume that $v > 0$ is sufficiently large (in terms of $\mu$),  for if $|v|$ is bounded   there is nothing to prove, and the case of negative $v$ is similar.   We may then focus our attention to the second term in \eqref{trafo}, since by  \eqref{bound-G} the first term, which  contains a factor $\mathcal{G}^+_{\mu}({\tt s} + \frac{1}{2} {\tt u})$, is exponentially decaying as $v \rightarrow +\infty$.

We shift the $\xi$-contour back to the far right, to $\Re \xi = A_2 > 0$, say, making sure that this line does not cross poles. This cancels all residues from \eqref{resi} that we picked up earlier, but introduces potentially new residues at $1 - {\tt s} -  \frac{1}{2} {\tt u} - \xi = -\mu_j-n$, $n \in \Bbb{N}_0$, $j \in \{1,2, 3\}$. We do not compute them explicitly, but observe that a priori the residues must be meromorphic in $({\tt u}, {\tt s}, {\tt w})$;  since we know already that $(\mathscr{V}_{{\tt s}, {\tt w}}^{-} \widehat{H})({\tt u})$  is holomorphic, their joint contribution must be holomorphic, too. Moreover, for $1 - {\tt s} -  \frac{1}{2} {\tt u} - \xi  = O(1)$ but off the poles, the integrand in \eqref{trafo} is 
\begin{equation*}
   \ll v^{O(1)} e^{-\pi v/2}
 \end{equation*}  
    by \eqref{bound-g} and \eqref{bound-G};   by Cauchy's theorem this bound also holds for  the residues at $1 - {\tt s} -  \frac{1}{2} {\tt u} - \xi = -\mu_j-n$, which is   in agreement with \eqref{boundJ}.  This time we used the exponential decay of $G^-(\xi)$ in the second term of \eqref{trafo} under our current assumption $v > 0$.

      It remains to estimate the remaining integral on the line $\Re \xi  = A_2$. We split the integral smoothly into two pieces as follows. Let  
    $w : \Bbb{R} \rightarrow [0, 1]$ be a smooth function that is constantly 1 on $[-1, 1]$ and vanishes outside $[-2, 2]$. Then it suffices to estimate
$$I^-_1 := \int_{(A_2)}  \widehat{H}(3{\tt s} - {\tt w} - 1 + {\tt u} + 2\xi) G^-(\xi) \mathcal{G}_{-\mu}^+(1 - {\tt s} - \textstyle\frac{1}{2} {\tt u} - \xi)\mathcal{G}_{\mu}^{-}({\tt s} + \frac{1}{2} {\tt u})  \displaystyle w\left(\frac{\Im \xi}{V}\right)\frac{d\xi}{2\pi i}$$
and 
$$I^-_2 := \int_{(A_2)}  \widehat{H}(3{\tt s} - {\tt w} - 1 + {\tt u} + 2\xi) G^-(\xi) \mathcal{G}_{-\mu}^+(1 - {\tt s} - \textstyle\frac{1}{2} {\tt u} - \xi)\mathcal{G}_{\mu}^{-}({\tt s} + \frac{1}{2} {\tt u}) \displaystyle\left(1 - w\left(\frac{\Im \xi}{V}\right)\right) \frac{d\xi}{2\pi i},$$
where $1 \leq V \leq v$ will be chosen in a moment. 
 We estimate $I_1^-$ trivially  by \eqref{assumption-c}, \eqref{bound-g} and \eqref{bound-G} getting  \begin{equation*}
I^-_1 \ll  v^{\Re(3{\tt s} - {\tt w} - 2 + {\tt  u} - A_2) } V^{A_2 + \frac{1}{2}}.
\end{equation*}
Since we are free to choose $A_2$ as large as we wish, this is   admissible for \eqref{boundJ} provided $V \leq v^{1-\delta}$ for some fixed $\delta > 0$, which we assume from now on. For the estimation of $I_2^-$, we first observe that we can   restrict to the branch  $ \Im \xi > 0$, as the branch $\Im \xi < 0$ can   be estimated  trivially by
\begin{equation}\label{exp-sr}
   \ll v^{O(1)} e^{-\frac{1}{2}\pi V},
 \end{equation}  
 which is admissible for \eqref{boundJ} provided $V > v^{\delta}$ for some fixed $\delta > 0$. 
The treatment of the remaining case $\Im \xi > 0$, which we assume from now on,  is the only point where properties \eqref{asymp-j} and \eqref{ad3}  of $\widehat{H}$ are required. 
We insert the asymptotic formula \eqref{asymp-j} along with the asymptotic formulae \eqref{asymp-g} and \eqref{asymp-G}. The error terms corresponding to \eqref{asymp-g} and  \eqref{asymp-G} save arbitrarily many powers of $V> v^{\delta}$ which is admissible for  \eqref{boundJ}. The error term corresponding to \eqref{asymp-j}  contributes at most
\begin{equation}\label{E1}
\int_{x \geq V} (v+x)^{-  \Re( {\tt w} + \frac{1}{2} {\tt u}) - A_2 - \frac{1}{2} - B} x^{A_2 - \frac{1}{2}} v^{3\Re({\tt s} + \frac{1}{2} {\tt u}) - \frac{3}{2}}dx \ll v^{\Re (3{\tt s} +   {\tt u}  - {\tt w}) - \frac{3}{2} - B}.
\end{equation}
For the main terms, 
we obtain an integral of the shape 
\begin{displaymath}
\begin{split}
\mathcal{G}_{\mu}^+({\tt s} + \textstyle \frac{1}{2}  {\tt u})&\displaystyle \int_0^{\infty}  (x + v  )^{-\frac{1}{2}  - A_2 - \Re {\tt w} - \frac{1}{2} \Re {\tt u}} x^{A_2 - \frac{1}{2}} \omega(x) e^{i\phi(x)} 
 \left(1 - w\left(  \frac{x}{V}\right)\right) dx
\end{split}
\end{displaymath}
with $x^j \omega^{(j)} (x) \ll_j 1$ for all $j \in \Bbb{N}_0$ and
\begin{equation}\label{phi}
\phi(x) = x \log \displaystyle\frac{|x|}{2\pi e} -  (x + v) \displaystyle\log \frac{|x + v|}{2\pi e}.
\end{equation}
We compute $$\phi'(x) = \log\Big| \frac{x} {x+v}\Bigr|, 
\quad \phi^{(j)}(x) \ll \frac{|v|(|v| + |x|)^{j-2}}{|x|^{j-1} |v+x|^{j-1}}, \quad j \geq 2. $$
 We attach another smooth partition of unity and decompose the integral smoothly into dyadic pieces supported on $Z< x < 4Z$ with $Z \geq V$. For each piece  we can apply Lemma \ref{integrationbyparts} with 
 $$\beta-\alpha \asymp U = Q = Z, \quad Y = \frac{vZ}{v + Z}, \quad R = \frac{v}{v + Z}, \quad X = ( v   + Z)^{-\frac{1}{2}  - A_2 - \Re {\tt w} - \frac{1}{2} \Re {\tt u}} Z^{A_2 - \frac{1}{2}}$$
 (note that for $ Z \asymp x \geq v^{\delta}$ we have  $\log x/(x+v) \asymp v/(v+Z)$) and bound it by
 $$\ll v^{3(\Re {\tt s} + \frac{1}{2} \Re {\tt u} - \frac{1}{2})}(  v  + Z)^{-\frac{1}{2}  - A_2 - \Re {\tt w} - \frac{1}{2} \Re {\tt u}} Z^{A_2 + \frac{1}{2}}  \Bigl(\sqrt{\frac{Zv}{v + Z}}\Bigr)^{-B/2}.$$
 We estimate $Zv/(v + Z) \gg \min(v, Z) \gg V$ and bound the previous display by
 $$\ll v^{3(\Re{\tt s} + \frac{1}{2} \Re {\tt u}  - \frac{1}{2})} V^{-\frac{1}{4}B} Z^{-\Re({\tt w} + \frac{1}{2} {\tt u})}.$$ 
 Summing this over  $Z = 2^{\nu}$, we get a convergent sum by the last condition in \eqref{region}, and choosing $V > v^{3/4}$, we obtain a total bound of
 \begin{equation}\label{E2}
 v^{3(\Re{\tt s} + \frac{1}{2} \Re {\tt u}  - \frac{1}{2}) - \frac{3}{16}B}.
 \end{equation} 
 
 We proceed to treat $(\mathscr{V}_{{\tt s}, {\tt w}}^{+} \widehat{H})({\tt u})$, which similar, but slightly simpler.  Let initially $({\tt u}, {\tt s}, {\tt w})$ be in the region \eqref{cond}.  We straighten the   $\xi$-contour and shift it to the far left to $\Re \xi = -A_1 \not\in \Bbb{N}$ satisfying \eqref{a1}. By \eqref{resi} we do not pick up any poles. As before this gives analytic continuation to the region \eqref{region} and meromorphic continuation to the region \eqref{mero} with poles at most at ${\tt u}/2 + {\tt w} = 0$ and $ {\tt u}/2 + {\tt s} + \mu_j = 0$, $j = 1, 2, 3$. In the present case, we do not shift back to the right, but estimate the integral on the line $\Re \xi = -A_1$. We continue to write $v = \frac{1}{2} \Im {\tt u} > 0$ and assume that $v$ is sufficiently large. We focus now on the first term in \eqref{trafo}, since the second is exponentially decreasing in $v$. Again we split the $\xi$-integral into two pieces and consider 
  $$I^+_1 := \int_{(-A_1)}  \widehat{H}(3{\tt s} - {\tt w} - 1 + {\tt u} + 2\xi) G^+(\xi) \mathcal{G}_{-\mu}^-(1 - {\tt s} - \textstyle\frac{1}{2} {\tt u} - \xi)\mathcal{G}_{\mu}^{-}({\tt s} + \frac{1}{2} {\tt u})  \displaystyle w\left(\frac{v+\Im \xi}{V}\right)\frac{d\xi}{2\pi i}$$
and 
$$I^+_2 := \int_{(-A_1)}  \widehat{H}(3{\tt s} - {\tt w} - 1 + {\tt u} + 2\xi) G^+(\xi) \mathcal{G}_{-\mu}^-(1 - {\tt s} - \textstyle\frac{1}{2} {\tt u} - \xi)\mathcal{G}_{\mu}^{-}({\tt s} + \frac{1}{2} {\tt u}) \displaystyle\left(1 - w\left(\frac{v+\Im \xi}{V}\right)\right) \frac{d\xi}{2\pi i},$$
that is, we distinguish between $\Im ({\tt u} + 2\xi)$ small or large. As before we assume $v^{\delta} < V < v^{1-\delta}$. We estimate $I_1^+$ trivially getting
\begin{equation}\label{err1}
I_1^+ \ll v^{3\Re({\tt s} + \frac{1}{2} {\tt u} - \frac{1}{2})} v^{-A_1 - \frac{1}{2}} V^{\frac{3}{2} - \frac{1}{2}\Re{\tt u} - \Re{\tt w}  +A_1} . 
\end{equation}
For the estimation of $I_2^+$ we can restrict ourselves to the branch  $v + \Im \xi < -V$, for the other branch satisfies \eqref{exp-sr}. Again we now use the asymptotic formula \eqref{asymp-j} which contributes an error term  of at most
\begin{equation}\label{err2}
\begin{split}
&\int_{-\infty}^{-V-v} |v+x|^{-\Re( \frac{1}{2}{\tt u} + {\tt w}) - \frac{1}{2} + A_1-B} |x|^{-A_1 - \frac{1}{2}} v^{3\Re({\tt s} + \frac{1}{2}{\tt u}) - \frac{3}{2}} dx\\
& \ll  v^{\Re (3{\tt s} -   {\tt u}  - {\tt w}) - \frac{3}{2} - B} + v^{3\Re( {\tt s} + \frac{1}{2} {\tt u}) - \frac{3}{2} - A_1 - \frac{1}{2}} V^{-\Re( \frac{1}{2}{\tt u} + {\tt w}) + \frac{1}{2} + A_1 - B}.
\end{split}
\end{equation}
The main term produces an integral of the shape
\begin{displaymath}
\begin{split}
\mathcal{G}_{\mu}^-({\tt s} + \textstyle \frac{1}{2} {\tt u})&\displaystyle \int_{-\infty}^{0}  (x + v  )^{-\frac{1}{2}  + A_1 - \Re( {\tt w} + \frac{1}{2}   {\tt u})} x^{-A_1 - \frac{1}{2}} \omega(x) e^{i\phi(x)} 
 \left(1 - w\left(  \frac{v+x}{V}\right)\right) dx
\end{split}
\end{displaymath}
 with the same phase function \eqref{phi} as before and $|x+v|^j \omega^{(j)}(x) \ll_j 1$. We restrict the integral to smooth dyadic ranges $|x+v| \asymp Z \geq V$, so $|x| \asymp Z + v$. 
 By the same application of Lemma \ref{integrationbyparts}   each such dyadic region can be bounded by
$$Z^{\frac{1}{2}  + A_1 - \Re( {\tt w} + \frac{1}{2}   {\tt u})} (Z+v)^{-A_1 - \frac{1}{2}}  v^{3\Re({\tt s} + \frac{1}{2} {\tt u}) - \frac{3}{2}} \Bigl(\sqrt{\frac{Zv}{v+Z}}\Bigr)^{-B/2}\ll v^{3\Re({\tt s} + \frac{1}{2} {\tt u}) - \frac{3}{2}}  Z^{- \Re( {\tt w} + \frac{1}{2}   {\tt u})} V^{-\frac{1}{4}B},$$
 and summing this over $Z = 2^{\nu}$, we obtain a total bound of 
 \begin{equation}\label{err3}
   v^{3\Re({\tt s} + \frac{1}{2} {\tt u}) - \frac{3}{2}}    V^{-\frac{1}{4}B}.
 \end{equation}
We now choose $V = v^{1/2}$, $A_1 = \frac{1}{2} (A + \Re( s -  w) - 2)$ in agreement with \eqref{a1}. Then the bounds \eqref{err1}, \eqref{err2}, \eqref{err3} become $O(v^E)$ with
\begin{displaymath}
\begin{split}
E = \max& \Big(\textstyle \frac{1}{4} \Re(11{\tt s} + 5{\tt u} - {\tt w} - 3) - \frac{1}{4} A , \,\,\Re (3{\tt s} -    {\tt u}  - {\tt w}) - \frac{3}{2} - B, \\
& \textstyle\frac{1}{4}(\Re(11{\tt s} + 5{\tt u} - {\tt w} - 5) - \frac{1}{4}A - \frac{1}{2}B,    \,\,   3\Re({\tt s} + \frac{1}{2} {\tt u}) - \frac{3}{2} - \frac{1}{8}B \Big).
\end{split}
\end{displaymath}
 Combining this with the bounds \eqref{E1} and \eqref{E2}, we complete the proof of \eqref{boundJ}. 
 
\section{A preliminary reciprocity formula}\label{prelim-formula}

As outlined in the introduction, Theorem \ref{thm1} is a consequence of the five-step procedure Kuznetsov-Voronoi-reciprocity-Voronoi-Kuznetsov. In the following proposition we consider the middle triplet Voronoi-reciprocity-Voronoi. 

For $s, w \in \Bbb{C}$ with $\Re s, \Re w > 3/2$, $q, \ell \in \Bbb{N}$  coprime and a function $H$ satisfying $H(x) \ll x^{2/3}$, we define the absolutely convergent expression 
$$\mathcal{E}^{\pm}_{q, \ell}(s, w; H) := \sum_{\ell \mid r} \sum_{q \mid n_1c} \sum_{n_2} \frac{A(n_2, n_1)S(\pm r, n_2, c)}{n_2^{s  }n_1^{2s} c r^{w}}H\left(\frac{\sqrt{rn_2}}{c}\right). $$

\begin{prop}\label{prelim} Let $H$ be an admissible function of type $(500, 500)$.  Let $3/2 < \Re s < 2$, $4 < \Re w < 5$ and suppose that $q, \ell \in \Bbb{N}$ are coprime. Then
\begin{equation}\label{rec-prelim}
\mathcal{E}_{q, \ell}^{+}(s, w; H) =  \mathcal{N}^{(1)}_{q, \ell}(s, w; H)  - \mathcal{N}^{(2)}_{q, \ell}(s, w; H) + \sum_{\pm} \mathcal{E}_{\ell, q}^{\pm}\Big(s', w'; \widecheck{\big(\mathscr{V}^{\pm}_{s', w'}\widehat{H}\big)}\Big) \end{equation}
where $s', w'$ are as in \eqref{new}, the right hand side employs the notation \eqref{mellin1}, \eqref{mellin2} and \eqref{trafo},  and the ``main terms'' $\mathcal{N}^{(1)}_{q, \ell}(s, w; H)$ and $ \mathcal{N}^{(2)}_{q, \ell}(s, w; H)$ are given in \eqref{res1} and \eqref{N2}; they vanish if $F$ is cuspidal, and if $F=\E_0$,  they have meromorphic continuation to an $\varepsilon$-neighbourhood of $\Re w \geq \Re s \geq 1/2$ and satisfy the bounds
 \begin{equation}\label{boundsN}
 \begin{split}
& \big((s - \textstyle\frac{1}{2})(w-s)(s+w-1)\big)^4\mathcal{N}^{(1)}_{q, \ell}(s, w; H) \ll_{s, w} (q\ell)^{-1+\varepsilon}, \\
 &   \big((s - \textstyle\frac{1}{2})(w-s)(s+w-1)\big)^6\mathcal{N}^{(2)}_{q, \ell}(s, w; H) \ll_{s, w} (\ell q)^{\varepsilon} \ell^{-1}
 \end{split}
 \end{equation}
Moreover, the function $\mathscr{V}^{\pm}_{s', w'}\widehat{H}$  is holomorphic in  
\begin{equation}\label{prop-1}
0 < \Re s < 2, \quad 0 < \Re w < 5, \quad \max(2\theta - 2 \Re s', -2\Re w') < \Re u  < 15
\end{equation} and in this region bounded by 
\begin{equation}\label{prop-2}
   \mathscr{V}^{\pm}_{s', w'}\widehat{H}(u) \ll_{s, w} (1 + |u|)^{-15}.
 \end{equation}  
In addition,  it is meromorphic in 
\begin{equation}\label{mero2}
0 < \Re s < 2, \quad 0 < \Re w < 5, \quad \max(-2+2\theta - 2 \Re s', -2-2\Re w') < \Re u  < 15
\end{equation}
with poles at most at $u \in \{-2w', -2s' - 2\mu_1, -2s' - 2\mu_2, -2s' - 2\mu_3\}$. 

  \end{prop} 

\textbf{Proof.}  We first note that if $H$ is admissible  of type $(A, B)$, then   by Lemma \ref{admis1}a) and \eqref{mellin}  we have  $H(x) \ll x^{2/3}$, so that $\mathcal{E}^+_{q, \ell}(s, w; H)$ makes sense.  By Mellin inversion and \eqref{alt} we have 
$$ \mathcal{E}^+_{q, \ell} (s, w; H)  = \int_{(-1 )} \widehat{H}(u) \mathcal{D}^+_{q, \ell}(s, u, w) \frac{du}{2\pi i}.$$
By \eqref{boundD} and \eqref{assumption-c} the $u$-integral is absolutely convergent provided that $$ \Re u < 0, \quad \Re(3 s - u/2) > 3, \quad \Re(3 s + u/2) > 3/2,$$
which is automatically satisfied, provided \eqref{variables} holds. We shift the $u$-contour to $\Re (s + u/2) = -2/3$ (for the moment any negative number would suffice, but later we need $\Re (s + u/2) < -1/2$). This is still in agreement with   \eqref{variables}. On the way we may pick up a pole at $u = 2 - 2s$, which contributes
\begin{equation}\label{res1}
\mathcal{N}^{(1)}_{q, \ell}(s, w; H) := \underset{u = 2-2s}{\text{res}} \widehat{H}(u) \mathcal{D}^+_{q, \ell}(s, u, w) =\begin{cases} \sum_{j = 1}^3 \frac{1}{(j-1)!}\mathcal{R}_{q, \ell;j}(s, w) \widehat{H}^{(j-1)}(2-2s), & F = \E_0,\\
0, & F \text{ cuspidal}\end{cases}
\end{equation}
with the notation as in \eqref{resR}. Lemma \ref{pole1} provides analytic continuation of this term as well as  the first bound in \eqref{boundsN}.

Having shifted the contour to  $\Re (s + u/2) = -2/3$,  we can insert the first functional equation in \eqref{Dfunc}  and apply the definition \eqref{defPhi} since we are in the region of absolute convergence.  In this way we conclude
\begin{equation}\label{int}
\begin{split}
\mathcal{E}_{q, \ell}^{+}(s, w; H) =  \mathcal{N}^{(1)}_{q, \ell}(s, w; H) + & \int_{(-\frac{4}{3}- 2\Re s)} \widehat{H}(u) \sum_{\pm} \mathcal{G}_{-\mu}^{\mp}(1 - s - u/2)  \\
& \sum_{(c,  d) = 1} \sum_{\substack{\ell \mid md \\ q \mid mc}}   \sum_n \frac{A(m, n) e(\pm n \bar{d}/c)}{c^{3s+u/2 -1} m^{s+w} d^{w+u/2} n^{1 - s - u/2}}  \frac{du}{2\pi i}.
\end{split}
\end{equation}
At this point we insert artificially a factor
$$1 = e\left(\mp \frac{n}{cd}\right) e\left(\pm \frac{n}{cd}\right) = e\left(\mp \frac{n}{cd}\right) \int_{\mathcal{C}}  G^{\pm}(\xi)\left(\frac{n}{cd}\right)^{-\xi} \frac{d\xi}{2\pi i},$$
where $G^{\pm}$ and  $\mathcal{C}$ were defined in \eqref{Gpm} resp.\ \eqref{contour}. 
The contour is designed so that  the integral is absolutely convergent, but the contour is to the right of the pole at $\xi = 0$. 
By the reciprocity formula \eqref{rec} the integral in \eqref{int} equals
\begin{displaymath}
\begin{split}
  \int_{(-\frac{4}{3}-2\Re s)}  \int_{ \mathcal{C}}\widehat{H}(u)\sum_{\pm} G^{\pm}(\xi) \mathcal{G}_{-\mu}^{\mp}(1 - s - u/2)   \sum_{(c,  d) =1}  \sum_{\substack{q \mid mc \\ \ell \mid md}}   \frac{\Phi(d, \mp c, m; 1 - s - u/2 + \xi) }{c^{3s+u/2 -1} m^{s+w}d^{w+u/2} (cd)^{-\xi} }   \frac{d\xi}{2\pi i} \frac{du}{2\pi i},
\end{split}
\end{displaymath}
and the entire expression is still absolutely convergent. Here we used that $\Re(1 - s - u/2 + \xi) \geq 16/15 > 1$ on the entire $\xi$-contour, so that we are in the region of absolute convergence of $\Phi$.  It is now convenient to interchange the $u$- and $\xi$-integration and to replace the straight $u$-contour  with a polygonal contour $\mathcal{C}(\xi)$ (depending on $\xi$) such that $\Re(1 - s - u/2 + \xi) = 16/15$. The resulting expression is still absolutely convergent.  We now introduce the new variables
\begin{equation}\label{new-var}
s' = \textstyle\frac{1}{2}(1 - s + w), \quad w' = \frac{1}{2}(3s+w-1),\quad u' = 3s-w -1 + u - 2\xi. 
\end{equation}
This has the following effect: the exponents $(3s + \frac{1}{2}u - 1- \xi, s+w, w+\frac{1}{2}u - \xi)$ of $(c, m, d)$ become $(w'+\frac{1}{2}u', s'+w', 3s' + \frac{1}{2}u' - 1)$, and the contour $C(\xi)$ given by $\Re(1 - s - u/2 + \xi) = 16/15$ becomes $s' + \frac{1}{2} u' = -\frac{1}{15}$ (independently of $\xi$). 
Recalling the definition \eqref{deftildeD}, we can recast the integral in \eqref{int} as 
\begin{displaymath}
\begin{split}
  \int_{ \mathcal{C}}  \int_{( - \frac{2}{15} - 2\Re s')} \widehat{H}(3s' - w' - 1 + u' + 2\xi) \sum_{\pm} G^{\pm}(\xi) \mathcal{G}_{-\mu}^{\mp}(1 - s' - \textstyle\frac{1}{2}u'-\displaystyle\xi) \widetilde{\mathcal{D}}^{\mp}_{\ell, q}(s', u', w')  \frac{du'}{2\pi i} \frac{d\xi}{2\pi i}.
\end{split}
\end{displaymath}
Now we shift the $u'$-integral the right to $\Re (s' + u'/2) = 1/2$. Again we may pick up a pole with residue
$$\mathcal{N}^{(2)}_{q, \ell}(s, w; H) := \int_{\mathcal{C}}\underset{u'= -2s'}{\text{res}} \widehat{H}(3s' - w' - 1 + u' + 2\xi) \sum_{\pm} G^{\pm}(\xi) \mathcal{G}_{-\mu}^{\mp}(1 - s' - \textstyle\frac{1}{2}u'-\displaystyle\xi) \widetilde{\mathcal{D}}^{\mp}_{\ell, q}(s', u', w')   \frac{d\xi}{2\pi i} $$
(which counts negative because of the right shift). 
 If $F$ is cuspidal, this vanishes, and if $F = \E_0$
, then by \eqref{tildeD} it equals
\begin{equation}\label{N2}
 \sum_{j=1}^3  \tilde{\mathcal{R}}_{\ell, q;j}(s', w') \sum_{\nu_1 + \nu_2 = j-1} \left( - \frac{1}{2}\right)^{\nu_2}  \sum_{\pm} \int_{(1/10)}   G^{\pm}(\xi)    \widehat{H}^{(\nu_1)}(-2s + 2\xi) (\mathcal{G}_{(0, 0, 0)}^{\mp})^{(\nu_2)}(1 -\displaystyle\xi)   \frac{d\xi}{2\pi i},
 \end{equation}
where we straightened  the contour $\mathcal{C}$. Note that the conditions $\Re w \geq \Re s \geq 1/2$ and $\Re w' \geq \Re s' \geq 1/2$ are equivalent and $(s' - 1/2)(w'-s')(s' +w' - 1) = (s-1/2)(w-s)(s+w-1)$, so Lemma \ref{pole2} and Lemma \ref{sec-prop} with $x = 1/10$, $c = -2s$, $d=1$   provide meromorphic continuation of this term as well as the second bound in \eqref{boundsN}. (Note that $\ell$ and $q$ are interchanged relative to Lemma \ref{pole2} and the $\xi$-integral contributes at most a triple pole at $s = 1/2$.)

Having shifted the contour to $\Re(s' + u'/2) = 1/2$, we apply the second functional equation in \eqref{Dfunc} and recall the definition \eqref{trafo} to obtain
$$\mathcal{E}_{q, \ell}^{+}(s, w; H) =  \mathcal{N}^{(1)}_{q, \ell}(s, w; H)  - \mathcal{N}^{(2)}_{q, \ell}(s, w; H) +\sum_{\pm} \int_{(1 - 2\Re s')}  (\mathscr{V}^{\pm}_{s', w'}\widehat{H})(u') \mathcal{D}_{\ell, q}^{\pm}(s', u', w') \frac{du'}{2\pi i}.$$ 
Lemma \ref{lem9}  gives us analytic continuation and decay conditions for  $\mathscr{V}^{\pm}_{s', w'}\widehat{H}$, so that we can continue to shift the contour to the right into the region of absolute convergence of $ \mathcal{D}_{\ell, q}^{\pm}(s', u', w')$ to, say, $\Re(s' + u'/2) = 3/2$. We write $ \mathcal{D}_{\ell, q}^{\pm}(s', u', w') $ in terms of its Dirichlet series representation  \eqref{alt}, and by Mellin inversion we obtain the formula \eqref{rec-prelim}. 

Now we observe that for $0 < \Re s < 2$, $0 < \Re w < 5$ we have $ -1/2 < \Re s', \Re w' < 5$. 
By \eqref{region} -- \eqref{boundJ} with $A, B \geq 500$, $|\Re {\tt s}|, |\Re{\tt w}| \leq 5$,  
   we obtain the holomorphicity of $\mathscr{V}_{s', w'}^{\pm}\widehat{H}$ in the region \eqref{prop-1} and the bound \eqref{prop-2}. In particular, if $3/2 < \Re s < 2$ and $4 < \Re w < 5$, then $\Re s' \geq 3/2$ and $\Re w' \geq 7/2$, so that \eqref{prop-1} and \eqref{prop-2} in combination with \eqref{mellin} imply that $\widecheck{\big(\mathscr{V}^{\pm}_{s', w'}\widehat{H}\big)}(x) \ll x^{2/3}$, so that the rightmost term in \eqref{rec-prelim} makes sense.   Meromorphic continuation of $\mathscr{V}^{\pm}_{s', w'}\widehat{H}(u)$ to \eqref{mero2} and the location of poles follows from the statement containing \eqref{mero}. This completes the proof. \\

We end this section by relating $\mathcal{E}^{\pm}(s, w; H)$ to spectral sums. By \eqref{kuz-all} we have for $\Re s, \Re w > 3/2$, $q, \ell$ coprime and $\mathfrak{h}$ admissible that 
\begin{equation}\label{formula1}
\sum_{d_1d_2 = q} \sum_{ \ell \mid r} \sum_{(n_1, q) = d_1} \sum_{n_2} \frac{A(n_2, n_1)}{n_2^sn_1^{2s}r^w} \mathcal{A}_{d_2}(r, n_2, \mathfrak{h}) = \sum_{\ell \mid n_2}\sum_{n_1} \frac{A(n_2, n_1)}{n_2^{s+w} n_1^{2s}} \mathscr{N}\mathfrak{h} + \mathcal{E}_{q, \ell}^+(s, w; \mathscr{K}^*\mathfrak{h}).
\end{equation}
Conversely, if ${\tt H}$ satisfies $x^j{\tt H}^{(j)}(x) \ll \min(x, x^{-3/2})$ for $0 \leq j \leq 3$ and if in addition $\widehat{{\tt H}}$ is holomorphic in $-2\vartheta-\varepsilon  <\Re u < 5$ and satisfies $\widehat{{\tt H}}(u) \ll (1 + |u|)^{-5}$, say, then  by 
  \eqref{kuz2} for $\Re s', \Re w' > 3/2$ we have
\begin{equation}\label{formula2}
\mathcal{E}_{\ell, q}^{\pm}(s', w';  {\tt H}) = \sum_{d_1d_2 = \ell} \sum_{q \mid r} \sum_{(n_1, q) = d_1} \sum_{n_2} \frac{A(n_2, n_1)}{n_2^{s'}n_1^{2s'}r^{w'}} \mathcal{A}_{d_2}(\pm r, n_2, \mathscr{L}^{\pm}{\tt H}),  
\end{equation}
where by \eqref{rho-bound}, \eqref{rho-cusp-bound}, \eqref{gl3-RS}, Lemma \ref{final-decay}a and Weyl's law, the various sums in \eqref{formula1} and \eqref{formula2} are absolutely convergent. The next section is devoted to relating the left hand side of \eqref{formula1} and the right hand side of \eqref{formula2} to $L$-functions.

\section{Local factors}\label{local}

\subsection{Local computations} For a prime $p$ let $\alpha_{f, \nu}(p)$ ($\nu = 1, 2$), $\alpha_{F, j}(p)$ ($j = 1, 2, 3$) denote the Satake parameters of $f$ and $F$ at   $p$ satisfying
\begin{equation}\label{satake3}
 \alpha_{F, 1}(p) \alpha_{F, 2}(p) \alpha_{F, 3}(p)
 = 1, \quad \alpha_{f, 1}(p) \alpha_{f, 2}(p) = \begin{cases} 1, & p \nmid \text{cond}(f).\\ 0, & p \mid \text{cond}(f). \end{cases}
 \end{equation}
We have
\begin{equation*}
\lambda_f(p^{\nu}) = \frac{\alpha_{f, 1}(p)^{\nu+1}  - \alpha_{f, 2}(p)^{\nu+1} }{\alpha_{f, 1}(p) - \alpha_{f, 2}(p)}
\end{equation*}
and 
\begin{equation}\label{satake2}
 A(p^{\nu}, p^{\mu}) =  \det\left(\begin{matrix} \alpha_{F, 1}(p)^{\nu+\mu+2} & \alpha_{F, 2}(p)^{\nu+\mu+2} & \alpha_{F, 3}(p)^{\nu+\mu+2}\\ \alpha_{F, 1}(p)^{\mu+1} & \alpha_{F, 2}(p)^{\mu+1} & \alpha_{F, 3}(p)^{\mu+1}\\1 & 1 & 1\end{matrix}\right) V_F(p)^{-1}
\end{equation}
where 
$V_F(p)= \det\left(\begin{matrix} \alpha_{F, 1}(p)^{ 2} & \alpha_{F, 2}(p)^{ 2} & \alpha_{F, 3}(p)^{2}\\ \alpha_{F, 1}(p) & \alpha_{F, 2}(p) & \alpha_{F, 3}(p)\\1 & 1 & 1\end{matrix}\right); 
$
 in particular
\begin{equation}\label{satake4}
\begin{split}
& \lambda_f(p) = \alpha_{f, 1}(p) + \alpha_{f,2}(p), \quad A(p, 1) =  \alpha_{F, 1}(p)+ \alpha_{F, 2}(p)+ \alpha_{F, 3}(p), \\
 &  A(1, p) =  \alpha_{F, 1}(p) \alpha_{F, 2}(p) + \alpha_{F, 1}(p) \alpha_{F, 3}(p)  + \alpha_{F, 2}(p) \alpha_{F, 3}(p).  
 \end{split}
 \end{equation}
 Note that \eqref{satake2} remains formally true for $\nu$ or $\mu = -1$ if we define $A(p^{\nu}, p^{\mu}) = 0$ in this case. 
Let $$L_p(p^{-s}, f \times F) 
= \prod_{j=1}^3\prod_{\nu=1}^2 \left(1 - \frac{\alpha_{F, j}(p) \alpha_{f, \nu}(p)}{p^s}\right)^{-1}$$
denote the local Rankin-Selberg factor at $p$ and correspondingly
\begin{equation*}
L(s, f \times F) := \prod_p L_p(p^{-s}, f \times F)
\end{equation*}
in $\Re s > 1$ (this may differ from the corresponding Langlands $L$-function by finitely many Euler factors). 
We start with the following combinatorial lemma. 

\begin{lemma}\label{eulerlemma} For $p$ prime and $\Re s> \theta + \vartheta$, the following   identities hold:
\begin{displaymath}
\begin{split}
&\sum_{\nu \in \Bbb{N}_0} \frac{A(p^{\nu}, 1)\lambda_f(p^{\nu})}{p^{\nu s}} = \begin{cases} \left(1 - \frac{A(1, p)}{p^{2s}} + \frac{\lambda_f(p)}{p^{3s}}\right) L_p(p^{-s}, f \times F), & p \nmid \text{{\rm cond}}(f),\\ L_p(p^{-s}, f \times F), &p \mid \text{{\rm cond}}(f),\end{cases} \\
&\sum_{\nu \in \Bbb{N}_0} \frac{A(p^{\nu+1}, 1)\lambda_f(p^{\nu})}{p^{\nu s}} = \left(A(p, 1) - \frac{A(1, p)\lambda_f(p)}{p^{s}} + \frac{\lambda_f(p^2)}{p^{2s}}\right)L_p(p^{-s}, f \times F),\\
&\sum_{\nu_1, \nu_2 \in \Bbb{N}_0} \frac{A(p^{\nu_2}, p^{\nu_1+1}) \lambda_f(p^{\nu_2})}{p^{s(\nu_2 + 2\nu_1)}}  = \left(A(1, p) - \frac{\lambda_f(p)}{p^s}\right) L_p(p^{-s},  f \times F), \quad p \nmid \text{{\rm cond}}(f).\\
\end{split}
\end{displaymath}
\end{lemma}

\textbf{Proof.} A simple computation based on \eqref{satake2} and geometric series shows the power series identity
\begin{displaymath}
\begin{split}
\sum_{\nu \in \Bbb{N}_0} A(p^{\nu}, 1) \lambda_f(p^{\nu}) X^{\nu} &= \big[1 - ( \alpha_{F, 1}(p) \alpha_{F, 2}(p) +  \alpha_{F, 1}(p) \alpha_{F, 3}(p) +\alpha_{F, 2}(p) \alpha_{F, 3}(p))\alpha_{f, 1}(p) \alpha_{f, 2}(p) X^2  \\
& +  \alpha_{F, 1}(p) \alpha_{F, 2}(p) \alpha_{F, 3}(p)(\alpha_{f, 1}(p)+ \alpha_{f, 2}(p) )\alpha_{f, 1}(p) \alpha_{f, 2}(p) X^3\big] L_p(X, f \times F) 
\end{split}
\end{displaymath}
which by \eqref{satake3} and  \eqref{satake4} with $X = p^{-s}$  gives the first formula. The other two are proved in the same way.\\

\textbf{Remarks:} 1) Note that $\theta+\vartheta \leq 5/14 + 7/64  < 1/2$. This is useful numerical coincidence, but $f \times F$ is known to correspond to an automorphic representation on ${\rm GL}(6)$ \cite{KiSh}, therefore  such a relation   follows at every place from   general bounds towards the Ramanujan conjecture on ${\rm GL}(n)$ \cite{LRS}. \\
2) We conclude from the first formula that the Dirichlet series expansion of $L(s, f \times F)$ in $\Re s > 1$ is given by
\begin{equation}\label{L23}
\sum_{\substack{n, m\\ (\text{cond}(f) , m) = 1}} \frac{A(n, m) \lambda_f(n)}{n^s m^{2s}}.
\end{equation}
This holds also with $(F, \lambda_f)$ replaced with $(E_{t, \chi}, \lambda_{t, \chi})$.\\

We need a similar technical lemma that we apply later for the contribution of Eisenstein series. 

\begin{lemma}\label{euler2} Let $M, d, g_1, g_2, q \in \Bbb{N}$ with $g_1 \mid g_2$  and $d, M \mid q$. 
Then the series
$$\sum_{\substack{c, f, n \mid q^{\infty}\\ (c, g_1) = (n, g_2) = 1}} \frac{A(cfM, nd)
}{c^u f^v n^{u+v}} \prod_{p \mid q} \prod_{j=1}^3 \left(1 - \frac{\alpha_{F, j}(p) }{p^u}\right) \left(1 - \frac{\alpha_{F, j}(p) }{p^v}\right), $$
initially defined for $\Re u, \Re v > \theta$ as an absolutely convergent series, has a holomorphic extension to an $\varepsilon$-neighbourhood of $\Re u, \Re v \geq 0$, $\Re (u+v) \geq 1/2$ and is bounded by $O(q^{\varepsilon} (dM)^{\theta})$ in this region. 
\end{lemma}


\textbf{Proof.} The sum factorizes into a product of $p\mid q$, and it suffices to consider each factor separately. Let $m = v_p(M)$, $k=v_p(d)$, and put $X =  p^{-u}$, $Y =   p^{-v}$. 
The local $p$-factor  equals
$$E(p) := \sum_{\substack{\beta, \gamma, \delta \in \Bbb{N}_0\\ (p^{\beta}, g_1) = (p^{\delta}, g_2) = 1}} A( p^{\beta+\gamma + m}, p^{\delta + k} )X^{\beta} Y^{\gamma} (XY)^{\delta} \prod_{j=1}^3 (1 - \alpha_{F, j}(p) X)(1 - \alpha_{F, j}(p) Y).$$
By \eqref{hecke} we have
 $$A( p^{\beta+\gamma + m}, p^{\delta + k} ) = A( p^{\beta+\gamma + m}, 1) \overline{A(p^{\delta + k}, 1 )} -  A( p^{\beta+\gamma + m-1}, 1) \overline{A(p^{\delta + k-1}, 1 )}$$
 (with the above convention $A(p^{-1}, 1) = 0$). We treat the first summand on the right hand side, the second one is similar. Depending on whether (i) $p \mid g_1$, (ii)  $p \mid g_2$, but $p \nmid g_1$ or (iii) $p \nmid g_2$, we need to compute one, two or three geometric series. In  cases (i) and (ii) we obtain  
 $$E(p) = \overline{{A}(p^{k}, 1) }\sum_{i= 0}^2 \sum_{j=0}^3  \frac{P_{i, j}(\alpha_{F, 1}(p),\alpha_{F, 2}(p), \alpha_{F, 3}(p))}{ V_F(p)}  
 X^i Y^j$$
where $P_{i, j}$ is a homogeneous polynomial of degree $3 +m +i+j$. 
Since    $E(p)$ is an entire function in $\alpha_{F, j}$ 
for $X, Y$ sufficiently small,  each $P_{i, j}$ must be divisible $V_F(p)$, 
and we obtain 
$$\overline{{A}(p^{k}, 1) } \sum_{i= 0}^2 \sum_{j=0}^3  \tilde{P}_{i, j}(\alpha_{F, 1},\alpha_{F, 2}, \alpha_{F, 3}) 
   X^i Y^j$$
with a homogeneous polynomial $\tilde{P}_{i, j}$ 
 of degree $m  +i+j$. In case (iii) a similar argument shows 
$$E(p) = \sum_{i= 0}^2 \sum_{j=0}^2  \frac{ \tilde{Q}_{i, j}(\alpha_{F, 1}(p),\alpha_{F, 2}(p), \alpha_{F, 3}(p))  \tilde{R}_{i, j}(\overline{\alpha_{F, 1}(p)},\overline{\alpha_{F, 2}(p)}, \overline{\alpha_{F, 3}(p)}) }{(1 - \overline{\alpha_{F, 1}(p)}XY)(1 - \overline{\alpha_{F, 2}(p)}XY)(1 - \overline{\alpha_{F, 3}(p)}XY)}  
 X^i Y^j$$
 for homogeneous polynomials $\tilde{Q}_{i, j}$,  $\tilde{R}_{i, j}$
 of degrees $m  +i+j$ and $k+i+j$ respectively. 

 Since $\max_j |\alpha_{F, j}(p)| 
 \leq p^{\theta}$, we can continue each $E(p)$ to an $\varepsilon$-neighbourhood of $\Re u, \Re v \geq 0$, $\Re (u+v) \geq 1/2$ and bound it by  $O(p^{\theta(k+m) + \varepsilon})$ in this region. 
 This completes the proof.  \\ 

Finally, using the formula for Bump's double Dirichlet series \cite[Proposition 6.6.3]{Go}, we have 
\begin{equation}\label{bump-double}
  \sum_{\ell \mid n_2}\sum_{n_1} \frac{A(n_2, n_1)}{n_2^{s+w} n_1^{2s}} =  \frac{L(s+w, F) L(2s, \bar{F})}{\ell^{s+w}\zeta(3s+w)} \sum_{n_1, n_2 \mid \ell^{\infty}} \frac{A(\ell n_2, n_1)}{n_2^{s+w} n_1^{2s}}\Biggr(\sum_{n_1, n_2 \mid \ell^{\infty}} \frac{A( n_2, n_1)}{n_2^{s+w} n_1^{2s}}\Biggl)^{-1},
\end{equation}
which provides analytic continuation of the left hand side, initially defined in $\Re s, \Re w > 1/2$ to the region $\Re (3s + w) > 1$, $\Re s, \Re w > \theta$ (with polar divisors at most at $s = 1/2$, $s + w = 1$ if $F = \E_0$), and it also provides the bound 
\begin{equation}\label{bound-bump}
 O_{s, w}(\ell^{-\Re (s+ w)+ \theta + \varepsilon})
 \end{equation} in this region, 
 away from polar divisors.

\subsection{The cuspidal case} We start by considering
\begin{displaymath}
\begin{split}
&\sum_{d_1d_2 = q} \sum_{ \ell \mid r} \sum_{(n_1, q) = d_1} \sum_{n_2} \frac{A(n_2, n_1)}{n_2^sn_1^{2s}r^w} \mathcal{A}^{\text{Maa{\ss}}}_{d_2}(\pm r, n_2,h) =\sum_{d_0 \mid q} \sum_{f \in \mathcal{B}^{\ast}(d_0)}  \epsilon_f^{(1 \mp 1)/2} h(t_f) \mathcal{S}_{q, \ell}(s, w; f) \end{split}
\end{displaymath}
with
$$  \mathcal{S}_{q, \ell}(s, w; f) := \sum_{d_1d_2 = q}\sum_{ M \mid d_2/d_0}    \sum_{ \ell \mid r} \sum_{(n_1, q) = d_1} \sum_{n_2}     \frac{A(n_2, n_1) \rho_{f, M, d_2}( r) \overline{\rho_{f, M, d_2}(n_2)}}{n_2^sn_1^{2s}r^w}
$$
for $f \in \mathcal{B}^{\ast}(d_0)$. We insert \eqref{rho-cusp} and obtain
\begin{displaymath}
\begin{split}
& \frac{L(1, \text{Ad}^2f)}{\prod_{p \mid d_0}(1 -p^{-2})}  \mathcal{S}_{q, \ell}(s, w; f) \\
= &  \sum_{d_1d_2 = q}\sum_{\ell \mid r} \sum_{(n_1, q) = d_1}   \sum_{n_2}  \sum_{d_0M \mid d_2}  \sum_{\delta_1, \delta_2 \mid M}\frac{\xi_{f}(M, \delta_1)  \xi_{f}(M, \delta_2) }{d_2\nu(d_2)}   \frac{\delta_1\delta_2}{M} \frac{A(n_2, n_1)\lambda_{f}(r/\delta_1) \lambda_{f}(n_2/\delta_2)}{n_2^{s  }n_1^{2s} r^{w}}\\
 = & \frac{1}{q \ell^w } \sum_{d_1d_2 = q}  \sum_{Md_0 \mid d_2}  \sum_{\delta_1, \delta_2 \mid M}   \frac{ \xi_{f}(M, \delta_1)  \xi_{f}(M, \delta_2) \delta_1^{1-w}\delta_2^{1-s} d_1^{1-2s} }{M  \nu(d_2)   }  \sum_{  r} \sum_{(n_1, d_2) = 1}   \sum_{n_2}\frac{A(n_2\delta_2, n_1d_1)\lambda_{f}(\ell r) \lambda_{f}(n_2)}{n_2^{s  }n_1^{2s} r^{w}}. 
  \end{split}
 \end{displaymath}
We open  $\lambda_f(\ell r)$ using the Hecke relations and recognize the $r$-sum as $\Lambda_f(\ell; w) L(w, f)$ with the notation as in \eqref{Lambda}. We also recognize the $n_1, n_2$-sum as $L(s, f \times F)$ up to Euler factors at primes dividing $q$. Hence
\begin{displaymath}
\begin{split}
\mathcal{S}_{q, \ell}(s, w; f)  =   \frac{\phi(q)}{q^2} \frac{L(w, f) L(s, f \times F)}{ \ell^w L(1, \text{Ad}^2f)}  \Lambda_f(\ell; w) \tilde{L}_q(s, w, f \times F)
 \end{split}
 \end{displaymath}
with  
\begin{equation}\label{local-q-cusp}
 \begin{split}
\tilde{L}_q(s, w, f\times F) = & \frac{q}{\phi(q)}\prod_{p \mid q}  \prod_{j=1}^3 \prod_{\nu=1}^2 \left(1 - \frac{\alpha_{F, j}(p) \alpha_{f, \nu}(p)}{p^s}\right)  \sum_{d_1d_2 = q}  \sum_{Md_0 \mid d_2}  \sum_{\delta_1, \delta_2 \mid M}   \delta_1^{1-w}\delta_2^{1-s} d_1^{1-2s}    \\
&\frac{ \xi_{f}(M, \delta_1)  \xi_{f}(M, \delta_2)  }{M  \nu(d_2)   } \sum_{\substack{(n_1, d_2) = 1\\ n_1 \mid q^{\infty}}}   \sum_{n_2 \mid q^{\infty}}\frac{A(n_2\delta_2, n_1d_1)  \lambda_{f}(n_2)}{n_2^{s  }n_1^{2s} }\prod_{p \mid d_0}(1 -p^{-2}).
 \end{split}
 \end{equation}
 In particular, using the notation \eqref{moment} we have
 \begin{equation}\label{moment1}
   \sum_{d_1d_2 = q} \sum_{ \ell \mid r} \sum_{(n_1, q) = d_1} \sum_{n_2} \frac{A(n_2, n_1)}{n_2^sn_1^{2s}r^w} \mathcal{A}^{\text{Maa{\ss}}}_{d_2}(\pm r, n_2, h)  = \mathcal{M}_{q, \ell}^{\text{Maa{\ss}}, \pm}(s, w; (h, h^{\text{hol}})) 
 \end{equation}
for any $h^{\text{hol}}$.  Estimating trivially (using \eqref{indiv} and \eqref{xi-arithmetic}), we obtain
 $$\tilde{L}_q(s, w, f\times F) \ll q^{\varepsilon} \sum_{d_1d_2 = q}  \sum_{Md_0 \mid d_2}  \sum_{\delta_1, \delta_2 \mid M}  (\delta_1\delta_2)^{1/2} \left(\frac{M^2}{\delta_1\delta_2}\right)^{\vartheta} \frac{1}{M} (d_1\delta_2)^{\theta} \ll q^{\theta  + \varepsilon}$$ 
for $\Re s, \Re w \geq 1/2$.  This proves \eqref{local-bound} when $f$ is Maa{\ss}. 
For $f \in \mathcal{B}^{\ast}(q)$ we have $d_0 = q$, hence $d_2=q$, $d_1 = M = \delta_1 = \delta_2 = 1$, and the $n_1, n_2$-sum over powers of primes dividing $q$ in \eqref{local-q-cusp} can be computed explicitly using the first formula in Lemma \ref{eulerlemma}. In this way  one obtains $\tilde{L}_q(s, w, f\times F) = 1$ for Maa{\ss} newforms $f$ of  level $q$. With slightly more computational effort and the the other two formulae in Lemma \ref{eulerlemma} one shows 
\begin{equation}\label{extension}
\tilde{L}_q(1/2, 1/2,  f\times F)  =    \frac{\phi(d_0)}{qd_0}  \prod_{p \mid \frac{q}{d_0}} \left(\frac{2 + A(1, p) + A(p, 1)}{1 +  \lambda_f(p)p^{-1/2} + p^{-1}}-1\right)
\end{equation}
if $q$ is squarefree and $f$ is a newform of level $d_0 \mid q$. 

The same computation holds verbatim for $f \in \mathcal{B}^{\ast}_{\text{hol}}(d_0)$, which completes the proof Lemma \ref{lem2} in the cuspidal case. 



\subsection{The Eisenstein case} Similarly as in the previous subsection we consider
\begin{displaymath}
\begin{split}
&\sum_{d_1d_2 = q} \sum_{ \ell \mid r} \sum_{(n_1, q) = d_1} \sum_{n_2} \frac{A(n_2, n_1)}{n_2^sn_1^{2s}r^w} \mathcal{A}^{\text{Eis}}_{d_2}(\pm r, n_2, h) = \sum_{d_0^2 \mid  q} \sum_{\substack{\chi \text{ (mod } d_0)\\ \text{primitive}}} \int_{\Bbb{R}}    \mathcal{S}_{q, \ell}(s, w; (t, \chi)) h(t) \frac{dt}{2\pi}  \end{split}
\end{displaymath}
with
$$  \mathcal{S}_{q, \ell}(s, w; (t, \chi)) := \sum_{d_1d_2 = q} \sum_{d_0 \mid M_1 \mid d_0^{\infty}} \sum_{\substack{ (M_2, d_0) = 1 \\ d_0 M_1M_2 \mid d_2}} \sum_{\ell \mid r} \sum_{(n_1, q) = d_1}   \sum_{n_2}     \frac{A(n_2, n_1)\rho_{\chi, d_0M_1M_2, d_2}(r, t) \overline{\rho_{\chi, d_0M_1M_2, d_2}(n_2, t)}}{n_2^{s  }n_1^{2s} r^{w}}$$
for $t\in \Bbb{R}$ and $\chi$ a primitive Dirichlet character modulo $d_0$. We define 
\begin{equation}\label{lchitd2}
L(\chi, t, d_2) := \prod_{p \mid d_2} \left(\left(1 - \frac{\chi^2(p)}{p^{1 + 2it}}\right) \left(1 - \frac{\bar{\chi}^2(p)}{p^{1 - 2it}}\right)\right)^{-1},
\end{equation}
insert \eqref{rho-eis} and recast $  |L(1 + 2it, \chi^2)|^2\mathcal{S}_{q, \ell}(s, w; (t, \chi))$ as
\begin{displaymath}
\begin{split}
 \sum_{d_1d_2 = q} &\frac{L(\chi, t, d_2) }{d_2\nu(d_2)}  \sum_{d_0 \mid M_1 \mid d_0^{\infty}} \sum_{\substack{ (M_2, d_0) = 1 \\ d_0 M_1M_2 \mid d_2}} \sum_{\delta_1, \delta_2 \mid M_2} \frac{M_1\delta_1\delta_2 \mu(M_2/\delta_1) \mu(M_2/\delta_2)\bar{\chi}(\delta_1) \chi(\delta_2)   }
 { \tilde{\mathfrak{n}}(d_0M_1M_2)^2M_2}  \\
 &  \sum_{(n_1, q) = d_1}   \sum_{\substack{c_1,f_1\\ (c_1, \frac{d_2}{d_0M_1M_2}) = 1}}   \sum_{\substack{\ell \mid c_2f_2\\ (c_2, \frac{d_2}{d_0M_1M_2}) = 1}}    \frac{A(c_1f_1M_1\delta_1, n_1)\bar{\chi}(c_1f_2)  \chi(c_2f_1) (c_2f_2M_1\delta_2)^{it} (c_1f_1M_1\delta_1)^{-it}  }{(c_2/c_1)^{2it}(c_1f_1M_1\delta_1)^{s  }n_1^{2s} (c_2f_2 M_1 \delta_2)^{w}}.
\end{split}
\end{displaymath}
We slightly re-write this as
\begin{displaymath}
\begin{split}
\frac{1}{q} \sum_{d_1d_2 = q} &\frac{d_1^{1-2s}}{ \nu(d_2)}L(\chi, t, d_2) \sum_{d_0 \mid M_1 \mid d_0^{\infty}} \sum_{\substack{ (M_2, d_0) = 1 \\ d_0 M_1M_2 \mid d_2}} \sum_{\delta_1, \delta_2 \mid M_2} \frac{M_1^{1-s-w}\delta_1^{1-s-it}\delta_2^{1-w+it} \mu(M_2/\delta_1) \mu(M_2/\delta_2)   }{ \tilde{\mathfrak{n}}(d_0M_1M_2)^2M_2}  \\
 & \bar{\chi}(\delta_1) \chi(\delta_2) \sum_{(n_1, d_2) = 1}   \sum_{\substack{c_1,f_1\\ (c_1, \frac{d_2}{d_0M_1M_2}) = 1}}   \sum_{\substack{\ell \mid c_2f_2\\ (c_2, \frac{d_2}{d_0M_1M_2}) = 1}}    \frac{A(c_1f_1M_1\delta_1, n_1d_1)\bar{\chi}(c_1f_2)  \chi(c_2f_1)   }{(c_2f_1/c_1f_2)^{it}(c_1f_1 )^{s  }n_1^{2s} (c_2f_2  )^{w}}.
\end{split}
\end{displaymath}
Next we dissolve the condition $\ell \mid c_2f_2$ by using the formula
$$\sum_{x \mid  ab} F(a, b) = \sum_{r_1 r_2 r_3 \mid x} \mu(r_3)  \sum_{a, b} F(r_1 r_3 a, r_2 r_3 b)$$ 
for any function $F$ (as long as the $a, b$-sum is absolutely convergent), which can be seen by first summing over $b \equiv 0$ (mod $x/(x, a)$), then sorting the $a$-sum by the gcd of $x$ and $a$ and finally removing the coprimality condition by M\"obius inversion. In this way we see that the previous display equals
\begin{displaymath}
\begin{split}
\frac{1}{q \ell^w}& \sum_{\lambda_1\lambda_2\lambda_3 = \ell}  \frac{\mu(\lambda_3)\chi(\lambda_1 \overline{\lambda_2})}{\lambda_3^w(\lambda_1/\lambda_2)^{it}} \sum_{d_1d_2 = q}  \frac{d_1^{1-2s}}{ \nu(d_2)}L(\chi, t, d_2)\sum_{d_0 \mid M_1 \mid d_0^{\infty}} \sum_{\substack{ (M_2, d_0) = 1 \\ d_0 M_1M_2 \mid d_2}} \\ 
& \sum_{\delta_1, \delta_2 \mid M_2} \frac{M_1^{1-s-w}\delta_1^{1-s-it}\delta_2^{1-w+it} \mu(M_2/\delta_1) \mu(M_2/\delta_2)  \bar{\chi}(\delta_1) \chi(\delta_2)  }{ \tilde{\mathfrak{n}}(d_0M_1M_2)^2M_2}  \\
 &  \sum_{(n_1, d_2) = 1}   \sum_{\substack{c_1,f_1\\ (c_1, \frac{d_2}{d_0M_1M_2}) = 1}}   \sum_{\substack{c_2,f_2\\ (c_2, \frac{d_2}{d_0M_1M_2}) = 1 }}    \frac{A(c_1f_1M_1\delta_1, n_1d_1)\bar{\chi}(c_1f_2)  \chi(c_2f_1)   }{(c_2f_1/c_1f_2)^{it}(c_1f_1 )^{s  }n_1^{2s} (c_2f_2  )^{w}}.
\end{split}
\end{displaymath}
We recognize the $\lambda_1, \lambda_2, \lambda_3$-sum as $\Lambda_{(t, \chi)}(\ell; w)$, defined in \eqref{Lambda}, and the $c_1, c_2, f_1, f_2$-sum  as $L(w +  i t, \chi) L(w-it, \bar{\chi}) L(s + it, F \times \chi) L(s - it, F \times \bar{\chi}) $ up to Euler factors at primes dividing $q$. Thus,  by brute force, we write the previous display as 
$$\frac{\phi(q)}{ q^2\ell^w}\Lambda_{(t, \chi)}(\ell; w) L(w +  i t, \chi) L(w-it, \bar{\chi}) L(s + it, F \times \chi) L(s - it, F \times \bar{\chi})\tilde{L}_q(s, w; E_{t, \chi} \times F),$$
where $\tilde{L}_q(s, w; E_{t, \chi} \times F)$ is defined by 
 \begin{equation}\label{local-q-eis}
\begin{split}
  & \frac{q}{\phi(q)}\sum_{d_1d_2 = q}  \frac{d_1^{1-2s}}{ \nu(d_2)}L(\chi, t, d_2) \sum_{d_0 \mid M_1 \mid d_0^{\infty}} \sum_{\substack{ (M_2, d_0) = 1 \\ d_0 M_1M_2 \mid d_2}}  \sum_{\delta_1, \delta_2 \mid M_2} \frac{ \mu(M_2/\delta_1) \mu(M_2/\delta_2)  \bar{\chi}(\delta_1) \chi(\delta_2)  }{ \tilde{\mathfrak{n}}(d_0M_1M_2)^2M_2}  \\
& M_1^{1-s-w}\delta_1^{1-s-it}\delta_2^{1-w+it} \sum_{\substack{(n_1, d_2) = 1 \\ n_1 \mid q^{\infty}}}   \sum_{\substack{c_1 f_1\mid q^{\infty} \\ (c_1, \frac{d_2}{d_0M_1M_2}) = 1}}      \frac{A(c_1f_1M_1\delta_1, n_1d_1)\bar{\chi}(c_1 )  \chi(f_1)   }{(f_1/c_1)^{it}(c_1f_1 )^{s  }n_1^{2s} } \\
& \prod_{p \mid \frac{d_2}{d_0M_1M_2}} \left(1 - \frac{\chi(p)}{p^{w + it}}\right)  \left(1 - \frac{\bar{\chi}(p)}{p^{w - it}}\right)\prod_{p \mid q} \prod_{j=1}^3 \left(1 - \frac{\alpha_{F, j}(p) \chi(p)}{p^{s+it}}\right) \left(1 - \frac{\alpha_{F, j}(p) \bar{\chi}(p)}{p^{s-it}}\right) 
 \end{split}
\end{equation}
with the notation \eqref{lchitd2}.  In particular, we can write
\begin{equation}\label{moment2}
\sum_{d_1d_2 = q} \sum_{ \ell \mid r} \sum_{(n_1, q) = d_1} \sum_{n_2} \frac{A(n_2, n_1)}{n_2^sn_1^{2s}r^w} \mathcal{A}^{\text{Eis}}_{d_2}(\pm r, n_2, h)    = \mathcal{M}_{q, \ell}^{\text{Eis}}(s, w; (h, h^{\text{hol}}))
\end{equation} 
for any $h^{\text{hol}}$. For $\Re s, \Re w  \geq 1/2$, $t \in \Bbb{R}$ we estimate trivially (using \eqref{indiv})
\begin{displaymath}
\begin{split}
 \tilde{L}_q(s, w; E_{t, \chi} \times F)  \ll q^{\varepsilon} \sum_{d_1d_2 = q} \sum_{d_0 \mid M_1 \mid d_0^{\infty}} \sum_{ d_0 M_1M_2 \mid d_2}  \sum_{\delta_1, \delta_2 \mid M_2}  \frac{ (\delta_1\delta_2)^{1/2 }}{M_2}  (M_1 \delta_1 d_1)^{\theta}  \ll q^{\theta  + \varepsilon}
\end{split}
\end{displaymath}
confirming \eqref{local-bound} in the case of Eisenstein series.  For an application in the next section we will also need   analytic continuation of $\tilde{L}_q(s, w; E_{t, \text{triv}} \times F) $ to certain complex values of $t$ for the  trivial character $\chi = \text{triv}$ modulo 1, i.e.\ the constant function with value 1. 

\begin{lemma}\label{cor-eis} The functions $\tilde{L}_q(s, w; E_{\pm (1-s)/i, \text{{\rm triv}}} \times F) $, $\tilde{L}_q(s, w; E_{\pm (1-w)/i, \text{{\rm triv}}} \times F) $, initially defined in $\Re s, \Re w > 1$ as absolutely convergent series, have meromorphic continuation to an $\varepsilon$-neighbourhood of $\Re s, \Re w \geq 1/2$ with polar divisors at most at $s = 1/2$, $ w = 1/2$ and satisfy the bounds
\begin{equation}\label{bound-cor-eis}
 \begin{cases}(s - 1/2) \tilde{L}_q(s, w; E_{\pm (1-s)/i,\text{{\rm triv}}} \times F) \\(w - 1/2) \tilde{L}_q(s, w; E_{\pm (1-w)/i, \text{{\rm triv}}} \times F) \end{cases} \ll_{s, w} q^{\theta + \varepsilon }
 \end{equation}
for $1/2 -\varepsilon \leq \Re s, \Re w < 1$. 
\end{lemma}

\textbf{Proof.} The possible polar divisors at $s = 1/2$ or $w = 1/2$ come from $L(\text{triv}, t, d_2)$ at $t = \pm(1-s)/i$ or $\pm (1-w)/i$ defined in \eqref{lchitd2}. An application of  Lemma \ref{euler2}    shows that the rest can be continued holomorphically to an $\varepsilon$-neighbourhood of $\Re s, \Re w \geq 1/2$. For $it \in \{\pm (1-s), \pm(1-w)\}$ and $\delta_1, \delta_2 \mid M_2$ we have 
$|\delta_1^{1-s-it} \delta_2^{1-w+it}| \leq \max(M_2^{2 - 2\min(\Re s, \Re w)}, M_2^{|\Re w - \Re s|})$, so that \eqref{bound-cor-eis} follows from \eqref{local-q-eis} and 
$$q^{\varepsilon} \sum_{d_1d_2 = q} \sum_{d_0 \mid M_1 } \sum_{   d_0 M_1M_2 \mid d_2}  \sum_{\delta_1, \delta_2 \mid M_2} \frac{M_2^{2 - 2\min(\Re s, \Re w)} + M_2^{|\Re w - \Re s|}}{M_2}  (M_1 \delta_1d_1)^{\theta}   \ll  q^{\theta + \varepsilon }$$
for $1/2 -\varepsilon \leq \Re s, \Re w < 1$.  This completes the proof.\\

We end this section by combining \eqref{moment1} (and the corresponding formula for the holomorphic case) and \eqref{moment2} with \eqref{moment-together} to obtain
\begin{equation}\label{moment-comb}
 \sum_{d_1d_2 = q} \sum_{ \ell \mid r} \sum_{(n_1, q) = d_1} \sum_{n_2} \frac{A(n_2, n_1)}{n_2^sn_1^{2s}r^w} \mathcal{A}_{d_2}(\pm r, n_2, \mathfrak{h}) = \mathcal{M}^{\pm}_{q, \ell}(s, w; \mathfrak{h}). 
\end{equation}

\section{Proof of Theorem \ref{thm1}}\label{proof1}

We have now prepared the scene for a quick proof of Theorem \ref{thm1}.  Let initially be $3/2 < \Re s < 2$, $4 < \Re w < 5$ and let $\mathfrak{h} = (h, h^{\text{hol}})$ be admissible.  Then by definition, $\mathscr{K}^{\ast}\mathfrak{h}$ is admissible of type $(500, 500)$, and moreover $3/2 < \Re s' < 2$, $7/2 < \Re w' < 5$.  Combining \eqref{rec-prelim}, \eqref{formula1}, \eqref{formula2},  \eqref{moment-comb}, we obtain
\begin{equation}\label{weobtain}
\begin{split}
& \mathcal{M}^{+}_{q, \ell}(s, w; \mathfrak{h}) = \sum_{\ell \mid n_2}\sum_{n_1} \frac{A(n_2, n_1)}{n_2^{s+w} n_1^{2s}} \mathscr{N}\mathfrak{h} - \sum_{j \in \{1, 2\}} (-1)^j  \mathcal{N}^{(j)}_{q, \ell}(s, w; \mathscr{K}^*\mathfrak{h}) 
+ \sum_{\pm}  \mathcal{M}^{\pm}_{\ell, q}\big(s', w', \mathscr{T}^{\pm}_{s', w'} \mathfrak{h}\big)
 \end{split}
\end{equation}
where
\begin{equation}\label{final-trafo}
 \mathscr{T}_{s', w'}^{\pm} \mathfrak{h} :=   \mathscr{L}^{\pm} \widecheck{\big(\mathscr{V}^{\pm}_{s', w'}\widehat{\mathscr{K}^*\mathfrak{h}}\big)}
 \end{equation}
 using the notation \eqref{kast},  \eqref{trafo},  \eqref{Hback} and \eqref{mellin1}, \eqref{mellin2}. Note, however, that we cannot simply insert the various formulas into each other to compute the transform on the right of \eqref{final-trafo}, because there is a process of analytic continuation in Lemma \ref{lem9}. By \eqref{prop-1}, \eqref{prop-2} and \eqref{mellin} we conclude that \eqref{formula2} is indeed applicable under the current assumption $3/2 < \Re s' < 2$, $7/2 < \Re w' < 5$. 
 
 Write $\mathscr{T}_{s', w'}^{\pm} \mathfrak{h}  = ({\tt h}_{s', w', \pm}, {\tt h}^{\text{hol}}_{s', w', \pm})$. Combining \eqref{prop-1} -- \eqref{prop-2} with Lemma 
   \ref{final-decay}a, we conclude that $\mathscr{T}_{s', w'}^{\pm} \mathfrak{h}$ is weakly admissible as in \eqref{weakly} 
provided  that
$$\Re s' > \theta + \vartheta, \quad  \Re w' > \vartheta, \quad 0 < \Re s < 2, \quad 0 < \Re w < 5$$
(clearly ${\tt h}_{s', w', \pm}(t)$ is even in $t$).  Since $\theta + \vartheta < 1/2$, this includes in particular the region $1/2 \leq \Re s \leq \Re w < 3/4$. Moreover,  combining \eqref{mero2} with   Lemma 
   \ref{final-decay}b, we see that ${\tt h}_{s', w', \pm}(t)$ is meromorphic in an $\varepsilon$-neighbourhood of $|\Im t | \leq 1/2$ with   poles at most at $\pm it \in \{w', s' + \mu_1, s'+\mu_2, s' + \mu_3\}$. We will need this observation in a moment after having proved the next lemma.   

It remains to continue all terms in \eqref{weobtain} to a domain containing \eqref{final-region}. 
The analytic continuation of the cuspidal contribution of the terms $\mathcal{M}^{+}_{q, \ell}(s, w; \mathfrak{h}) $ and $\mathcal{M}^{\pm}_{\ell, q}(s', w', \mathscr{T}^{\pm} \mathfrak{h})$ 
is clear. 
For the Eisenstein contribution, we appeal to the following lemma. 
 \begin{lemma}\label{analyticcont}  Let $\mathfrak{h} = (h, h^{\text{{\rm hol}}})$ be weakly admissible, and suppose that $h$ has a meromorphic continuation to an $\varepsilon$-neighbourhood of $|\Im t | \leq 1/2$ with at most finitely many poles. Then the term $ \mathcal{M}^{\text{{\rm Eis}}}_{q, \ell}(s, w; \mathfrak{h}) $, initially defined in $\Re s, \Re w > 1$ continues meromorphically to an $\varepsilon$-neighbourhood of $\Re w, \Re s \geq 1/2$ with at most finitely many polar divisors.   If  $1/2 \leq \Re s, \Re w < 1$,  its analytic continuation is given by  
$ \mathcal{M}^{\text{{\rm Eis}}}_{q, \ell}(s, w; \mathfrak{h}) + R_{q, \ell}(s, w; \mathfrak{h})$ 
where $R_{q, \ell}(s, w; \mathfrak{h})$ is defined as
$$\frac{\phi(q)}{q^2} \sum \underset{\substack{t = \pm i(1-s)\\ t =  \pm i(1-w)}}{\text{{\rm res}}} (\pm i) \frac{L(s + it, F)L(s - it, F) \zeta(w + it) \zeta(w - it)}{\zeta(1 + 2it)\zeta(1 - 2it)} \tilde{L}_q(s, w, E_{t, \text{{\rm triv}}} \times F) \frac{\Lambda_{t, \text{{\rm triv}}}(\ell; w)}{\ell^w} h(t).$$ 
 \end{lemma}

\textbf{Proof.}  For $t \in \Bbb{R}$ choose   $0 < \sigma(t) < 1/4$ in a continuous way such that $L(1 - 2 \sigma + 2it, \chi) \not= 0$ for $0 \leq \sigma < 2\sigma(t)$ and all primitive Dirichlet characters $\chi$ of conductor $c_{\chi}$ such that $c_{\chi}^2 \mid q$ and in addition $h(t - i\sigma)$ is pole-free for $0 \leq \sigma < 2\sigma(t)$. Let initially $1 <  \Re s < 1 + \sigma(\Im s)$, $1 < \Re w  < 1 + \sigma(\Im w)$. In the defining integral of $ \mathcal{M}^{\text{{\rm Eis}}}_{q, \ell}(s, w; \mathfrak{h}) $ we shift the   $t$-contour down to $\Im t = -\sigma(\Re t)$. We pick up a pole at  $w-it = 1$ if $\chi = \text{triv}$ and a (triple) pole at $s - it = 1$ if in addition $F = \E_0$. The remaining integral is holomorphic in $1- \sigma(\Im s) <  \Re s < 1 + \sigma(\Im s)$, $1-  \sigma(\Im w) < \Re w  < 1 + \sigma(\Im w)$. Now choosing $s, w$ with $1- \sigma(\Im s) <  \Re s < 1  $, $1-  \sigma(\Im w) < \Re w  < 1  $, we shift the $t$-contour back to $\Im t = 0$ picking up a pole at $w + it = 1$  if $\chi = \text{triv}$ and a (triple) pole at $s + it = 1$ if in addition $F = \E_0$. This proves the desired formula for $1- \sigma(\Im s) <  \Re s < 1$, $1-  \sigma(\Im w) < \Re w  < 1$ , but then it follows for $s, w$ in an $\varepsilon$-neighbourhood of $\Re w, \Re s \geq 1/2$ by analytic continuation, using Lemma \ref{cor-eis} for the term $ \tilde{L}_q(s, w, E_{t, \text{{\rm triv}}} \times F) $. 
This completes the proof. \\

Note that this lemma is applicable both for the admissible function $\mathfrak{h}$ and for $\mathscr{T}_{s', w'}^{\pm} \mathfrak{h}$, since the latter satisfies the assumption of meromorphic continuation.  We recall that  condition \eqref{final-region} implies $1/2 \leq \Re s' \leq \Re w' < 1$. 
We conclude that the analytic continuation of $ \mathcal{M}^{\text{{\rm Eis}}}_{q, \ell}(s, w; \mathfrak{h}) $ and $\mathcal{M}^{\text{Eis}}_{\ell, q}\big(s', w', \mathscr{T}^{\pm}_{s', w'} \mathfrak{h}\big)$  infers two extra main terms, and we obtain the reciprocity formula \eqref{final-formula} with
\begin{equation}\label{main-term-final}
\begin{split}
 & \mathcal{N}_{q, \ell}(s, w; \mathfrak{h}) =\\
  &   \sum_{\ell \mid n_2}\sum_{n_1} \frac{A(n_2, n_1)}{n_2^{s+w} n_1^{2s}} \mathscr{N}\mathfrak{h} - \sum_{j \in \{0, 1\}}(-1)^j  \mathcal{N}^{(j)}_{q, \ell}(s, w; \mathscr{K}^*\mathfrak{h}) -  R_{q, \ell}(s, w; \mathfrak{h}) +  R_{\ell, q}(s', w';  \mathscr{T}_{s', w'}^{\pm} \mathfrak{h}).
  \end{split}
\end{equation}
A (meromorphic) continuation of  these terms follows from 
\eqref{bump-double} 
and Proposition \ref{prelim} and Lemma \ref{analyticcont}, and 
the bounds   \eqref{bound-cor-eis}, \eqref{bound-bump}, \eqref{boundsN} confirm \eqref{required-bound}, away from polar divisors, in an $\varepsilon$-neighbourhood of \eqref{final-region}. While the individual main terms may have polar lines, their joint contribution must be holomorphic, because the rest of terms in the reciprocity formula are holomorphic. Hence the various (possible) polar divisors 
must cancel, so that by standard complex analysis (e.g.\ Cauchy's integral formula in the $s$ and $w$ variable) the estimates \eqref{required-bound} are valid in the entire region \eqref{final-region}. This completes the proof.

\section{Proof of Theorem \ref{thm2}}\label{sec-new}


Before we start with the proof, we recall that a standard application of the large sieve (cf.\ e.g.\ \cite[Theorm 7.35]{IK}) shows that
\begin{equation}\label{largesieve}
\sum_{f \in \mathcal{B}^{\ast}(q)}  L(1/2, f)^4 |h(t_f)|+ 
    \sum_{f \in   \mathcal{B}_{\text{hol}}^{\ast}(q)}  L(1/2, f)^4  |h^{\text{hol}}(k_f) | \ll q^{1+\varepsilon}
    \end{equation}
    for any $q \in \Bbb{N}$ (not necessarily prime), whenever $\mathfrak{h} = (h, h^{\text{hol}})$ satisfies \eqref{weakly}.   Moreover, by another standard application of the large sieve (\cite[Theorem 7.34]{IK}) we have 
   \begin{equation}\label{largesieve2}
   \sum_{c_{\chi}^2 \mid q} \int_{-\infty}^{\infty} |L(1/2 + it, \chi)|^8  |h(t)| dt \leq   \sum_{c_{\chi}\leq q^{1/2} } \int_{-\infty}^{\infty} |L(1/2 + it, \chi)|^8    |h(t)| dt  \ll q^{1+\varepsilon}.
      \end{equation} 
       We will also frequently use the standard bounds
\begin{equation}\label{hl}
\begin{split}
(q(1 + |t_f|))^{-\varepsilon} & \ll L(1, \text{Ad}^2f) \ll (q(1 + |t_f|))^{\varepsilon}, \quad f \in \mathcal{B}^{\ast}(q),\\
 (q k_f)^{-\varepsilon}  &\ll L(1, \text{Ad}^2f) \ll (q k_f)^{\varepsilon}, \quad\quad\quad\,\,\,\,\, f \in \mathcal{B}_{\text{hol}}^{\ast}(q),\\
 (q(1+|t|))^{-\varepsilon} & \ll |L(1 + 2 it, \chi)|, \quad\quad\quad\quad\quad\quad\quad t \in \Bbb{R}, c_{\chi}^2 \mid q.
\end{split}
\end{equation}
Finally we recall that the central values $L(1/2, f)$ for $f \in \mathcal{B}^{\ast}(q) \cup \mathcal{B}^{\ast}_{\text{hol}}(q)$ are non-negative \cite{KZ, KaS}, for arbitrary $q \in \Bbb{N}$.

Now let $F = \E_0$ be the minimal parabolic Eisenstein series with trivial spectral parameters. By Theorem \ref{thm1}   with $\theta = 0$ we have for $(q, \ell) = 1$, $\mathfrak{h}$ admissible that 
\begin{equation}\label{moeb}
  \sum_{a b =  \ell}  \mathcal{M}^+_{q, a}(1/2, 1/2,  \mathfrak{h})  \left(\frac{a}{b }\right)^{\frac{1}{2}} \ll (\ell q)^{\varepsilon}\sum_{a b =  \ell}  \left(\frac{a}{b}\right)^{\frac{1}{2}} \Bigr(\frac{1}{a} + \frac{1}{q} + \sum_{\pm}\big |\mathcal{M}^{\pm}_{a, q}(1/2, 1/2;   \mathscr{T}_{1/2, 1/2}^{\pm}\mathfrak{h})\big|\Bigr).
\end{equation}
For $\Lambda_f(a; w)$  as in \eqref{Lambda} we have
$$ \sum_{a b = \ell} \frac{\Lambda_f(a; 1/2)}{a^{1/2}}  \left(\frac{a}{b}\right)^{1/2} =   \lambda_f(\ell).$$
Now let $q$ be squarefree. Then we can use  \eqref{largesieve2} with $q = 1$ together with \eqref{local-bound}  and the last bound in \eqref{hl} to estimate the Eisenstein contribution in \eqref{moment}. In this way we see that 
 \begin{displaymath}
 \begin{split}
\sum_{a b =  \ell}  \mathcal{M}^+_{q, a}(1/2, 1/2, \mathfrak{h})& \left(\frac{a}{b }\right)^{1/2}  =   
\frac{\phi(q)}{q^2} \sum_{d_0 \mid q} \sum_{f \in \mathcal{B}^{\ast} (d_0)} \frac{L(1/2, f)^4}{L(1, \text{Ad}^2 f)} \tilde{L}_q(1/2, 1/2, f \times \E_0) \lambda_f(\ell) h(t_f) \\
&  +  \frac{\phi(q)}{q^2}  \sum_{d_0 \mid q} \sum_{f \in  \mathcal{B}^{\ast}_{\text{hol}}(d_0)} \frac{L(1/2, f)^4}{L(1, \text{Ad}^2 f)}  \tilde{L}_q(1/2, 1/2, f \times \E_0)  \lambda_f(\ell) h^{\text{hol}}(k_f) + O(q^{-1} \ell^{\varepsilon}).
  \end{split}
  \end{displaymath}
On the other hand, again by \eqref{moment}, \eqref{moment-together},   \eqref{Lambda} with $q$ in place of $\ell$ and \eqref{local-bound} with $\theta = 0$ and $a$ in place of $q$   
we have 
 \begin{displaymath}
 \begin{split}
 \mathcal{M}^{\pm}_{a, q}(1/2, 1/2; &\mathscr{T}_{1/2, 1/2}^{\pm}\mathfrak{h}) \ll  a^{\varepsilon-1} q^{\vartheta - 1/2+\varepsilon}  \left(\sum_{a_0 \mid a}\sum_{f\in \mathcal{B}^{\ast}(a_0)} \frac{|L(1/2, f)|^4}{L(1, \text{Ad}^2 f)} (1+|t_f|)^{-15} \right.\\
& \left. + \sum_{a_0 \mid a}\sum_{f\in \mathcal{B}^{\ast}_{\text{hol}}(a_0)} \frac{|L(1/2, f)|^4}{L(1, \text{Ad}^2 f)} k_f^{-15} + \sum_{\chi: c_{\chi}^2 \mid a} \int_{\Bbb{R}} \frac{L(1/2 + it, \chi)|^8}{|L(1 + 2it, \chi)|^2}(1+|t|)^{-15} \frac{dt}{2\pi} \right).
 \end{split}
  \end{displaymath}
By 
\eqref{largesieve} -- \eqref{hl} and the non-negativity of $L(1/2, f)$  we conclude
$$ \mathcal{M}^{\pm}_{a, q}(1/2, 1/2; \mathscr{T}_{1/2, 1/2}^{\pm}\mathfrak{h}) \ll q^{\vartheta-1/2+\varepsilon} a^{\varepsilon},$$ 
so that the right hand side of \eqref{moeb} is $\ll (q\ell)^{\varepsilon} ( q\ell^{-1/2} +q^{1/2 + \vartheta}\ell^{1/2})$. 
We have shown
\begin{prop}\label{prop-new} For $q$ squarefree, $\varepsilon > 0$, $(\ell, q) = 1$ and $\mathfrak{h}$ admissible, we have
\begin{equation}\label{squarefree}
\begin{split}
&\sum_{d_0 \mid q} \sum_{f \in \mathcal{B}^{\ast} (d_0)} \frac{L(1/2, f)^4}{L(1, \text{{\rm Ad}}^2 f)} \tilde{L}_q(1/2, 1/2, f \times \E_0) \lambda_f(\ell) h(t_f) \\
&  +   \sum_{d_0 \mid q}\sum_{f \in  \mathcal{B}^{\ast}_{\text{{\rm hol}}}(d_0)} \frac{L(1/2, f)^4}{L(1, \text{{\rm Ad}}^2 f)}  \tilde{L}_q(1/2, 1/2, f \times \E_0)  \lambda_f(\ell) h^{\text{{\rm hol}}}(k_f)  \ll  (q\ell)^{\varepsilon}\left( \frac{q}{\ell^{1/2}} +q^{1/2 + \vartheta}\ell^{1/2}\right).
\end{split}
\end{equation}
\end{prop}
If $q$ is prime,   by Lemma \ref{lem2} we have 
 $\tilde{L}_q(1/2, 1/2, f \times \E_0) = 1$ for 
 $f \in \mathcal{B}^{\ast}(q) \cup \mathcal{B}_{\text{hol}}^{\ast}(q)$ and 
$\tilde{L}_q(1/2, 1/2, f \times \E_0) \ll q^{\varepsilon}$ for $f \in \mathcal{B}(1) \cup \mathcal{B}_{\text{hol}}(1)$.
%
The contribution of the level 1 forms in \eqref{squarefree} is $O((q\ell)^{\varepsilon}   \ell^{\vartheta })$. 
This completes the proof of Theorem \ref{thm2}.  

\section{Proofs of Theorems \ref{cor1} and \ref{cor2}} 

\subsection{Proof of Theorem \ref{cor1}} Let $f \in \mathcal{B}^{\ast}(q)$.   By an  approximate functional equation\footnote{This is the only point where an approximate functional equation is used explicitly although it is implicit in \eqref{largesieve}.}, cf.\ \cite[Theorem 5.3]{IK},   
we have
\begin{equation}\label{approx}
L(1/2, f) =  (1 + \omega_f)\sum_{\ell} \frac{\lambda_f(\ell)}{\ell^{1/2}} W_f\Bigl(\frac{\ell}{q^{1/2}}\Bigr)
\end{equation}
where $\omega_f \in \{\pm 1\}$ is the root number and we can choose
$$W_f(x) =\frac{1}{2\pi i} \int_{(\varepsilon)}   \frac{L_{\infty}(1/2 + s, f)}{L_{\infty}(1/2, f)} \frac{G_f(s)}{G_f(0)} x^{-s} \frac{ds}{s}$$
with 
$L_{\infty}(s, f) = \pi^{-s}\textstyle \Gamma\left(\frac{1}{2}(  s + it_f)\right)\Gamma\left(\frac{1}{2}(  s - it_f)\right)$ 
and 
$$G_f(s) = \prod_{j=0}^{1000} \prod_{\epsilon_1, \epsilon_2 \in \{\pm 1\}} \left( \frac{1}{2} + \epsilon_1 s + i \epsilon_2 t_f + j\right).$$
The choice for the particular weight function $G_f$ will become clear in a moment.  
Since $L(1/2, f) \geq 0$  with equality  if $\omega_f = -1$, we can get rid of the root number in \eqref{approx} and obtain by \eqref{hl} that 
\begin{displaymath}
\begin{split}
 \sum_{f \in \mathcal{B}^{\ast}(q)}  &L(1/2, f)^5 e^{- t_f^2} \ll  \sum_{f \in \mathcal{B}^{\ast}(q)}  L(1/2, f)^4   e^{-   t_f^2} \frac{(q  ( 1+ |t_f|))^{\varepsilon}}{L(1, \text{Ad}^2f)}\cdot  2\frac{\lambda_f(\ell)}{\ell^{1/2}} W_f\Bigl(\frac{\ell}{q^{1/2}}\Bigr)\\
 & = 
   \frac{2q^{\varepsilon}}{2\pi i} \int_{(1)} \pi^{-s} q^{s/2} \sum_{\ell}  \ell^{-1/2-s} \sum_{f\in \mathcal{B}^{\ast}(q)} \frac{L(1/2, f)^4}{L(1, \text{Ad}^2 f)} \lambda_f(\ell) \frac{L_{\infty}(1/2 + s, f)}{L_{\infty}(1/2, f)} \frac{G_f(s)}{G_f(0)} e^{-  t_f^2} ( 1+ |t_f|)^{\varepsilon} \frac{ds}{s}.
 \end{split}
 \end{displaymath}
Shifting the contour to the far right, we can truncate the $\ell$-sum at $q^{1/2+\varepsilon}$ at the cost of an error $O(q^{-10})$. Note that now automatically $(q, \ell) = 1$. Having done this, we shift the contour back to $\Re s = \varepsilon$, and by the exponential  decay of $$\frac{L_{\infty}(1/2 + s, f)}{L_{\infty}(1/2, f)} e^{-t_f^2}$$ along vertical lines (uniformly in $t_f$) we can truncate the contour at $|\Im s| \leq (\log q)^2$ with the same error. Hence
$$ \sum_{f \in \mathcal{B}^{\ast}(q)} L(1/2, f)^5 e^{-   t_f^2} \ll q^{\varepsilon} \max_{|\tau| \leq (\log q)^2} \sum_{\ell \leq q^{1/2 + \varepsilon}} \frac{1}{\ell^{1/2}} \Bigl| \sum_{f\in \mathcal{B}^{\ast}(q)} \frac{L(1/2, f)^4}{L(1, \text{Ad}^2 f)} \lambda_f(\ell)h_{\tau}(t_f) \Bigr| + q^{-10}$$
where
$$h_{\tau}(t_f) = \frac{L_{\infty}(1/2 + \varepsilon + i\tau, f)}{L_{\infty}(1/2, f)} \frac{G_f(\varepsilon + i\tau)}{G_f(0)} e^{-  t_f^2} ( 1+ |t_f|)^{\varepsilon}.$$
By construction, the family 
 $\{\mathfrak{h}_{\tau}   =(h_{\tau}, 0)  : |\tau| \leq (\log q)^2\}$ is uniformly admissible.  Theorem \ref{thm2} now yields the desired bound.


\subsection{Proof of Theorem \ref{cor2}} Let $f_0 \in \mathcal{B}^{\ast}(q) \cup \mathcal{B}^{\ast}_{\text{hol}}(q)$. Let $q^{1/100} < L < q$ be a parameter to be determined later in terms of $q$. In the following all implied constants may depend on the archimedean parameter of $f_0$, and $p$ denotes a prime number. 
For $f \in \mathcal{B}(q) \cup \mathcal{B}_{\text{hol}}(q)$ we choose the following amplifier
$$A_f := \Bigl|\sum_{ \substack{p \leq L\\ p \nmid q} } \lambda_f(p) x(p)\Bigr|^2 +  \Bigl|\sum_{\substack{p \leq L\\ p \nmid q}} \lambda_f(p^2) x(p^2)\Bigr|^2, \quad x(n) = \text{sgn}(\lambda_{f_0}(n)).$$
Then  
\begin{equation}\label{lowerbound}
A_{f_0} = \Bigl(\sum_{\substack{ p\leq L\\ p \nmid q} } |\lambda_{f_0}(p) |\Bigr)^2 +  \Bigl(\sum_{\substack{ p\leq L\\ p \nmid q} } |\lambda_{f_0}(p^2)|  \Bigr)^2 \geq \frac{1}{2} \Bigl(\sum_{\substack{ p\leq L\\ p \nmid q} } |\lambda_{f_0}(p)| + |\lambda_{f_0}(p^2)| \Bigr)^2 \gg \frac{L^2}{\log L}
\end{equation}
by the prime number theorem and the Hecke relation $\lambda_{f_0}(p)^2 = 1 + \lambda_{f_0}(p^2)$. On the other hand,
$$A_f = \sum_{\substack{ p \leq L \\ p \nmid q}} (x(p)^2 + x(p^2)^2)+ \sum_{ \substack{p_1, p_2 \leq L \\ p_1p_2 \nmid q}} \left(x(p_1)x(p_2) + \delta_{p_1 = p_2} x(p_1^2) x(p_2^2)\right) \lambda_f(p_1p_2) + \sum_{ \substack{p_1, p_2 \leq L \\ p_1p_2 \nmid q}}x(p_1^2) x(p_2^2) \lambda_f(p_1^2p_2^2).$$
Now suppose that $q$ is squarefree and $(6, q) = 1$. Then by  \eqref{extension} with $A(1, p) = A(p, 1) = \tau_3(p) = 3$ and $|\lambda_f(p)| \leq p^{\vartheta} + p^{-\vartheta}$ we see that 
$$\tilde{L}_q(1/2, 1/2, f \times \E_0) > 0, \quad f \in \mathcal{B}(q) \cup \mathcal{B}_{\text{{\rm hol}}}(q).$$
We employ now 
the admissble function $\mathfrak{h}_{\text{pos}} = (h_{\text{pos}}, h_{\text{pos}}^{\text{hol}})$ with $a = 1000$, $b = 400$,  the  
  non-negativity   of central values and \eqref{hl} to conclude
 \begin{displaymath}
 \begin{split}
 A_{f_0} L(1/2, f_0)^4  \ll & q^{\varepsilon}   \sum_{f \in \mathcal{B}^{\ast}(q)} \frac{L(1/2,  f)^4  }{L(1, \text{{\rm Ad}}^2 f)}  \tilde{L}_q(1/2, 1/2, f \times \E_0)  A_f  h_{\text{pos}}(t_f)\\
 & +  q^{\varepsilon}  \sum_{f \in \mathcal{B}^{\ast}_{\text{{\rm hol}}}(q)} \frac{L(1/2,    f  )^4 }{L(1, \text{{\rm Ad}}^2 f)} \tilde{L}_q(1/2, 1/2, f \times \E_0)   A_f  h_{\text{pos}}^{\text{hol}}(k_f) .
 \end{split}
 \end{displaymath}
We insert the lower bound  \eqref{lowerbound} on the left hand side and the exact formula for $A_f$ on the right hand side. By Proposition \ref{prop-new} with $\ell \in\{ 1, p_1p_2, p_1^2p_2^2\}$, we obtain
$$L^2 L(1/2, f)^4 \ll q^{\varepsilon}(q L + q^{1/2 + \vartheta} L^4).$$
Choosing $L = q^{(1 - 2\vartheta)/6}$ completes the proof.

\end{document}